%% file: schrod_md_arxiv.tex
\documentclass[a4paper]{amsart}
\usepackage[hmargin=1in,vmargin=1in]{geometry}
\usepackage{amssymb,amsfonts,amsmath,amsthm}
\usepackage{amstext}
\usepackage{graphicx}
\usepackage{subfigure}
\usepackage{algorithm}
\usepackage{algorithmic}
\usepackage{xcolor}

\numberwithin{equation}{section}
\newtheorem{theorem}{Theorem}[section]
\newtheorem{lemma}[theorem]{Lemma}
\newtheorem{remark}[theorem]{Remark}

\newcommand{\rset}{\mathbb{R}}      %
\newcommand{\tpsi}{\psi} %
\newcommand{\hpsi}{\phi}

\def\REF#1{{~\cite{{#1}}}}
\def\REFS#1{{~\cite{{#1}}}}
\input ppmacros_2.tex

\begin{document}

  \title[Accuracy of molecular dynamics for crossings of potential surfaces]%
        { %
        Computational error estimates for Born-Oppenheimer molecular dynamics with nearly crossing potential surfaces}
 \author[C. Bayer]{Christian Bayer}
 \address{Weierstrass Institute for Applied Analysis and Stochastics,
 Mohrenstrasse 39,
10117 Berlin,
Germany}
 \email{christian.bayer@wias-berlin.de}
 
  \author[H. Hoel]{H{\aa}kon Hoel}
  \address{Department of Mathematics,
  University of Oslo,
  Postboks 1053 Blindern, 0316 Oslo, Norway}
  \email{haakonah@math.uio.no}

  \author[A. Kadir]{Ashraful Kadir}
  \address{Department of Mathematics,
     Kungl. Tekniska H\"ogskolan,
    100 44 Stockholm,
    Sweden}
  \email{smakadir@kth.se}
 
  \author[P. Plech\'{a}\v{c}]{Petr Plech\'{a}\v{c}}
  \address{Department of Mathematical Sciences, 
     University of Delaware, 
     Newark, DE 19716, 
     USA}
  \email{plechac@math.udel.edu}
 
  \author[M. Sandberg]{Mattias Sandberg}
  \address{Department of Mathematics,
    Kungl. Tekniska H\"ogskolan,
    100 44 Stockholm,
    Sweden}
  \email{msandb@kth.se}
 
  \author[A. Szepessy]{Anders Szepessy}
  \address{Department of Mathematics,
    Kungl. Tekniska H\"ogskolan,
    100 44 Stockholm,
    Sweden}
  \email{szepessy@kth.se}

\begin{abstract}
The difference of the values of observables for the time-independent Schr\"odinger equation, with matrix valued potentials, and the values of  
observables for {\it ab initio} Born-Oppenheimer molecular dynamics, of the ground state, depends on the probability to be in excited states and the electron/nuclei mass ratio.
The paper first proves an error estimate (depending on the electron/nuclei mass ratio and the probability to be in excited states) for this difference of microcanonical observables, assuming that 
molecular dynamics space-time averages converge, with a rate related to the maximal Lyapunov exponent. The error estimate is uniform in the number of particles and
the analysis does not  assume a uniform lower bound on the spectral gap of the electron operator
and consequently the probability to be in excited states can be large. 
A numerical method to determine the probability to be in excited states is then presented, based on 
Ehrenfest molecular dynamics and stability analysis of a perturbed eigenvalue problem.
\end{abstract}


\maketitle
\tableofcontents

\section{Motivation for error estimates in ab initio molecular dynamics}

Molecular dynamics is a computational method
to study molecular systems in materials science, chemistry and molecular biology.
The simulations are used, for example, in designing and understanding
new materials or for
determining biochemical reactions in drug design \REF{frenkel}.

The wide popularity of molecular dynamics simulations
relies on the fact that in many cases it agrees very well with experiments.
Indeed  when we have experimental data it is often possible to verify the correctness of the method
by comparing with experiments at certain parameter regimes.
However, if we want the simulation to predict something that has no comparing experiment,
we need a mathematical estimate of the accuracy of the computation.
In the case of molecular systems with few particles, such studies are made
by directly solving the Schr\"odinger equation. 
A fundamental and still open question in classical molecular dynamics simulations 
is how to verify the accuracy computationally,
i.e., when the solution of the Schr\"odinger equation
is not a computational alternative.

The aim of this paper 
is to determine quantitative error estimates for molecular dynamics, including the case with %
nearly crossing electron 
potential surfaces that can yield large errors, without solving the Schr\"odinger equation but by combining mathematical stability analysis of eigenvalue problems 
with quantitative numerical Ehrenfest molecular dynamics computations of perturbations. %
Having molecular dynamics error estimates opens, for instance, the possibility 
of systematically evaluating which density functionals or
empirical force fields are good approximations and under what conditions the approximation
properties hold.
Computations with such error estimates could also give improved understanding of
when quantum effects are important and when they are not, in particular 
in cases when the Schr\"odinger equation is too computationally complex to solve.

\smallskip
{\it The first step to check the accuracy} of a molecular dynamics simulation is
to know what to compare with.
Here we compare  the value of any {\it observable} $g:\rset^{3N}\rightarrow  \rset$, of nuclei positions $X$,
for the {\it time-independent Schr\"odinger} eigenvalue equation,
$\HOPER \Phi=E\Phi$,
with the corresponding molecular dynamics observable,
defined by a Hamiltonian $\hat H$, the wave function $\Phi:\rset^{3N}\rightarrow\mathbb C^J$ and the eigenvalue $E\in\rset$.
The approximation error we study is therefore the microcanonical setting
\begin{equation}\label{g_def}
  \frac{\int_{\rset^{3N}} \langle \Phi(X),g(X)\Phi(X)\rangle \,dX}{ \int_{\rset^{3N}} \langle \Phi(X),\Phi(X)\rangle \,dX} - 
  \lim_{T\to\infty}\frac{1}{T}\int_0^T g(X_t) dt\COMMA
\end{equation}
based on a molecular dynamics path $X_t$, with total the energy equal to the Schr\"odinger eigenvalue $E$.
The inner product $\langle\cdot,\cdot\rangle$ is defined in $\mathbb C^J$, where $\Phi\in\mathbb C^J$ corresponds
to $J$ discrete (electron) states. 
The observable can for instance be the local potential energy, used in \REF{erik} to 
determine phase-field partial differential equations from molecular dynamics simulations.
For the sake of simplicity we first explain the basic idea for observables depending only on position variables $X$.
The general case of position and momentum dependent observables is treated in Section \ref{sec_scalar} using the Weyl quantization. 
The time-independent Schr\"odinger equation has a remarkable property of accurately predicting 
experiments in combination with no unknown data, thereby forming the foundation of computational chemistry.
However, the drawback is the high dimensional solution space for nuclei-electron systems with several particles, 
restricting numerical solution to small molecules.
In this paper we study the {\it time-independent} setting of the Schr\"odinger equation as the reference.
The proposed approach has the 
advantage of not requiring any initial data as input, %
on the other hand, an assumption on convergence rates of time averages of molecular dynamics observables is needed.

\smallskip
{\it The second step to check the accuracy} is to derive error estimates.  
We have three types of error: time discretization error, sampling error and modeling error.
The time discretization error comes from approximating the differential equation for molecular dynamics with a
numerical method, based on replacing time derivatives with difference quotients 
and time steps $\Delta t$.
The sampling error is due to truncating the infinite $T$ and using a finite value of $T$ in the integral 
in \eqref{g_def}.
The modeling error (also called coarse-graining error) originates from eliminating the electrons in
the Schr\"odinger nuclei-electron system and replacing the nuclei dynamics with their classical paths;
this approximation error was first analyzed by Born and Oppenheimer in their seminal paper \REF{BO}.

The time discretization and truncation error components 
are in some sense simpler to handle by comparing simulations
with different choices of $\Delta t$ and $T$, although it can, of course, be difficult to know that 
the behavior does not change with even smaller $\Delta t$ and larger $T$ due to metastability, see \REF{tony_book, cances}.
The modeling error is more difficult to check since a direct approach would require to solve 
the Schr\"odinger equation.
Currently the Schr\"odinger partial differential equation can only be solved with few particles, therefore
it is not an option to solve the Schr\"odinger equation in general.
The reason for using  molecular dynamics is precisely to avoid solving  the Schr\"odinger equation.
Consequently, the modeling error requires mathematical error analysis. 
Egorov's theorem, cf. \REFS{robert_book,robert},  provides such error estimates and is used also here.
However, in the literature there seems to be  no error analysis that
is precise, simple and constructive enough so that a molecular dynamics simulation can use it in practice to
assess the modeling error also in the case when the electron operator has eigenvalues (i.e., potential surfaces) which are not well separated
and may cause large modeling error. For instance, crossing eigenvalue surfaces  can form so called conical intersections, which provide the mechanism for many chemical reactions, 
see \REFS{teller,kemist}. 
If the excited electron energy levels %
are well separated from the ground state energy, molecular dynamics based on the  ground state energy  is a good approximation, as first analyzed by Born and Oppenheimer \REF{BO}.
On the other hand, in a quantum system with two electron energy levels  that are not well separated,
the dynamic  transition probability
from one state to the other state for a moving particle can be substantial, as determined by Landau and Zener 
\REF{Z} and generalized to a similar dynamic case in two space dimensions in \REF{lasser_teufel}.
We denote the transition probability to go from the ground state to excited states for a time-dependent problem 
"dynamic transition probability", as in the Landau-Zener model, to distinguish it from the "probability to be in excited states" for the time-independent 
Schr\"odinger eigenvalue problem.

Our alternative
error analysis presented  here relates computable dynamic transition probabilities
to the probability to be in excited states for the time-independent Schr\"odinger equation and
is developed with the aim to 
give a different point of view that could help to construct
algorithms that estimate the modeling error in molecular 
dynamics computations.
Our analysis differs from previous ones by combining analytical estimates with computations and it is based on
three main steps:
\begin{itemize}
\item[-] analyzing the time-independent Schr\"odinger equation as the reference model,  including
excited electron states with near crossing potential surfaces and the accuracy of observables as a function of the probability to be in excited states (in Section \ref{sec_scalar}),
\item[-] studying stochastic molecular dynamics, constrained to the manifold of constant energy, based on 
small stochastic perturbations of the standard  Born-Oppenheimer Hamiltonian dynamics,
\item[-]  using stability analysis of a perturbed eigenvalue problem to estimate the probability to be in excited electron states,
based on perturbations related to (Landau-Zener like) dynamic transition probabilities (in Section \ref{sec:eigen_p}),
\item[-] applying Ehrenfest molecular dynamics to computationally estimate  the %
dynamic transition probabilities (in Section \ref{sec_num}).
\end{itemize}
The estimation method is tested on one and two dimensional problems but to conclude how useful the method is for realistic chemistry problems will require more work.

The next section introduces the Schr\"odinger and molecular dynamics models
and ends with a more detailed outline. 

\subsection{The Schr\"odinger  and molecular dynamics models with an outline of the main results}\label{sec:models}

In deriving the approximation of the solutions to the full Schr\"odinger equation,
the heavy particles are often treated within classical mechanics, i.e., by defining
the evolution of their positions and momenta by equations of motions of classical mechanics.
Therefore we denote by $X_t:[0,\infty) \to \R^{3N}$ 
the time-dependent
function of the nuclei positions, %
with time derivatives denoted by
$$
  \DT{X_t} = \frac{d X_t}{dt}\COMMA\;\; \DDT{X_t} = \frac{d^2 X_t}{dt^2}\PERIOD
$$
We denote the Euclidean scalar product on $\R^{3N}$ by
$$
  \EPROD{X}{Y} = \sum_{n=1}^{N}\sum_{i=1}^3 X^n_i Y^n_i\PERIOD 
$$
Furthermore, we use the notation $\GRADX \psi(X) = (\GRAD_{X^1}\psi(X),\dots,\GRAD_{X^N}\psi(X))$,
and as customary \break $\GRAD_{X^n}\psi = (\partial_{X_1^n}\psi, \partial_{X_2^n}\psi, \partial_{X_3^n}\psi)$.
The notation $\psi(X)=\BIGO(M^{-\alpha})$ is also used for complex valued functions, meaning that
$|\psi(X)|=\BIGO(M^{-\alpha})$ holds uniformly in $X$ for $\psi(X)\in \mathbb C^J$.

The {\it time-independent  Schr\"odinger} equation 
\begin{equation}\label{schrodinger_stat}
  \HOPER(X) \Phi(X) = E\Phi(X)
\end{equation}
models many-body (nuclei-electron) quantum systems and  is
obtained from minimization of the energy in the solution space of wave functions, see 
\REFS{schrodinger,schiff,berezin,tanner,lebris,lebris_hand}.
It is an eigenvalue problem for the energy  $E\in\rset$ of the system, 
described by the Hamiltonian operator $\HOPER(X)$  
\begin{equation}\label{V-definition}
   \HOPER(X)= \VOPER(X) - \frac{1}{2} M^{-1}I\sum_{n=1}^N\Delta_{X^n} \COMMA
\end{equation}
where $I$ denotes the $J\times J$ identity matrix,
and by the wave functions, $\Phi(X)$, belonging
 to a set of permissible electron states which we assume, for simplicity, to be finite, $\Phi(X)\in \mathbb C^J$.
We use a normalized solution $\int_{\rset^{3N}}|\Phi(X)|^2 dX=1$ and
without loss of generality we assume that all nuclei have the same mass $M$.
If this is not the case, we can introduce new coordinates
$M_1^{1/2}\tilde X^k=M^{1/2}_kX^k$, which transform the Hamiltonian to the form we want
${\VOPER(M_1^{1/2}M^{-1/2}\tilde X)}- (2M_1)^{-1}I\sum_{k=1}^N\Delta_{\tilde X^k}$. 
In computational chemistry, the  electron operator Hamiltonian $\bar V$, is independent of $M$ and it is
precisely determined by the sum of the
kinetic energy of electrons and the Coulomb interaction between nuclei and electrons.
The wave function is then in a subspace of $L^2(dXdx)$, with electron coordinates $x$, see \cite{lebris}. 
By introducing a finite dimensional basis $\{\phi_i\}_{i=1}^n$ %
of the electron solution space, minimization of
the energy
$\int_{\rset^{3N}}\int_{\rset^n} \Phi(X,x)^* (\bar V(X) - \frac{1}{2} M^{-1}I\sum_{n=1}^N\Delta_{X^n})\Phi(X,x)dxdX$
under the constraint $\int_{\rset^{3N}}\int_{\rset^n}\Phi(X,x)\Phi(X,x)^*dxdX=1$
leads to the eigenvalue problem \eqref{schrodinger_stat} for the solution $\Phi(X,x)=\sum_{j=1}^J \Phi_j(X)\phi_j(x)$, with the $J\times J$ matrix $V:=S^{-1}\bar{\bar V}$ defined by the matrix components
$S_{ij}:=\int_{\rset^n}\phi_i(x)\phi_j(x)^*dx$ and $\bar{\bar V}_{ij}:= \int_{\rset^n}\bar V(X,x)\phi_i(x)\phi_j(x)^*dx$.
To be able to compute with 
several electrons, i.e. $n/3\gg1$,  the electron solution space is often simplified in the
form of Hartree-Fock or density functional approximations, which also lead to  eigenvalue problems
which now are non linear,
see \cite{lebris}.
We assume  that the electron operator $\VOPER(X)$ is linear and self-adjoint,
 in a finite dimensional complex valued Euclidian space $\mathbb C^{J}$ (for simplicity),
with the usual inner product $\LPROD{\cdot}{\cdot}$, %
and acts as a matrix multiplication. %
An essential feature of  the partial differential equation \VIZ{schrodinger_stat}
is the high computational complexity of finding  the solution in a %
subset of the Sobolev space $H^1(\R^{3N})^J$.
The mass of the nuclei, $M$, are much greater than one (electron mass).  

In contrast to the Schr\"odinger equation, 
a {\it molecular dynamics} model of $N$ nuclei $\tilde X:[0,\infty)\to\rset^{3N}$, 
with a given potential $V_p:\rset^{3N}\to \rset$,
can be computationally studied for large $N$  by solving the ordinary differential equations 
$Md^2\tilde X_\tau/d\tau^2=-\GRADX V_p(\tilde X_\tau)$
in the fast time scale.  We will use the slow time scale $t=M^{-1/2}\tau$
with positions $X_t=\tilde X_\tau$ and scaled momenta $P_t:=\dot X_t$ satisfying 
\begin{equation}\label{md_eq}
 \ddot X_t=- \GRADX V_p(X_t)\PERIOD
\end{equation}
In the slow time scale the nuclei move $\mathcal O(1)$ in unit time, since $\dot P_t=-\GRADX V_p(X_t)$.
This computational and conceptual simplification motivates the study to determine 
the potential and its implied accuracy compared with the
the Schr\"odinger equation,  as started  already
in the 1920's with the Born-Oppenheimer approximation \REF{BO}. 
The purpose of our work is to contribute to the current understanding
of such derivations by showing convergence rates under new assumptions.
The precise aim in this paper is to estimate the error 
\begin{equation}\label{approximation}
\frac{\int_{\rset^{3N}} \langle \Phi(X), g(X)\Phi(X)\rangle \,dX}{\int_{\rset^{3N}} \langle\Phi(X),\Phi(X)\rangle \,dX}
 - \lim_{T\to\infty}\frac{1}{T}\int_0^T g(X_t)\, dt
\end{equation}
for scalar smooth observables $g$ of the time-indepedent 
Schr\"odinger equation \eqref{schrodinger_stat} approximated 
by the corresponding molecular dynamics observable 
$\lim_{T\to\infty}T^{-1}\int_0^T g(X_t) \, dt$,
which is computationally cheaper to evaluate in the case with several nuclei; Section \ref{sec_scalar} includes also the general case with position and momentum dependent observables.

The main step to relate the Schr\"odinger wave function 
and the molecular dynamics solution is the so-called zero-order Born-Oppenheimer approximation, 
where $X_t$ solves  the classical
{\it ab initio} molecular dynamics (\ref{md_eq}) with
the potential $V_p:\rset^{3N}\to \rset$ determined as an eigenvalue of the electron Hamiltonian
$\VOPER(X)$ for a given nuclei position $X$. That is,  $V_p(X)=\lambda_{0}(X)$
and 
\begin{equation}\label{eq:ground}
  \VOPER(X)\Psi_{\BO}(X)=\lambda_{0}(X)\Psi_{\BO}( X)\COMMA
\end{equation}
for a normalized  electron eigenvector $\Psi_{\BO}(X)\in \mathbb C^J$ (here the ground state).
The initial data $(X_0,P_0)$ is chosen to have its energy equal to the eigenvalue $E$, i.e. $|P_0|^2/2+\lambda_0(X_0)$.
The Born-Oppenheimer expansion \REF{BO} is an approximation of the solution of the time-independent
Schr\"odinger equation
which is shown in \REFS{hagedorn_egen,martinez} 
to solve the time-independent Schr\"odinger equation approximately when the electron operator $\VOPER$
have eigenvalues $\lambda_j(X)$ that are well separated, satisfying 
\begin{equation}\label{c_low_def}
\mbox{$\min_X\big(\lambda_1(X)-\lambda_0(X)\big)>c\ $  for a positive constant $c$ independent of $M$.}
\end{equation}
This expansion, analyzed by the methods of multiple scales, pseudo-differential operators
and spectral analysis, e.g.,  in \REFS{hagedorn,hagedorn_egen,martinez,fefferman},
can be used to study the approximation error
\eqref{approximation}. %
The work \REF{hagedorn,hagedorn_egen,martinez,martinez2}  prove  asymptotic expansions for the 
eigenfunction. In the literature one can also find precise statements      
on the error for the setting of the time-dependent Schr\"odinger equation,  
e.g.,  in \REFS{martinez2,robert,spohn_egorov,lasser}, which in some sense is easier
since the stability issue is more subtle in the eigenvalue setting.  
The aim of our work is to present a method to analyze this stability of the eigenvalue problem, using
the dynamic transition probability estimated from Ehrenfest dynamics simulations as quantitative input,
without assuming an $M$-uniform lower bound $c$ on the spectral gap.

We present error estimates for molecular dynamics approximation of time-independent observables in quantum mechanics, valid also 
when the eigenvalues $\lambda_n$ of the electron operator $V$ may nearly intersect and %
thereby increase the probability 
of being in excited states. 
For the Schr\"odinger solution $\Phi$ to \eqref{schrodinger_stat}, with the electron ground state $\Psi_0$ satisfying \eqref{eq:ground},  the probability $p_{ex}$ to be in excited states is
\begin{equation}\label{p_e_def1}
p_{ex}:=\int_{\rset^{3N}} \langle\Phi^\perp(X),\Phi^\perp(X)\rangle dX\, ,
\end{equation}
where $\Phi^\perp(X):=\Phi(X)-\langle \Psi_0(X),\Phi(X)\rangle\Psi_0(X)$.
The derivation of the approximation error is divided into three steps which defines the outline of the paper:
\begin{itemize}
\item[I.] The first step, in Section \ref{sec_scalar},  
    proves that observables based on the time-independent  Sch\"odinger equation
    can be approximated by stochastic molecular dynamics on the constant energy manifold. The convergence rate depends on the nuclei/electron mass ratio $M$ and the probability $p_{ex}$ to be in excited states. The proofs use Egorov's theorem and assume
     that  the expected value of space-time averages of the molecular dynamics observable converge in distributional sense with a rate related to the maximal Lyapunov exponent.
    \item[II.] The second step, in  Section \ref{sec:eigen_p},  presents a numerical method to
    determine the probability to be in excited states from stability analysis of a perturbed eigenvalue problem.  The perturbation is determined from a time-dependent excitation problem related to the Landau-Zener model (with $p_d$ denoting the dynamic transition probability to go to the excited state). Resonances of the eigenvalue problem then determines the larger probability  $p_{ex}$ (to be in excited states for the Schr\"odinger eigenvalue problem) from the perturbations $p_d$ in the dynamic problem.
\item[III.] The final step, in Section \ref{sec_num}, 
on numerical results, illustrate that the dynamic transition probability $p_d$ can be
    determined from Ehrenfest molecular dynamics simulations applied to two simple model problems.
\end{itemize}

Section \ref{sec:weyl} proves a lemma on the regularity of the expected value of the molecular dynamics solution
and establishes in three other lemmas estimates on remainder terms in Moyal expansions used in the theorems.
Section \ref{sec:caustics} formulates WKB solutions in the case with caustics.

The main ingredient in step I is Egorov's theorem
which has been used extensively  to study semiclassical limits both in the time-dependent and time-independent case,
 cf. \REFS{robert_book,robert,teufel_rev}.
We derive two estimates: Theorem \ref{thm:potential} proves an estimate of a weighted difference of the Schr\"odinger and molecular dynamics observables and Theorem \ref{ergod_sats} derives an estimate of the difference of the Schr\"odinger and molecular dynamics observables, without the weight. %
The idea to assume ergodic classical dynamics to prove convergence of observables, in \REFS{schnirelman_rep,schnirelman_book}, initiated the activity on quantum ergodicity \REF{ZZ}
and related results on semiclassical limits for scalar potentials, \REF{martinez_helffer}, and matrix
potentials with distinct eigenvalues \REF{glaser}. The work \REF{teufel_rev} includes a recent review of semiclassical limit results in the case of well
separated electron eigenvalue surfaces. In this case with well separated electron eigenvalues, satisfying \eqref{c_low_def}, 
the probability to be in the excited state is proved to be small $p_{ex}=\mathcal O(M^{-1})$. 
We want to include the case when molecular dynamics simulations become inaccurate  in practice.
For an actual  molecular dynamics simulation the mass $M$ is given and the electron potential $V$ does
not depend in $M$, so a positive spectral gap $c$ does not depend on $M$. However,
if the mass is not large enough  the molecular dynamics approximation may be inaccurate.
Our analysis is asymptotic for large $M$ and to include the situation with a small gap for a certain mass 
we therefore consider a model when the spectral gap 
\begin{equation}\label{gap_bound}
\min_X\big(\lambda_1(X)-\lambda_0(X)\big)=:\delta_M>0,
\end{equation}
can be close to zero, 
i.e. when it is not uniformly bounded from below with respect to $M$,
 and  then $p_{ex}$ can be large as seen e.g. in the Landau-Zener probability \eqref{eq:p_LZ} and Figure 
 \ref{fig:computed-estimated}.
 There are four main new 
 ingredients in our theorems compared to previous studies. It is well known that
 the standard proof of Egorov's theorem, in
 the time-dependent case, combined with an assumption on the molecular dynamics convergence rate towards its ergodic limit yields an error estimate for the molecular dynamics observables as compared to the Schr\"odinger 
observables.  The precise form of the assumption on the convergence to the ergodic limit is important for the conclusion of the approximation error compared to Schr\"odinger observables: (1) for instance, assuming an algebraic convergence rate 
$\mathcal O(T^{-1/2})$, for molecular dynamics time-averages of length $T$, would only give 
logarithmic error estimates $\mathcal O(1/\log^{1/2} M)$ (also in the case with a uniform spectral gap);
(2) it is difficult to verify ergodicity for Hamiltonian dynamics, both theoretically and numerically.
The \underline{first new ingredient} here is to use stochastic dynamics constrained on the constant energy surface instead of Hamiltonian dynamics. Introducing small stochastic perturbations on the constant energy manifold have the advantage
to guarantee ergodicity with respect to the microcanonical ensemble 
and to provide observables with higher regularity.
Our \underline{second ingredient} is to
use a more precise convergence rate assumption based on space-time averages,
which can be tested numerically (see Figure \ref{fig313})
and yields a better error estimate $\mathcal O(M^{-\bar\gamma})$, for a positive $\bar\gamma$, %
The \underline{third ingredient} is to prove uniform convergence rates in the number of particles $N$ by careful use of localized mollifiers. Previous convergence proofs require for instance, for a certain constant $C$, the quantity
$\max_{z\in \rset^{6N}, |\alpha|\le CN} |\partial_z^\alpha g(z)| M^{-1}$ to be small but typically this is a large number for systems with many particles unless the mass ratio $M$ is very large and increasing with the number of particles $N$, see \cite{martinez_helffer} and Theorem 15.4 in \cite{zworski}.
The \underline{fourth ingredient} is that we include non uniformly spectral gap bounds \eqref{gap_bound}:
Theorems \ref{thm:potential} and \ref{ergod_sats} prove the error estimate in step I in terms of the probability $p_{ex}$ to be in excited states, which is not a typical {\it a priori} information.
Consequently the computational approximation of $p_{ex}$ in steps II and III is an important complement.  
Although to evaluate the practical value of the suggested method to estimate the  computational error in molecular dynamics these two steps need to be tested on more realistic chemistry problems.

The alternative to consider molecular dynamics at constant temperature using the canonical ensemble 
 is important in practice and it
also has the advantage  of proved ergodic limits.
The canonical ensemble is in a sense
the average of the microcanonical ensemble over all energies weighted with the Gibbs density.
Therefore the results here are useful for ongoing analysis in the constant temperature setting. 

\section{Molecular dynamics approximation of
Schr\"odinger observables for matrix-valued potentials}\label{sec_scalar} %
For every mass $M\gg 1$, we consider  $L^2(\rset^{3N})$ wave functions $\Phi:\rset^{3N}\rightarrow \rset^J$
and eigenvalues $E\in \rset$, solving the Schr\"odinger
eigenvalue equation
\begin{equation}\label{schrodthm1}
\hat H\Phi=E\Phi,
\end{equation}
where
\begin{equation}\label{schrodthm2}
\hat H:= -\frac{1}{2M} \Delta_X + V(X)\PERIOD
\end{equation}
Here, the symmetric smooth matrix valued potential $V:\rset^{3N}\rightarrow\rset^{J^2}$ has (electron) eigenvalues  
$\{\lambda_n\}_{n=0}^{J-1}$, ordered increasingly, with normalized eigenvectors $\{\Psi_n\}_{n=0}^{J-1}$
satisfying \[V(X)\Psi_n(X)=\lambda_n(X)\Psi_n(X)\, \]
and the spectral gap condition \eqref{gap_bound}. The eigenvalues $\{\lambda_n,\ n\ge 1\}$ may be degenerate but $\lambda_0$ is not. 
Both $\Phi$ and $E$ depend on the mass $M$ and we will only consider those eigenfunctions where the corresponding eigenvalue is bounded, i.e., the $\Phi$ for which $E=\mathcal O(1)$
as $M\rightarrow\infty$. That is, 
we assume that for every $M\gg1$ 
there are solutions $\Phi=\Phi_E\in L^2(\rset^{3N})$ with eigenvalue $E=\mathcal O(1)$ to \eqref{schrodthm1}, as $M\rightarrow\infty$.
On the other hand, the matrix $V$ does  not depend on $M$ explicitly  and consequently
also the eigenvectors and eigenvalues $\Psi_n$ and $\lambda_n$ do not depend on $M$ explicitly.
However to include the case of a small spectral gap, we assume that the positive spectral gap, $\delta$, in \eqref{gap_bound} is not
uniform in $M$ and we study estimates depending on the two parameters $M$ and $\delta$.
To handle ergodicity for the Born-Oppenheimer dynamics, we will below assume that
the smallest eigenvalue $\lambda_0(X)$ is smooth and tends to infinity as $|X|\rightarrow\infty$.
This assumption in fact also implies that the spectrum of $\hat H$ is discrete, see \cite{dallara}.

We consider a given smooth scalar observable $g:\rset^{6N}\rightarrow \rset$  
in the Schwartz space constructed
as 
\begin{equation}\label{g_cond}
\begin{split}
&g=\bar g*\phi_\eta\, ,\\
&\sup_{z\in\rset^{6N}} \sum_{|\beta|\le 4}|\partial^\beta_z \bar g(z)|^2 \mbox{ is uniformly bounded with respect to $N$.}
\end{split}
\end{equation}
 where for some $\eta>4M^{-1/2}$ the mollifier $\phi_\eta:\rset^{6N}\to\rset$ is defined by 
$\phi_\eta(z):=(2\pi\eta)^{-3N}e^{-|z|^2/(2\eta)}$ and $\bar g:\rset^{6N}\to\rset$ is smooth and compactly supported.
We assume that each solution $\Phi=\Phi_E$ to \eqref{schrodthm1} is normalized, i.e. $\int_{\rset^{3N}}\langle\Phi(X),\Phi(X)\rangle dX=1$, and
define for each $\Phi$ the observable
\begin{equation}\label{gs_def}
g_{\SCH}:=\int_{\rset^{3N}} \langle \Phi(X),\hat g\Phi(X)\rangle dX\COMMA
\end{equation}
based on the {\it Weyl quantization}, see \cite{martinez_book},
\[
\hat g\Phi(X)= \OP[g]\Phi(X):=(2\pi M^{-1/2})^{-3N}\int_{\rset^{6N}} e^{iM^{1/2}(X-Y)\cdot P} g\left(\frac{X+Y}{2}, P\right) \Phi(Y) dY dP,
\]
using the inner products
\[
\langle w,v\rangle :=\sum_{j=1}^J w_j^* v_j \quad \mbox{ for } w,v\in \mathbb C^J\COMMA\;\;\;
X\cdot P := \sum_{n=1}^{3N} X^nP^n\quad \mbox{ for } X,P\in \rset^{3N}\PERIOD
\]
We note that in the case $g(X,P)=g(X)$ is a function of $X$ only, then the Weyl quantization
 acts as a multiplication operator, i.e. $\hat g\Phi(X)=g(X)\Phi(X)$.
We also define the {\it Wigner transform}
\begin{equation}\label{wigner}
W(X,P) := (2\pi M^{-1/2})^{-3N} \int_{\rset^{3N}} e^{iM^{1/2} Y\cdot P} \left\langle \Phi\left(X +\frac{Y}{2}\right),
 \Phi\left(X -\frac{Y}{2}\right)\right\rangle dY\, ,
 \end{equation}
 and 
 \[
W_{jk} (X,P) := (2\pi M^{-1/2})^{-3N} \int_{\rset^{3N}} e^{iM^{1/2} Y\cdot P}  \Phi_j^*\left(X +\frac{Y}{2}\right)
 \Phi_k\left(X -\frac{Y}{2}\right) dY\COMMA
 \]
 and note that $\Phi=\Phi_E$ is different for different eigenvalues $E$.
 The change of variables
 \[
 X'=X+ \frac{Y}{2} \COMMA\;\;\;\;
 Y'=X- \frac{Y}{2}
 \]
 shows that for any $g:\rset^{6N}\rightarrow\rset$ and $A:\rset^{6N}\rightarrow \rset^{J^2}$
 \[
 \begin{split}
 \int_{\rset^{3N}} \langle \Phi(X), \hat g \Phi(X)\rangle dX 
 &= \int_{\rset^{6N}} g(X,P) W(X,P) dXdP \, , \\ %
  \int_{\rset^{3N}} \langle \Phi(X), \hat A\Phi(X)\rangle dX 
& = \int_{\rset^{6N}} \sum_{j=1}^{J}\sum_{k=1}^{J}A_{jk}(X,P) W_{jk}(X,P) dXdP\, . \\ %
 \end{split}
 \]
 Our estimates of remainder terms use the corresponding function $W^{(s)}$ defined for $s\in [0,1]$ by
 \begin{equation}\label{w_s_def}
 W^{(s)}(X,P) := (2\pi M^{-1/2})^{-3N} \int_{\rset^{3N}} e^{iM^{1/2} Y\cdot P} \left\langle \Phi\left(X +sY\right),
 \Phi\left(X -(1-s)Y\right)\right\rangle dY\, ,
 \end{equation}
 which is the "Wigner"-function corresponding to the alternative quantization
 \begin{equation}\begin{split}\label{s_kvant}
& (2\pi M^{-1/2})^{-3N}\int_{\rset^{9N}} \big\langle\Phi(X), g\big(X+s(Y-X),P\big) \Phi(Y)\big\rangle e^{iM^{1/2}(X-Y)\cdot P} dY dXdP\\
& =\int_{\rset^{6N}} g(X,P) W^{(s)}(X,P) dXdP\, .
 \end{split}
 \end{equation}

 We will approximate the Schr\"odinger observable  by dynamics related to {\it ab initio}  molecular dynamics 
 \begin{equation}\label{MD_0}
 \begin{split}
 \dot X_t &= P_t\, ,\\
 \dot P_t &= -\nabla_X\lambda_0(X_t)\, , %
 \end{split}
 \end{equation}
 based on the Hamiltonian
 \[
 H_0(X,P):= \frac{|P|^2}{2} + \lambda_0(X)\, .
 \]
 We assume that $\lambda_0$ is smooth 
 and coercive in the sense that 
 \begin{equation}\label{coercive}
 \lim_{|X|\rightarrow+\infty}\lambda_0(X)=+\infty\, .
 \end{equation}
  In order to obtain precise estimates of the approximation error which are uniform in $N$
  we consider a mollified Hamiltonian
  \[
  \begin{split}
  H_\eta(X,P)&:=|P|^2/2 + \lambda_\eta(X)\, ,\\
  \lambda_\eta(X)&:=\lambda_0*\phi_\eta(X)\, ,
  \end{split}
  \]
  where $\phi_\eta(X)=(2\pi \eta)^{-3N/2} e^{-|X|^2/(2\eta)}$ is a standard mollifier on the small scale $\sqrt{\eta}$. To obtain good approximation of observables compared to the Schr\"odinger case we use
  $\eta=4M^{-1/2}$.
  To ensure ergodic molecular dynamics and still approximate the Born-Oppenheimer dynamics
   we will compare the Schr\"odinger observable to the observable for
 stochastic dynamics restricted to the manifold 
 \[
 \Sigma_E:=\{(X,P)\in \rset^{6N}\ |\ H_\eta(X,P)=E\}\, ,
 \]
 for the constant energy $E$ equal to the Schr\"odinger eigenvalue, 
 with a small stochastic perturbation of \eqref{MD_0}.
 Let $Z=(X,P)\in \rset^{6N}$ denote the phase space variable, then
 one example of such stochastic Stratonovich dynamics on $\Sigma_E$ takes the form 
  \begin{equation}\label{sde_proj}
  dZ_t = J\nabla H_\eta(Z_t) dt + \sqrt{2\epsilon\tau} \mathbb P(Z_t) \circ d\tilde W_t\, ,
\end{equation}
where  $\epsilon<M^{-1/2}$ is a small positive parameter, 
$\tilde W$ is the standard Wiener process in $\rset^{6N}$ with independent components,
the skew symmetric matrix $J$ is defined by $J\nabla_Z H_\eta(X,P)=(P,-\nabla_X\lambda_\eta(X))$, the parameter $\tau\sim 1$ is positive,
and  the projection $\mathbb P(Z)$ onto the tangent space of $\Sigma_E$ at  $Z\in\Sigma_E$
 reads
\[
\mathbb P(Z):= {\rm Id} - \hat n(Z)\otimes \hat n(Z)
\]
with the normal $\hat n(Z):=\nabla H_\eta(Z)/|\nabla H_\eta(Z)|$ to $\Sigma_E$. %
We assume that $\nabla H_\eta|_{\Sigma_e}\ne 0$
for all $e\in\rset$. 
Since also $H_\eta$ is assumed to be smooth and coercive,
the set $\Sigma_e$ is a smooth compact codimension one manifold in $\rset^{6N}$ for every $e\in\rset$.
The dynamics \eqref{sde_proj} is the projection of the Ito differential equation
\begin{equation}\label{sde_ito}
dz_t= \big(J\nabla H_\eta(z_t) -\epsilon\nabla H_\eta(z_t)\big) dt + \sqrt{2\epsilon \tau} d\tilde W_t
\end{equation}
to $\Sigma_E$ and $z_t$ has a unique  equilibrium measure, which is the Gibbs measure \[e^{-H_\eta(z)/\tau} dz/\int_{\rset^{6N}} e^{-H_\eta(z)/\tau} dz\, .\] 
The diffusion parameter $\epsilon$ can, for instance, also be a semi-positive definite constant diagonal matrix, so that
\eqref{sde_ito} includes the Langevin equation. %
The work \cite{tony_faou} proves that the projected dynamics \eqref{sde_proj}   also is 
ergodic and the equilibrium measure
is the microcanonical measure %
\[
d\nu(Z):=(\int_{\Sigma_E} \frac{d\Sigma}{|\nabla H_\eta|})^{-1} \frac{d\Sigma(Z)}{|\nabla H_\eta(Z)|}, 
\]
i.e. $\lim_{T\rightarrow\infty} T^{-1} \int_0^T g(Z_t) dt = \int_{\Sigma_E} g(z)d\nu(z)$
where $d\Sigma$ is the surface measure on $\Sigma_E$ induced by the Lebesgue measure in $\rset^{6N}$.
Alternative ergodic dynamics on $\Sigma_E$ sampling the microcanonical measure are presented in \cite{tony_faou}.

The Ito dynamics corresponding to \eqref{sde_proj} reads, see \cite{tony_faou} and \cite{tony_eric_cicco},
\begin{equation}\label{ito_eq}
  dZ_t = \underbrace{\mathbb P(Z_t)\big(J\nabla H_\eta(Z_t)  -\epsilon\nabla H_\eta(Z_t)\big)}_{=J\nabla H_\eta(Z_t)} dt  - \epsilon\tau \kappa(Z_t)\hat n(Z_t)dt
  + \sqrt{2\epsilon\tau} \mathbb P(Z_t) d\tilde W_t\, ,
  \end{equation}
  where $\kappa(z):={\rm div} \,\hat n(z)$ is the mean curvature at $z$ on $\Sigma_E$.
   The Kolmogorov equation for the expected value $u(t,z):=\mathbb E[g(\hat Z_T)\ |\ Z_t=z]$, for $t\le T$, of paths $Z$ satisfying \eqref{sde_proj}  %
  becomes
   \begin{equation}\label{kb}
   \begin{split}
   \partial_t u(t,z) + \underbrace{\big(J\nabla H_\eta(z)-\epsilon\tau \kappa(z)\hat n(z)\big)\cdot \nabla_z u(t,z) + \epsilon\tau\, {\rm trace}\big(
   \mathbb P(z)\nabla^2u(t,z)\big)}_{=Lu(t,z)}&=0\quad t<T,\\
 u(T,\cdot) &= g,\\
 \end{split}   
  \end{equation}
  where $\nabla^2u(t,Z)$ denotes the Hessian with respect to $z$ , so that 
  \[{\rm trace}\big(\mathbb P(z)\nabla^2u(t,z)\big)= \sum_{i,j} \mathbb P_{ij}(z)\partial_{z_i}\partial_{z_j} u(t,z) \mbox{ and }\mathbb P_{ij}=\delta_{ij} - \hat n_i\hat n_j.
  \]
   This Kolmogorov equation has an intrinsic  definition based only on the coordinates on $\Sigma_E$, see \cite{tony_faou}.
   The expected value $u$ is smooth and the compact support of $\bar g$ implies that $u$ decays rapidly as $|Z|\rightarrow\infty$, so that $u$ is in the Schwartz class.
 The work \cite{tony_faou} also proves an exponential convergence rate to equilibrium: let the initial data $g$ be smooth
 in a neighborhood of $\Sigma_E$ then there is a positive constant $\gamma_0\sim \epsilon$ such that
 \begin{equation}\label{kb_rate}
 \|u(t,\cdot)- \int_{\Sigma_E}u(T,Z)d\nu(Z)\|_{L^2(d\nu(\Sigma_E))}=\mathcal O( e^{-\gamma_0(T-t)}), \ t<T\, .
 \end{equation}
 The upper and lower bounds in \cite{herau}  
 on the convergence rate exponent $\gamma_0$ for Langevin dynamics in the full $\rset^{6N}$ phase space prove that
 the rate is bounded from below by a constant times $\epsilon$ and from above by a constant times $\log(1/\epsilon)$ where $\epsilon\ll 1$ is the damping factor, as in \eqref{sde_ito},  for fixed positive temperature $\tau$. 
 Let $S_{ts}(Z) := \big(Z_t\ | \ Z_s=Z\big)$, for $t>s$ be 
 the stochastic flow of the dynamics \eqref{sde_proj}, i.e. the solution $Z_t$ that starts in $Z$ at time $s$,  
 for a realization $\tilde W(\cdot,\omega)$ of the Wiener process.
 Then we have $u(t,Z)=\mathbb E[g\circ S_{Tt}(Z)]$. We write $S_t:=S_{t0}$ so that $u(0,Z)=\mathbb E[g\circ S_T(Z)]$.
 This function satisfies the  exponential growth 
 \begin{equation}\label{hatc_first}
 \sup_{Z\in\rset^{6N}} (\sum_{|\beta|\le n} 
 |\partial_{Z_0}^\beta\mathbb E[g\circ S_t(Z_0)]|^2)^{1/2}
 \le e^{\hat C (t+1)}\delta^{\min(0,-n+1)}\,, \quad n\le 4,
 \end{equation}
 where $\hat C$ is independent of $M$, $\delta$ and $\epsilon$, as
 derived in Lemma \ref{lemma_time} under assumption \eqref{cdeupp} including weak near crossing of eigenvalues.

 Define the molecular dynamics microcanonical observable
 \begin{equation}\label{g_limit}
 g_{\MD}(Z_0):= \lim_{T\rightarrow\infty} \frac{1}{T}\int_{0}^T g\big(S_t(Z_0)\big) dt= \int_{\Sigma_E} g(z) d\nu(z)\, ,
 \end{equation}
 depending on the initial energy $E=H_\eta(Z_0)$.
 Our main assumption on the dynamics
 is that we assume that  for some $\epsilon< M^{-1/2}$ %
 there is a  positive constants $\gamma$, %
 independent of $\epsilon$ and $M$, 
such that 
 \begin{equation}\label{g_bar}
 \begin{split}
 & \int_{\rset^{6N}} \frac{2}{T}\int_{T/2}^{T} \big(\mathbb E[g\circ S_s](Z_0) -g_{\MD}(Z_0)\big) \, ds\ W(Z_0) dZ_0
 \le e^{-\gamma T} + \mathcal O(M^{-1})
 , \mbox{ for } T \le  \log M/(\hat C+\gamma)\, .%
\end{split}
\end{equation}
The constant $\hat C$ depends on the observable $g$ so that
 an observable related to a macroscopic quantity may have better convergence rate in the case of a large system.

From \eqref{kb_rate} we know that there is a constant $\gamma\sim\epsilon$ satisfying
\eqref{g_bar}. Assumption \eqref{g_bar} says that the decay with respect to time $T$ of 
\[
\frac{2}{T}\int_{0}^{T/2}\big(u(t,Z_0)-\int_{\rset^{6N}}u(T,\cdot)d\nu\big) dt\]
integrated over the Wigner distribution is exponential with the decay rate $\gamma\sim 1$ (for some $\epsilon<M^{-1}$ and  $T<\log M/(\hat C+\gamma)$) 
 plus possibly an additional constant small term $M^{-1}$ independent of $t$. 
The ergodic assumption 
\[
\lim_{T\rightarrow\infty} \frac{1}{T}\int_0^T g\circ S_{t0}(Z_0)dt=\int_{\rset^{6N}} g(z)d\nu(z)\, ,
\]
for the Hamiltonian dynamics with $\epsilon=0$
is used in \cite{schnirelman_rep, schnirelman_book}, \cite{martinez_helffer} and \cite{glaser} to prove convergence of Schr\"odinger observables. 
This ergodic assumption is difficult to verify numerically, since
a numerical approximation of the dynamics will always perturb the dynamics and it is theoretically unclear
how this perturbation effects the dynamics over infinite time, see \cite{tupper}.  
To find a method to prove ergodicity for general Hamiltonian systems also remains a challenge, cf. \cite{tupper}.
An advantage with 
assumption \eqref{g_bar} is that it can be tested numerically for some initial points $Z_0$, see Figure \ref{fig313},
since the time discretization error can be made small compared to the bound $\mathcal O(M^{-1})$, by taking
a small time step $\Delta t<M^{-1}$, and the simulation time $T<\log M/\hat C$ is finite. 
Not all dynamics satisfy \eqref{g_bar}.
For instance there exists billiard dynamics 
in   "stadium" domains in $\rset^2$ that is proven to be ergodic for almost all initial data but
is non ergodic  for some data with corresponding concentrated Schr\"odinger eigenfunctions, see \cite{hassell}.
Small perturbations of this  data with non ergodic dynamics
will need very long time to reach asymptotically ergodic behavior and consequently also the stochastic
regularization above is likely to have 
$\gamma$ in \eqref{g_bar} depending on $\epsilon$ so that our assumption would not hold.

Having introduced the necessary notation and terminology, we are now ready to present an estimate for the weighted difference between the Schr\"odinger and stochastic molecular dynamics observables that is expressed in terms of the spectral gap, the mass and the excitation probability. Theorem \ref{ergod_sats}, presented in the end of this section, then estimates the difference $g_S-g_{MD}$ of the observables  without the weight.
 \begin{theorem}\label{thm:potential}  
For every $M\gg 1$ consider the Schr\"odinger eigenvalue problem \eqref{schrodthm1}-\eqref{schrodthm2}
with a solution $\Phi\in L^2(\rset^{3N})$, for a corresponding eigenvalue $E=\mathcal O(1)$, 
 and assume:
 \begin{itemize}
 \item[(i)] the observable $g$ satisfies \eqref{g_cond} 
 and
 the potential $V$ %
 is smooth 
 with $\sup_{X\in\rset^{3N}}\sum_{1\le |\sigma|\le 3} |\partial_X^\sigma V(X)|^2$ bounded, uniformly in $N$, where $|\cdot|$ is the matrix $2$-norm in $\mathbb C^J$,%
 \item[(ii)] the molecular dynamics limit $g_{\MD}$ in \eqref{g_limit} satisfies the weak convergence rate 
  \eqref{g_bar}, %
  \item[(iii)] the electron eigenvector $\Psi_0$ and eigenvalue $\lambda_\eta$ are smooth and have bounded derivatives 
  satisfying, 
  for some $\delta=\delta_M>0$ defined in \eqref{gap_bound} (which may depend on $M$) and any $|\beta|\le 2$,
  and for any $|\alpha|\le 4$, 
  \[
  \begin{split}
  |\partial_X^\beta\Psi_0(X)|&=\mathcal O(\delta^{-|\beta|})\, ,\\
  |\partial_X^{\alpha}\lambda_\eta(X)|&=\mathcal O(\delta^{-|\alpha|+1})\, ,
  \end{split}
  \]
  \item[(iv)] the Hamiltonian $H_\eta$ is smooth, satisfies the coercivity \eqref{coercive}
  and $\nabla H_\eta|_{\Sigma_e}\ne 0$ for all $e\in\rset$,
\item[(v)] the exponential bound \eqref{hatc_first} holds with the constant $\hat C$ uniformly bounded in $N$.
 \end{itemize}
  Let the probability to be in excited states be denoted by 
  $p_{ex}=p_{ex}(\Phi)$, as defined in \eqref{p_e_def1}, and define the Wigner transform $W=W(\Phi)$ by \eqref{wigner}.
 Then, for $M^{-1/2}\delta^{-2}=\mathcal O(1)$,
   the molecular dynamics observable $g_{\MD}$ approximates the  Schr\"odinger observable $g_{\SCH}$, defined  in \eqref{gs_def},   with the weighted error estimate
\begin{equation}\label{md_approx}
\begin{split}
|\int_{\rset^{6N}} \big( g_{\SCH} - g_{\MD}(Z_0)\big) W(Z_0) dZ_0|
&\le  C%
\big( e^{\hat C T} (M^{-1}\delta^{-4} +C_N\epsilon\delta^{-1}+\eta^{}\delta^{-2}+p_{ex}^{1/2})  + e^{-\gamma T}\big) \COMMA \\
\end{split}
\end{equation}
for any $T>0$, and by choosing $T$ such that $e^{-(\hat C+\gamma)T}=(M^{-1/2}\delta^{-2}+p_{ex}^{1/2})$, $\eta=4M^{-1/2}$ and $\epsilon$ sufficiently small
\begin{equation}\label{gammaC}
|\int_{\rset^{6N}} \big( g_{\SCH} - g_{\MD}(Z_0)\big) W(Z_0) dZ_0|
\le C(M^{-1/2}\delta^{-2}+p_{ex}^{1/2})^{\frac{\gamma}{\hat C+\gamma}}\, ,
\end{equation}
for a constant $C$,  independent of $M\gg 1, \ p_{ex}\ll 1, \ \delta$, $\gamma$, $N$  and $\hat C$.
   \end{theorem}

To prove this theorem, we 
modify the proof of Egorov's theorem \cite{robert_book} to include stochastic dynamics and the assumption \eqref{g_bar} on the space-time convergence rate and to replace the uniform spectral gap bound for
matrix valued perturbation potentials $V(X)$, in \cite{glaser, teufel_rev}, with Assumption (iii).
 Assumption (iii) with $\delta>0$ excludes
a conical intersection but avoided crossings are allowed with a spectral gap in \eqref{gap_bound}
that is not uniformly bounded from below
as a function of $M$.
The Landau-Zener probability \eqref{eq:p_LZ} in Section \ref{zener_ehren} (combined with \eqref{p_E_def} and \eqref{prop2})
illustrate for instance that $p_{ex}$ can be close to one for an avoided crossing, %
with a small spectral gap $\min_X(\lambda_1(X)-\lambda_0(X))=\delta=M^{-1/4}$.
That $p_{ex}^{1/2}$  can be of the same order as $M^{-1}\delta^{-4}$
 in the estimate \eqref{md_approx} motivates the study of the %
approximation error as a function of $p_{ex}$  and $M$. 
To have $p_{ex}$ large seems to require smaller gaps $\delta\lesssim M^{1/4}$, see Section \ref{sec:eigen_p}.
It would be desirable to have analytic estimates of $p_{ex}$, since it is determined by $M$ and $V$. 
Such estimates are known, $p_{ex}=o(1)$ as $M\to\infty$, for the case when spectral gap 
$\delta>c>0$ is uniformly bounded from below while in the case with no uniform spectral gap the required smoothness of projection to the ground state
is assumed, see \cite{teufel_book}.  The lack of analytic estimates of $p_{ex}$ is our motivation for
using $p_{ex}$ in addition to $M$ as parameter in the Theorems 
and presenting the computational method to determine $p_{ex}$ in this work. %

We note that stochastic dynamics has a regularizing effect since $\mathbb E[g\circ S_s]$, 
with $\epsilon>0$, is in
general more regular than $g\circ S_s$, with $\epsilon=0$. In fact, we do not explicitly need the stochastic flow $S$:
the proof only uses the deterministic value $u(t,\cdot)=\mathbb E[g\circ S_{Tt}(\cdot)]$ and its flow $u(t,\cdot)=:\bar S_{tT} u(T,\cdot)$, for $t<T$, with the generator $\partial_t u(t,\cdot)=- Lu(t,\cdot)$. The exponential bound in assumption (v) is proved
in Lemma \ref{lemma_time} under an assumption that allows weak near crossing of eigenvalues away from regions where $P$ vanishes.

\begin{proof}

Consider the solution operator $e^{iM^{1/2}t \hat H}$ of the time-dependent Schr\"odinger equation
\[
iM^{-1/2} \partial_t \Psi(t,X) = \hat H \Psi(t,X)
\]
defined
by $\Psi(t,\cdot)= e^{-iM^{1/2}t \hat H} \Psi(0,\cdot)$.
We will compare the observable
\[
\begin{split}
g_{\SCH} &= \int_{\rset^{3N}} \langle e^{-iM^{1/2}t \hat H}\Phi(X), \hat g e^{-iM^{1/2}t \hat H}\Phi(X)\rangle dX\\
&= \int_{\rset^{3N}} \langle \Phi(X), e^{iM^{1/2}t \hat H}\hat g e^{-iM^{1/2}t \hat H}\Phi(X)\rangle dX
\end{split}
\]
 to the expected value
\[
\frac{2}{T} \int_{T/2}^T \int_{\rset^{6N}} \mathbb E[g\circ S_t(Z_0)] W(Z_0) dZ_0 dt
=\frac{2}{T} \int_{T/2}^T \int_{\rset^{3N}} \langle \Phi(X_0), \widehat{\mathbb E[g\circ S_t]}\Phi(X_0)\rangle dX_0 dt\, .
\]
By construction we have that $S_t(Z_0)$ is in the compact manifold $\Sigma_{H_\eta(Z_0)}$.
The coercivity \eqref{coercive} of $\lambda_0$ and the decay of $g$ in phase space
imply therefore that $g\big(S_t(Z_0)\big)$ decays sufficiently fast for large $|Z_0|$.
The integrals above are therefore well defined.
Since $V$ is smooth, elliptic regularity implies that the solution $\Phi$ to the Schr\"odinger eigenvalue problem also is smooth. Let $\bar u$ be the solution to 
\begin{equation}\label{kbbar}
\begin{split}
\partial_t \bar u(t,z)-L\bar u(t,z)&=0\, \quad t>0\\
\bar u(0)&=\bar g\, ,
\end{split}
\end{equation}
then  $\mathbb E[\bar g\circ S_t(Z_0)]=\bar u(t,Z_0)$ and $g=\bar g*\phi_\eta$.
We have, using the commutator $[\hat A,\hat B]:=\hat A\hat B-\hat B\hat A$ and that $\hat H\Phi=E\Phi$,
\begin{equation}\label{gg_estimate}
\begin{split}
&\int_{\rset^{6N}} \bar u(t,\cdot)*\phi_\eta(Z) W(Z) dZ - g_{\SCH}\\
&=\int_{\rset^{3N}} \langle \Phi(X), \widehat{\bar u(t,\cdot)*\phi_\eta}- e^{iM^{1/2}t \hat H}\widehat{\bar u(0,\cdot)*\phi_\eta} 
e^{-iM^{1/2}t\hat H}\big]
\Phi(X) \rangle dX\\
&= 
\int_{\rset^{3N}} \langle \Phi(X), \int_0^t
\frac{\partial}{\partial s} \big( e^{iM^{1/2}(t-s)\hat H} \widehat{\bar u(s,\cdot)*\phi_\eta}
e^{-iM^{1/2}(t-s)\hat H}\big)ds
\Phi(X) \rangle dX\\
&=\int_{\rset^{3N}} \langle \Phi(X), \int_0^t
 e^{iM^{1/2}(t-s)\hat H} \Big( \frac{\partial}{\partial s}\OP\big[\bar u(s,\cdot)*\phi_\eta \big] 
 - iM^{1/2}  \big[\hat H, \widehat{\bar u(s,\cdot)*\phi_\eta}\big]\Big)
 e^{-iM^{1/2}(t-s)\hat H}
\Phi(X) \rangle ds dX\\
&= \int_0^t \int_{\rset^{3N}} \langle  e^{iM^{1/2}(t-s) \hat H}\Phi(X),
\OP\Big[\Big(\big(P\cdot\nabla_{X}- \nabla\lambda_\eta\cdot\nabla_{P}\big)\bar u(s,\cdot)\Big)*\phi_\eta(X,P)
\Big] - iM^{1/2} \big[\hat H, \widehat{\bar u(s,\cdot)*\phi_\eta}\big]\\
&\qquad+\OP\Big[\Big(-\frac{\epsilon\tau\kappa}{|\nabla_ZH_\eta|}\nabla_ZH_\eta
\cdot \nabla_Z\bar u(s,\cdot)*\phi_\eta
+\epsilon\tau {\rm trace} \big(\mathbb P\nabla_Z^2\big) \bar u(s,\cdot)\Big)*\phi_\eta(Z)\Big]  \\
&\qquad \times e^{-iM^{1/2}(t-s)\hat H}
\Phi(X) \rangle dX ds\\
&= \int_0^t \int_{\rset^{3N}} \langle  e^{iM^{1/2}(t-s) E}\Phi(X),
 \OP\Big[\Big(\big(P\cdot\nabla_{X}-\nabla\lambda_\eta(X)\cdot\nabla_{P}\big)\bar u(s,\cdot)\Big)*\phi_\eta(Z)\Big] 
 \\
 &\qquad +\mathcal O(\epsilon)   - iM^{1/2} \big[\hat H, \widehat{\bar u(s,\cdot)*\phi_\eta} \big]\Big]
e^{-iM^{1/2}(t-s)E}
\Phi(X) \rangle dX ds\\
&= \int_0^t \int_{\rset^{3N}} \big\langle  \Phi(X),
 \Big( \OP\big[\big(\big(P\cdot\nabla_{X}-\nabla\lambda_\eta(X)\cdot\nabla_{P}\big)
 \bar u(s,\cdot)\big)*\phi_\eta(Z)\big]
  \\
&\qquad 
 + \mathcal O(\epsilon)
 - iM^{1/2} \big[\hat H, \OP [\bar u(s,\cdot)*\phi_\eta ]\big]\Big)  
 \Phi(X) \big\rangle dX ds\, .\\
\end{split}
\end{equation}
In the fourth equality we have used \eqref{kbbar}
to conclude  that
\[
\begin{split}
\frac{\partial}{\partial s} \OP\big[\bar u(s,\cdot)*\phi_\eta(Z)\big]
&=\OP[\frac{\partial}{\partial s}\bar u(s,\cdot)*\phi_\eta(Z)]\\
&=\OP[(L\bar u(s,\cdot))*\phi_\eta(Z)]\\
&=\OP\Big[\Big(\big(J\nabla_{Z}H_\eta-\frac{\epsilon\tau\kappa}{|\nabla_ZH_\eta|}
\nabla_{Z}H_\eta\big)\cdot \nabla_{Z} 
+\epsilon\tau\, {\rm trace} \big(\mathbb P\nabla_{Z}^2 \big)\Big)
\bar u(s,\cdot)*\phi_\eta \Big]\\
&=\OP\Big[\big(J\nabla_{Z}H_\eta\cdot \nabla_{Z}\bar u(s,\cdot)\big)*\phi_\eta(Z)\Big] + 
\mathcal O(C_{N}\epsilon \delta^{-1} e^{\hat Ct})
\, .\\
\end{split}
\]
Here the estimate
\begin{equation}\label{epsL2}
|\int_{\rset^{3N}}\Big\langle\Phi,\OP\Big[\big(-\frac{\epsilon\tau\kappa}{|\nabla_ZH_\eta|}
\nabla_{Z}H_\eta\cdot \nabla_{Z} 
+\epsilon\tau\, {\rm trace} \big(P\nabla_{Z}^2 \big)
\bar u(s,\cdot)*\phi_\eta \Big]\Phi\Big\rangle dX|
\le C_{N}\epsilon \delta^{-1}e^{\hat Ct}
\end{equation} 
 follows by 
 Lemma \ref{lem_cn} and we note that the term is negligible  small by choosing $\epsilon$ sufficiently small,
 although the constant $C_N$ typically is large as it depends on   order $N$ derivatives of the symbol.

To estimate terms in the error representation \eqref{gg_estimate}  is now the remaining four steps in the proof.

{\bf 1.} The Moyal expansion for the commutator of two Weyl operators in \REFS{robert,martinez} and Lemma \ref{moyal_lemma}  is stated for scalars and %
here we apply Lemma \ref{moyal_lemma} to  %
each matrix component $jk$ of $\hat H$ to obtain
\begin{equation}\label{commut}
\begin{split}
&iM^{1/2} [\hat H, \widehat{\bar u(s,\cdot)*\phi_\eta}]_{jk} \Phi_k =
iM^{1/2} [\hat H_{jk}, \widehat{\bar u(s,\cdot)*\phi_\eta}] \Phi_k \\
&=\Big( \OP[\big( \nabla_{P}H_{jk}(X,P)\cdot \nabla_{X} -  \nabla_{X}H_{jk}(X,P) \cdot\nabla_{P}\big)
\big(\bar u(s,\cdot)*\phi_\eta\big)(X,P)] + \hat R_M\Big)\Phi_k(X)
\end{split}
\end{equation}
for the matrix components
\[
H_{jk}(X,P)= \frac{|P|^2}{2}\delta_{jk} + V_{jk}(X)\, .
\]
where the remainder takes the form, see \cite{robert} and Section \ref{sec:weyl},
\begin{equation}\label{R_exp}
R_M:=
\OP[\sum_{n=1}^m
2M^{-n} (2i)^{-n} \sum_{|\alpha|=2n+1} \frac{(-1)^{|\alpha|}}{\alpha!} \partial_{X}^{\alpha} V_{jk}(X)
\big(\partial_{P}^{\alpha}\bar u(s,\cdot)*\phi_\eta\big) (X,P) + M^{-(m+1)} r_m]\, 
\end{equation}
and $r_m$ is smooth. %
We will use %
$\hat R_M$ for $m=1$. 
Lemma \ref{moyal_lemma} in Section \ref{sec:weyl} shows that
\begin{equation}\label{shur}
\int_{\rset^{3N}}\big\langle\Phi,\hat R_M\Phi\big\rangle=\mathcal O( e^{\hat CT}\delta^{-2} M^{-1}).
\end{equation}

{\bf 2.} The normalization property $1=\int_{\rset^{6N}}W(Z)dZ$ implies that
\[
g_{\SCH}= \int_{\rset^{6N}} g_{\SCH} W(Z)dZ
\]
and we obtain by \eqref{gg_estimate} and \eqref{commut} %
\begin{equation}\label{r_c_estimate}
\begin{split}
&\int_{\rset^{6N}} \big(g_{\SCH} - \frac{2}{T}\int_{T/2}^T\bar u(t,\cdot)*\phi_\eta (Z) dt\big) W(Z)dZ\\
&=\frac{2}{T}\int_{T/2}^T \big( g_{\SCH} - \int_{\rset^{6N}}\bar u(t,\cdot)*\phi_\eta (Z)  \big) W(Z)dZ \big) dt\\
&=\frac{2}{T}\int_{T/2}^T \int_0^t\int_{\rset^{3N}}
\langle \Phi(X), \OP[
\big((\nabla\lambda_\eta(\cdot) I -\nabla_{X} V(X))\cdot\nabla_{P}\bar u(s,\cdot)\big)*\phi_\eta(X,P) ]
\Phi(X)\rangle dX ds dt \\
& \qquad + \underbrace{\frac{2}{T}\int_{T/2}^T\int_0^t \int_{\rset^{3N}}\langle \Phi(X), 
(\mathcal O(\epsilon)+\hat R_M)\Phi(X)\rangle dX ds dt}
_{=\mathcal O(e^{\hat CT} (M^{-1}\delta^{-2}+\epsilon C_N\delta^{-1}))} 
\COMMA
\end{split}
\end{equation}
where we use the bounds  \eqref{shur}  and \eqref{epsL2} in last term.

{\bf 3.} We have the desired quantity $\int_{\rset^{6N}}(g_{\SCH}-g_{\MD}) W dXdP$ in the left hand side above
by adding and subtracting the molecular dynamics observable 
$g_{\MD}$, to the second term in the left hand side, 
and using its convergence rate $e^{-\gamma T}+ M^{-1}$ from \eqref{g_bar} as follows
\[
\begin{split}
\int_{\rset^{6N}} \frac{2}{T}\int_{T/2}^T \bar u(t,\cdot)*\phi_\eta (Z)  dt W(Z)dZ
&=\int_{\rset^{6N}} g_{\MD}(Z) W(Z) dZ\\
&\qquad +\int_{\rset^{6N}} \frac{2}{T}\int_{T/2}^T \big(\bar u(t,\cdot)*\phi_\eta (Z) -g_{\MD}(Z)\big) dt\,  W(Z)dZ\\
&=\int_{\rset^{6N}} g_{\MD}(Z) W(Z) dZ + \mathcal O(e^{-\gamma T}+M^{-1}
)\, ,
\end{split}
\]
which yields the last term in the error bound \eqref{md_approx}.

{\bf 4.} It remains to estimate the first term in the right hand side of \eqref{r_c_estimate}. 
The
rule for the composition of Weyl quantizations, cf. \REF{martinez_book} and Lemma \ref{moyal_lemma} in Section \ref{sec:weyl},
 \begin{equation}\label{weyl_extra}
 \hat A \hat B= \widehat{AB} - iM^{-1/2}\OP[{\underbrace{\{A,B\}}_{:=\nabla_{P}A\cdot\nabla_{X}B-\nabla_XA\cdot\nabla_PB}}] + \mathcal O(M^{-1})\, ,
 \end{equation}
 for smooth scalar functions $A,B:\rset^{6N}\rightarrow\rset$. 
The first term in the right hand side of \eqref{r_c_estimate} can be written
\[
\begin{split}
&\int_{\rset^{3N}}\langle \Phi, \OP\big[\Big(\big(\nabla\lambda_\eta(\cdot) I -\nabla_{X} V(X)\big)\cdot\nabla_{P}
\bar u(s,\cdot)\Big)*\phi_\eta (X,P) \big]\Phi\rangle dX\\
\end{split}
\]
 and to estimate this term we use  \eqref{p_e_def1} 
and write  $\Phi= \Phi^\perp+\Pi_0\Phi$, where  
$\|\Phi^\perp\|_{L^2(dX)}=p_{ex}^{1/2}$ and $\Pi_0 w:=\langle w,\Psi_0\rangle\Psi_0$,  for $w\in\rset^J$, defines
the projection onto the ground state.
The  integral terms including the factor $\Phi^\perp$  is bounded using  Lemma \ref{moyal_lemma} 
\[
\begin{split}
&\int_{\rset^{3N}}\langle \Phi^\perp, \OP\big[ -\big(\nabla\lambda_\eta(\cdot) I -\nabla_{X} V(X)\big)\cdot \nabla_{P}
\bar u(s,\cdot)*\phi_\eta (X,P) \big]\Phi\rangle dX\\
&\qquad +\int_{\rset^{3N}}\langle \Phi, \OP\big[ -\big(\nabla\lambda_\eta(\cdot) I -\nabla_{X} V(X)\big)\cdot \nabla_{P}\bar u(s,\cdot)*\phi_\eta(X,P) \big]\Phi^\perp\rangle dX=\mathcal O(e^{\hat CT}p_{ex}^{1/2})\, .
\end{split}
\]
The main term related to the consistency of the approximation can be written as follows
\[
\begin{split}
&\OP\big[ -\big(\nabla\lambda_\eta(\cdot) I -\nabla_{X} V(X)\big)\cdot \nabla_{P}
\bar u(s,\cdot)*\phi_\eta (X,P) \big]\\
&=
\OP\big[ -\big(\nabla\lambda_\eta(X) I -\nabla_{X} V(X)\big)\cdot \nabla_{P}
\bar u(s,\cdot)*\phi_\eta (X,P) \big]\\
&\qquad + 
\OP\big[ -\big(\nabla\lambda_\eta(\cdot)  -\nabla \lambda_\eta(X)\big)\cdot \nabla_{P}
\bar u(s,\cdot)*\phi_\eta (X,P) \big]\, ,
\end{split}
\]
where Lemma \ref{u_lemma} shows that the last term is bounded by 
\[
\sup_{(X,P)\in\rset^{6N}} \underbrace{|\big(  \nabla \lambda_\eta(X)-\nabla\lambda_\eta(\cdot) \big)\cdot \nabla_{P}
\bar u(s,\cdot)*\phi_\eta (X,P) |}_{=\mathcal O(e^{\hat C s}\eta\delta^{-2})} + \mathcal O(\eta^{}e^{\hat C s}\delta^{-2})=
\mathcal O(e^{\hat C s}\eta^{}\delta^{-2})
\]
since the symbol satisfies
\[
\begin{split}
&\int_{\rset^{6N}}\big(\nabla\lambda_\eta(X)  -\nabla \lambda_\eta(X-X')\big)
\cdot \nabla_{P}
\bar u(s,X-X',P-P') \phi_\eta (X',P')dX'dP'\\
&=\int_{\rset^{3N}}\sum_{j=1}^{3N}\partial_{X_j}\nabla\lambda_\eta(X)  
\cdot \nabla_{P}
\bar u(s,X,P-P') \underbrace{\int_{\rset^{3N}}X'_j\phi_\eta (X',P')dX'}_{=0}\, dP' \\
&
\quad + \int_{\rset^{6N}}\sum_{j,k=1}^{3N}\int_0^1(1-t)
\Big(\partial_{X_jX_k}\nabla\lambda_\eta(X-tX') \cdot \nabla_{P} \bar u(s,X-tX',P-P') \\
&\qquad
+2\partial_{X_j}\nabla\lambda_\eta(X-tX') \cdot \partial_{X_k}\nabla_{P} \bar u(s,X-tX',P-P') \\
&\qquad
+\partial_{X_j}\big(\nabla\lambda_\eta(X)  -\nabla \lambda_\eta(X-tX')\big)\cdot \partial_{X_k}\nabla_{P} \bar u(s,X-tX',P-P')\Big) 
X'_jX'_k\phi_\eta (X',P')dX'dP' dt\\
&\le 
\|\sum_{j,k=1}^{3N}
\int_0^1(1-t)
\Big(|\partial_{X_jX_k}\nabla\lambda_\eta(X-tX') \cdot \nabla_{P} \bar u(s,X-tX',P-P')| \\
&\qquad
+2|\partial_{X_j}\nabla\lambda_\eta(X-tX') \cdot \partial_{X_k}\nabla_{P} \bar u(s,X-tX',P-P') |\\
&\qquad
+|\partial_{X_j}\big(\nabla\lambda_\eta(X)  
-\nabla \lambda_\eta(X-tX')\big)\cdot \partial_{X_k}\nabla_{P} \bar u(s,X-tX',P-P')|\Big)dt
\|_{L^\infty} \\
&\qquad \times\sup_{j,k}\|X'_jX'_k\phi_\eta (X',P')\|_{L^1} dt\, .\\
\end{split}
\]

The remaining term satisfies, using \eqref{weyl_extra} and the notation $\check V(X):=V(X)-\lambda_\eta(X)$,
\begin{equation}\label{extra}
\begin{split}
&\int_{\rset^{3N}}\langle \Pi_0\Phi, \OP\big[ \nabla_{X} \check V(X)\cdot \nabla_{P}\bar u(s,\cdot)*\phi_\eta(X,P) \big]\Pi_0\Phi\rangle dX\\
&=\int_{\rset^{3N}}\langle \Phi, \Pi_0\OP\big[ \nabla_{X} \check V(X)\cdot \nabla_{P}\bar u(s,\cdot)*\phi_\eta(X,P) \big]\Pi_0\Phi\rangle dX\\
&=\int_{\rset^{3N}}\langle \Phi, \Pi_0\frac{1}{2}\Big(
\nabla_{X} \check V(X)\cdot \OP\big[ \nabla_{P}\bar u(s,\cdot)*\phi_\eta(X,P) \big]\\
&\qquad+\OP\big[ \nabla_{P}\bar u(s,\cdot)*\phi_\eta(X,P) \big]\cdot \nabla_{X} \check V(X)\Big)
\Pi_0\Phi\rangle dX + \mathcal O(e^{\hat CT}M^{-1}\delta^{-4})\\
&=
\int_{\rset^{3N}}\langle \Phi, \underbrace{\Pi_0\nabla_{X} \check V(X)\Pi_0}_{=\mathcal O(\eta)} \cdot \OP\big[\nabla_{P}\bar u(s,\cdot)*\phi_\eta \big]
\Phi\rangle dX\\
&\qquad +
\int_{\rset^{3N}}\langle \Phi, \Pi_0\nabla_{X} \check V(X)\cdot \big( \OP\big[\nabla_{P}
\bar u(s,\cdot)*\phi_\eta(X,P) \big] \underbrace{(\Pi_0-\Pi_0*\phi_\eta)}_{=\mathcal O(\eta\delta^{-2})}\Phi\rangle dX\\
&\qquad -
\int_{\rset^{3N}}\langle \Phi, \Pi_0\nabla_{X} \check V(X)\cdot 
 \underbrace{(\Pi_0-\Pi_0*\phi_\eta)}_{=\mathcal O(\eta\delta^{-2})}
\OP\big[\nabla_{P}\bar u(s,\cdot)*\phi_\eta(X,P) \big] \big)
\Phi\rangle dX\\
&+
\int_{\rset^{3N}}\langle \Phi, \Pi_0\nabla_{X} \check V(X)\\
&\qquad \cdot\underbrace{\big( \OP\big[\nabla_{P}\bar u(s,\cdot)*\phi_\eta(X,P) \big] \Pi_0*\phi_\eta
- \Pi_0*\phi_\eta \OP\big[\nabla_{P}\bar u(s,\cdot)*\phi_\eta(X,P) \big] \big)}_{=\mathcal O(e^{\hat CT}M^{-1}\delta^{-4})}
\Phi\rangle dX
\end{split}
\end{equation}
where the first term is $\mathcal O(\eta)$
 since $\check V=(V-\lambda_0) + \lambda_0-\lambda_\eta$ and the gradient of the ground state condition $(V-\lambda_0)\Psi_0=0$ implies
\[
\langle \Psi_0,\nabla (V-\lambda_0)\Psi_0\rangle= -\langle \Psi_0,(V-\lambda_0)\nabla\Psi_0\rangle=
 -\langle (V-\lambda_0)\Psi_0,\nabla\Psi_0\rangle=0\COMMA
 \]
 and  the second term is  estimated %
by Lemma \ref{moyal_lemma}.
 Therefore the right hand side in \eqref{r_c_estimate} yields the contribution 
 $\mathcal O(e^{\hat CT}(M^{-1}\delta^{-4} +p_{ex}^{1/2}))$ to the error bound \eqref{md_approx}.
\end{proof}

The next result shows %
that the difference of the observables can be estimated without the weight $W$.
 \begin{theorem}\label{ergod_sats} Let $\bar{E} :=  \int_{\rset^{6N}} H_\eta(Z)W(Z) dZ$ and
 assume that 
 \begin{itemize}
 \item[(i)]the assumptions in Theorem \ref{thm:potential} hold, and
 \item[(ii)]%
 the function $\tilde g_{\MD}:\rset\rightarrow\rset$ is defined by
 \[
 g_{\MD}(Z) =: \tilde g_{\MD}\big( H_\eta(Z)\big)\, ,\mbox{ for   $Z\in \rset^{6N}$,} %
 \]
 \end{itemize}
 then 
 \[
 \int_{\rset^{6N}} g_{\MD}(Z)W(Z) dZ
 =\tilde g_{\MD}(\bar{E}) + \mathcal O( M^{-1}\delta^{-4}+ p_{ex}^{1/2})
 \]
 and
 \[
 g_{\SCH} - \tilde g_{\MD}(\bar{E})=\mathcal O\big((M^{-1/2}\delta^{-2} + p_{ex}^{1/2})^{\frac{\gamma}{\hat C+\gamma}}\big)\, ,
 \]
 where
 \[
 \bar{E}  %
 = E + \mathcal O(p_{ex}+\eta)\, .
 \]
 \end{theorem}

 \begin{proof}
 The function $\tilde g_{MD}(E)=\int_{\Sigma_E} g(z)d\nu(z)/\int_{\Sigma_E} d\nu(z)$
 can be written as the convolution
 \begin{equation}\label{gmd_bar}
 \tilde g_{MD}(E)=\int_{\Sigma_E} g(z)d\nu(z)
 =\int_{\rset^{6N}} \int_{\Sigma_E} \bar g(z-z')d\nu(z)\phi_\eta(z')dz' %
 \end{equation}
 and as
 \[
 \tilde g_{MD}(E)=\frac{\int_{\rset^{6N}} g(Z)\delta\big(H_\eta(Z)-E\big) dZ}{\int_{\rset^{6N}} \delta\big(H_\eta(Z)-E\big) dZ}\, .
 \]
  Using the assumption that $|\nabla_Z H_\eta(Z)|\big|_{\Sigma_E}$ does not vanish, 
  the function $\tilde g_{MD}$  is smooth since we have for $Z=(X,P)$ 
  \[
  \begin{split}
  \frac{d}{dE}\int_{\rset^{6N}} g(Z)\delta\big(H_\eta(Z)-E\big) dZ&=
  -\int_{\rset^{6N}} g(Z)\delta'\big(H_\eta(Z)-E\big) dZ\\
&= -\int_{\rset^{6N}} \frac{g(Z)}{|\nabla H_\eta(Z)|} \frac{\nabla H_\eta(Z)}{|\nabla H_\eta(Z)|}\cdot\nabla H_\eta(Z)\delta'\big(H_\eta(Z)-E\big) dZ\\
&= -\int_{\rset^{6N}} \frac{g(Z)}{|\nabla H_\eta(Z)|} \frac{\nabla H_\eta(Z)}{|\nabla H_\eta(Z)|}\cdot\nabla_Z \delta\big(H_\eta(Z)-E\big) dZ\\
&= \int_{\rset^{6N}}\mbox{div}\big( \frac{g(Z)\nabla H_\eta(Z)}{|\nabla H_\eta(Z)|^2}\big)  \delta\big(H_\eta(Z)-E\big) dZ\\
\end{split}
\]
and similarly for higher order derivatives.

 Taylor expansion implies
 \[
 \tilde g_{\MD}(\tilde E)= \tilde g_{\MD}(\bar E) + \tilde g_{\MD}'(\bar E)(\tilde E-\bar E)
 + \frac{1}{2} \tilde g_{\MD}''(\xi)(\tilde E-\bar E)^2\, ,
 \]
 for some $\xi=\xi(\tilde E)$ between $\tilde E$ and $\bar E$ satisfying 
 $\tilde g_{MD}''\big(\xi(\tilde E)\big):=2\int_0^1 \tilde g_{MD}''\big(s\tilde E+(1-s) E\big)(1-s)ds$,
and  integration  with $\tilde E=H_\eta(Z)$ yields the moment relation
 \begin{equation}\label{g_expand}
 \begin{split}
 & \int_{\rset^{6N}}  g_{\MD}(Z) W(Z) dZ\\
 &=
\int_{\rset^{6N}} \tilde g_{\MD}\big(H_\eta(Z)\big) W(Z) dZ\\
 &=
 \tilde g_{\MD}(\bar E)\underbrace{\int_{\rset^{6N}}W(Z)dZ}_{=1}\\
&\qquad + \tilde g_{\MD}'(\bar E)\underbrace{\int_{\rset^{6N}}\big(H_\eta(Z)-\bar E\big)W(Z)dZ}_{=0}\\
&\qquad +\frac{1}{2}\int_{\rset^{6N}}\tilde g_{\MD}''\Big(\xi\big(H_\eta(Z)\big)\Big)\big(H_\eta(Z)-\bar E\big)^2W(Z)dZ\, .
\end{split}
 \end{equation}
 The definition of $\bar E$ shows the error estimate
 \begin{equation}\label{E_est}
 \begin{split}
 \bar E&= \int_{\rset^{6N}} H_\eta(Z)W(Z) dZ\\
& = \int_{\rset^{3N}} \langle \Phi(X), \hat{H}_\eta\Phi(X)\rangle dX\\
& =\int_{\rset^{3N}} \langle \Phi(X), \hat{ H}\Phi(X)\rangle dX
 -\int_{\rset^{3N}} \langle \Phi(X),( V-\lambda_0)\Phi(X)\rangle dX 
 -\int_{\rset^{3N}} \langle \Phi(X), (\lambda_0-\lambda_\eta)\Phi(X)\rangle dX\\
& =E + \mathcal O( p_{ex})  +\mathcal O(\eta)\, ,
 \end{split}
 \end{equation}
 since $\hat H_\eta= \hat H-\hat{\check V}=\hat H-(V-\lambda_0)-(\lambda_0-\lambda_\eta)$. %
 
 The rule for composition of Weyl quantizations, in \eqref{rem_2_2} in Lemma \ref{moyal_lemma}, shows that
 \[
 \hat A \hat B= \widehat{AB} - iM^{-1/2}\OP[{\underbrace{\{A,B\}}_{=:\nabla_{P}A\cdot\nabla_{X}B-\nabla_XA\cdot\nabla_PB}}] + \mathcal O(M^{-1}\delta^{-4})\, ,
 \]
 for $A=H_\eta-\bar E$ and   $B(Z)=\tilde g_{MD}''(\xi(H_\eta(Z))$,  %
 which implies
 \[
 \OP[\tilde g_{\MD}''(\xi(H_\eta))(H_\eta-\bar E)^2]=
 \OP[H_\eta-\bar E]\OP[\tilde g_{\MD}''(\xi(H_\eta))]\OP[H_\eta-\bar E] +\mathcal O(M^{-1}\delta^{-4})\, ,
 \]
 since
we have 
 \[
 \begin{split}
 \{\tilde g_{\MD}''(\xi(H_\eta)),H_\eta-\bar E\}&=0\\
  \{H_\eta-\bar E,H_\eta-\bar E\}&=0\, .\\
  \end{split}
  \]
 Therefore, %
 the last term in \eqref{g_expand}  becomes
 \[
 \begin{split}
 &\int_{\rset^{6N}}\tilde g_{\MD}''\Big(\xi\big(H_\eta(Z)\big)\Big)\Big(H_\eta(Z)-\bar E\Big)^2W(Z)dZ
 +\mathcal O(M^{-1}\delta^{-4})\\
& =\int_{\rset^{3N}}\langle \widehat{(H_\eta-\bar E)}\Phi(X), \OP[\tilde g_{\MD}''\Big(\xi\big(H_\eta(Z)\big)\Big)]
 \widehat{(H_\eta-\bar E)}\Phi(X)\rangle dX\\
& =\int_{\rset^{3N}}\langle (\hat  H-\bar E -\check V)\Phi(X), \OP[\tilde g_{\MD}''\Big(\xi\big(H_\eta(Z)\big)\Big)]
 (\hat  H-\bar E -\check V)\Phi(X)\rangle dX\\
& =\int_{\rset^{3N}}\langle (E-\bar E -\check V)\Phi(X), \OP[\tilde g_{\MD}''\Big(\xi\big(H_\eta(Z)\big)\Big)]
 (E-\bar E -\check V)\Phi(X)\rangle dX\\
 & =\int_{\rset^{3N}}\langle \underbrace{(E-\bar E)}_{=\mathcal O(p_{ex}+\eta)}
 \Phi -\check V\Phi^\perp(X), \OP[\tilde g_{\MD}''\Big(\xi\big(H_\eta(Z)\big)\Big)]
 (E-\bar E -\check V)\Phi(X)\rangle dX\\
 &=\mathcal O(p_{ex}^{1/2}+\eta)
 \end{split}
 \]
 which combined with \eqref{g_expand} and Theorem \ref{thm:potential} proves the theorem.
 \end{proof}

\section{Determining the probability to be in excited states}\label{sec:eigen_p} %
The purpose of this section is to describe a numerical method to approximate the probability $p_{ex}$ to be in excited states for
the Schr\"odinger eigenvalue problem, without solving a Schr\"odinger eigenvalue problem and instead
 analyze a certain eigenvalue problem with respect to perturbations
related to the dynamic behavior for a time dependent problem.
This perturbation study is formal in the sense that precise
conditions for its validity is not presented here.  The perturbation model is first formulated in one spatial dimension
and then in multiple dimensions in Section \ref{sec_excite}. The perturbation uses the WKB-method to formulate transitions and
is presented in Section~\ref{sec:wkb}. 
The dynamic problem is related to the Landau-Zener model and Ehrenfest dynamics as described in Section~\ref{zener_ehren}.

\subsection{Construction of WKB-solutions}\label{sec:wkb}

To construct a numerical method for approximating the  probability $p_{ex}$ to be in excited states,
without solving a Schr\"odinger equation,  we will use a decomposition of
the Schr\"odinger solution into WKB-functions \eqref{wkb_decomp}. This section presents
a construction of such WKB-functions.

The singular perturbation $-(2M)^{-1}\sum_k\Delta_{X^k}$ of the matrix valued potential $\VOPER$ in the Hamiltonian \eqref{V-definition} introduces
an additional  small scale $M^{-1/2}$ of high frequency oscillations,
as shown by a WKB-expansion, see   \REFS{Rayleigh,jeffreys,sjostrand}.
We shall construct  solutions to (\ref{schrodinger_stat}) in  such a WKB-form 
\begin{equation}\label{wkb_form}
  \Phi(X)=\PSIA(X)\EXP{\Iunit M^{1/2}\theta(X)}\COMMA
\end{equation}
where the amplitude function $\PSIA: \rset^{3N} \to \C^J$ 
is complex valued, the phase $\theta: \rset^{3N}\to \rset$ 
is real valued in the classically allowed region  $\{X\in \rset^{3N}\, \big| \, \lambda_0(X)\le E\}$
(and purely imaginary  elsewhere, see Remark \ref{imaginary}), 
and the factor $M^{1/2}$
is introduced in order to have  well-defined limits of $\PSIA$ and $\theta$  
as  $M\to\infty$.
Note that it is trivially always possible to find functions $\PSIA$ and
$\theta$ satisfying  \VIZ{wkb_form}, even in the sense of a true equality.  Of
course, the ansatz only makes sense if $\PSIA$ and $\theta$ do not have strong
oscillations for large $M$.
The standard WKB-construction  \REFS{maslov,sjostrand} is
based on a series expansion in powers of $M^{1/2}$, which solves
the Schr\"odinger equation with arbitrarily high accuracy. 
Instead of an asymptotic solution,
we introduce  an actual local solution based on
a time-dependent Schr\"odinger transport equation. 
This transport equation
reduces to the formulation
in \REF{maslov}.

\subsubsection{A WKB-solution}\label{first_WKB}
\begin{figure}[htbp]
  \includegraphics[height=9cm]{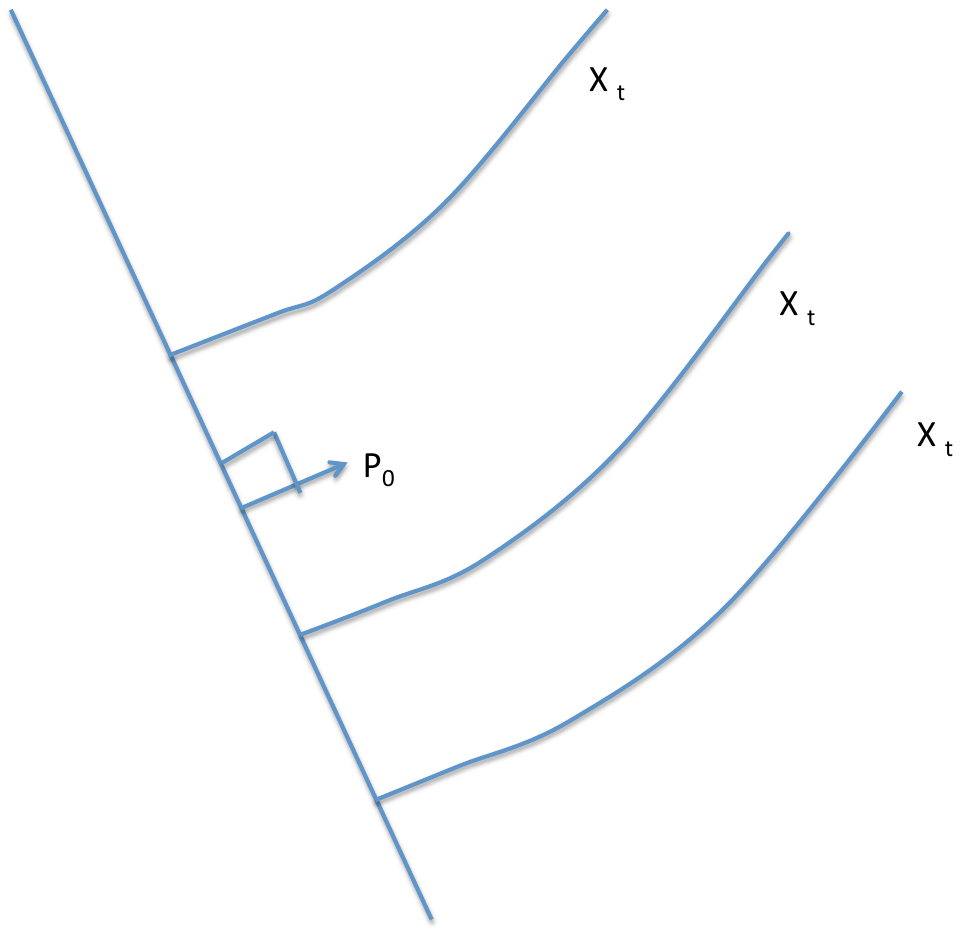}
  \caption{Paths $X_t$ with velocity $\dot X_t=P_t$ going out from the plane with fixed normal $P_0=\mbox{constant}$.}
  \label{char}
 \end{figure}

The WKB-function \eqref{wkb_form} satisfies the Schr\"odinger equation \VIZ{schrodinger_stat} provided that
\begin{equation}\label{wkb_eq}
  \begin{split}
     0 &=(\HOPER-E)\PSIA \, \EXPWKB\\
       &= \left( ( \frac{1}{2}|\GRADX\theta|^2 + \VOPER - E) \PSIA \right.
            - \frac{1}{2M} \Delta_{}\PSIA
              - \frac{\Iunit}{M^{1/2}}   (\EPROD{\GRADX\PSIA}{\GRADX\theta}
       \left. + \frac{1}{2}\PSIA\, \Delta_{}\theta)\right) \EXPWKB \PERIOD %
  \end{split}
\end{equation}
We shall see that only eigensolutions $\Phi$ that correspond to
dynamics without caustics correspond to such a single WKB-mode,
as for instance when the eigenvalue $E$  is inside an electron eigenvalue gap.
Solutions in the presence of caustics use a Fourier integral of such WKB-modes, see
\REF{maslov} and Appendix \ref{sec:caustics}. %
The purpose of the phase function $\theta$ is to generate an accurate
approximation in the limit as $M\to\infty$.
A possible definition, see \REF{maslov},  is 
the {\it Hamilton-Jacobi equation}, also called the {\it eikonal equation}
\begin{equation}\label{theta_eq}
    \frac{1}{2}|\GRADX\theta|^2 =E - \lambda_0 \PERIOD
\end{equation}
The solution to the Hamilton-Jacobi eikonal equation can be constructed locally from the
associated Hamiltonian system
\begin{equation}\label{hj_first}
   \begin{split}
     \DT{X}_t &= P_t\COMMA\\
     \DT{P}_t &= -\GRADX \lambda_0(X_t)\COMMA
     \end{split}
\end{equation}
through the characteristics path $(X_t,P_t)$ satisfying
$\GRADX\theta(X_t)=:P_t$. { For the time being, we fix some $P_0
  \in \rset^{3N}$ and vary $X_0$ over a hyperplane (denoted by $I$) in
  $\rset^{3N}$ \emph{orthogonal to $P_0$}, see Figure \ref{char}.
   Hence, the trajectories $X_t$ can
  (locally) cover the space. Later on, in Section~\ref{sec_excite} we
  shall argue that a \emph{superposition} of such solutions actually allows us
  to reconstruct the true (local) solution of the Schr\"{o}dinger equation. In
  particular, we will solve the system \eqref{hj_first} not just once, but for
  all different $X_0$.}
  
The amplitude function $\phi$ can be determined by requiring the ansatz \eqref{wkb_form} to be a solution to \eqref{wkb_eq}, which gives
\begin{eqnarray*}
     0 &=    &  (\HOPER - E)\PSIA \EXPWKB\\
       &=    & \Big( \underbrace{( \frac{1}{2} |\GRADX\theta|^2 + \lambda_0 - E) }_{=0} \PSIA \\
       &\quad&  %
       { - \frac{1}{2M} \LAP_{}\PSIA +(\VOPER - \lambda_0)\PSIA
                - \frac{\Iunit}{M^{1/2}} ( \EPROD{\GRADX\PSIA}{\GRADX\theta}
                + \frac{1}{2}\PSIA\LAP_{}\theta)\Big)} %
                \EXPWKB\PERIOD
\end{eqnarray*}
Thus by using  \eqref{theta_eq} we have
\begin{equation}\label{transp_2}
   - \frac{1}{2M} \LAP_{}\PSIA +(\VOPER - \lambda_0)\PSIA
   - \frac{\Iunit}{M^{1/2}} ( \EPROD{\GRADX\PSIA}{\GRADX\theta}
   + \frac{1}{2}\PSIA\LAP_{}\theta) =0 \PERIOD
\end{equation}
The usual method for determining $\PSIA$ from this so-called {\it transport equation}
uses an asymptotic expansion $\PSIA\simeq \sum_{k=0}^K M^{-k/2} \PSIA_k$, 
see \REFS{hagedorn_egen,martinez}. %
An alternative is to write it as a Schr\"odinger equation, as in \REF{maslov}:
we apply the characteristics in \eqref{hj_first} to write
\[
  \frac{d}{dt} \PSIA(X_t)= \EPROD{\GRADX\PSIA}{\DT{X}_t}=\EPROD{\GRADX\PSIA}{\GRADX\theta}\COMMA
\]
and define the weight function $G$ by
\begin{equation}\label{G_first}
  \frac{d}{dt} \log G_t= \frac{1}{2} \LAP_{}\theta(X_t)\COMMA
\end{equation}
and the  variable $\tpsi_t:=\PSIA(X_t) G_t$.
Then the transport equation \eqref{transp_2} becomes a Schr\"odinger-type equation
\begin{equation}\label{schrod_first}
  \Iunit M^{-1/2} \DT{\tpsi}_t =(\VOPER - \lambda_0)\tpsi_t 
                               - \frac{G_t}{2M} \LAP_X\left(\frac{\tpsi_t}{G_t}\right) =:\tilde\VOPER\tpsi_t\PERIOD
\end{equation}
{ As~\eqref{schrod_first} is (formally) a time-dependent
  Schr\"{o}dinger equation, we need to impose initial conditions. Once more,
  we refer to Section~\ref{sec_excite} for a detailed discussion of this
  subtle issue.}

{
  The last step in the construction of the WKB system $(\theta, \phi)$ is
  to patch together the different trajectories obtained in
  \eqref{hj_first}--\eqref{schrod_first}. A priori, $X_t$, $P_t$, $G_t$ and
  $\psi_t$ are functions of $(t, X_0)$, for $t \ge 0$ and $X_0 \in I$, where
  $I$ was a hyperplane in $\rset^{3N}$. By uniqueness of solutions
  to~\eqref{hj_first}, the map $(t,X_0) \mapsto X_t$ is locally
  injective. Considering $-P_0$ in addition to $P_0$, we obtain
    that $(t,X_0) \mapsto X_t$ is locally invertible. Hence, we may also
  interpret the functions $P = P(t,X_0)$, $G = G(t,X_0)$, $\psi = \psi(t,X_0)$
  as functions of the space variable alone. Abusing notation, we set
  $\nabla \theta (X) = P(t,X_0)$, $G(X) = G(t,X_0)$, $\psi(X) = \psi(t,X_0)$
  for $X = X(t,X_0)$ in a neighborhood of $I$.
}

In conclusion, equations \eqref{theta_eq}-\eqref{schrod_first} determine the WKB-ansatz
\eqref{wkb_form} to be a local solution to the Schr\"odinger equation \eqref{schrodinger_stat}
in the following sense.  Assume that the Hamilton-Jacobi equation 
$\frac{|\nabla\theta(X)|^2}{2} + \lambda_0(X)=E$
has a $\mathcal C^2$ solution $\theta:\NEIGH\rightarrow \rset$ in a domain $\NEIGH\subseteq\rset^{3N}$.
Let $\dot X_t=\nabla\theta(X_t)$ and $P_t=\nabla\theta(X_t)$, then $(X_t,P_t)$ solves the Hamiltonian system
\eqref{hj_first}, for $t\in [0,t_*]$ such that $X_t\in \NEIGH$.
Then 
\begin{equation}\label{wkb_met}
\Phi(X_t)= G^{-1}(X_t) \psi(X_t) e^{\Iunit M^{1/2}\theta(X_t)}
\end{equation} solves 
\eqref{schrodinger_stat} in $\NEIGH$, since both \eqref{theta_eq} and \eqref{transp_2} are satisfied
and consequently also \eqref{wkb_eq}.
It is well known that Hamilton-Jacobi equations in general do not have global $C^2$ solutions,
due to $X$-paths that collide %
and generate blow up in $\partial_{XX}\theta(X)$.
However if the domain is small enough and the data on the boundary is compatible (in the sense that $H_\eta(X,\GRADX\theta(X))=E$ on the boundary), noncharacteristic (in the sense that the normal derivative
$\partial_n\theta(X)\ne 0$  on the boundary) and $\lambda_0$ is smooth, then the converse property holds,
i.e. the characteristics generate a local solution to the Hamilton-Jacobi equation, see \REF{evans}. 
Maslov's method
to find a global asymptotic solution by patching together local solutions is described in
\REF{maslov}. %

\begin{remark}\label{imaginary}{\rm
In the classically forbidden region $\{X\in\rset^{3N}\, |\, E<\lambda_0(X)\}$
we can apply the same WKB method by replacing the phase function $\theta$ by the imaginary phase $i\theta$.
The Eikonal equation becomes $-|P|^2/2 + \lambda_0(X)=E$, which then has a solution, and
the transport equation becomes real valued %
\[ 
  M^{-1/2} \DT{\tpsi}_t =-(\VOPER - \lambda_0)\tpsi_t 
                               + \frac{G_t}{2M}
                               \LAP_X\left(\frac{\tpsi_t}{G_t}\right)\PERIOD 
\]
The matrix $\VOPER-\lambda_0$ is positive semi-definite, so that $\tpsi$ remains bounded
for  bounded time and approaches the ground state $\Psi_0$.
In the classical region we  instead use the oscillatory behavior to conclude that
$\tpsi$ tends to the ground state, see \REF{old_paper}. }%
\end{remark}

\subsection{The Landau-Zener model, transition probabilities and Ehrenfest
  dynamics}\label{zener_ehren}

{ Notice that the parameter $t$ in~\eqref{hj_first} --
\eqref{schrod_first} is just a numerical parameter, which, prima facie, has
no connection with physical time. Indeed, we are working in a
time-independent setting after all. On the other hand, if $t$ \emph{is}
interpreted as time, \eqref{schrod_first} corresponds to a \emph{time dependent}
Schr\"{o}dinger equation, and we may ask for dynamic transition
probabilities of the time-dependent model, denoted by $p_d(t) = p_d(X_t)$
and rigorously defined in~\eqref{p_E_def}. It turns out below that the
time-dependent transition probabilities are easier to analyze than the time-independent ones, $p_{ex}$, at least
under some simplifying assumptions. The link between the time-dependent
transition probabilities $p_d$ and the time-independent $p_{ex}$ is not
trivial and will be explored in detail in Section~\ref{sec_excite}. To have
some intuition, it might be helpful to think of $p_d(t)$ as a ``local''
excitation probability around $X_t$, with the idea that $p_{ex}$, in turn, is
 given as a time/space average of $p_d$.}

{ We start our discussion by looking at a simple special case
  of~\eqref{schrod_first}, the Landau-Zener model, for which the first results on transition probabilities with crossing or nearly crossing electron potentials
were obtained. It is given by}
\begin{equation}\label{eq:LZ-model}
iM^{-1/2}\dot\phi_t= \left[\begin{array}{cc}
P_0t & \delta\\
\delta & -P_0t\\
\end{array}\right]
\phi_t
\end{equation}
with a wave function $\phi:\rset\rightarrow \mathbb C^2$, constant positive parameters
$(M,P_0,\delta)$  and initial data $\lim_{t\rightarrow-\infty}\phi(t)= (1,0)$.
The transition probability 
\begin{equation}\label{eq:p_LZ}
p_{LZ}:= e^{-\pi\delta^2 M^{1/2}/P_0}\, ,
\end{equation}
is the so-called Landau-Zener probability, determined using Weber functions in
\REF{Z} and illustrated in Figure \ref{fig:md-1d}. { In this particular model, we note that $p_{LZ} =
  \lim_{t\rightarrow\infty}|\phi_2(t)|^2$.} 

{ In the context of \eqref{hj_first}--\eqref{schrod_first}, the
  Landau-Zener model can be seen as a special case of
}
\begin{equation}\label{wkb_ekv1}
\begin{split}
\dot X_t & = P_t\COMMA\\
\dot P_t &= -\nabla\lambda_0(X_t)\COMMA \\
iM^{-1/2}\dot \psi_t &= \underbrace{\big(V(X_t)-\lambda_0(X_t)\big)}_{=:\check V(X_t)}\psi_t -(2M)^{-1}G \Delta(\psi/G)\COMMA
\end{split}
\end{equation}
which by the WKB-method  \eqref{wkb_met}
determines a Schr\"odinger WKB-solution $\Phi_Q$ locally
and hence the transition probability 
\begin{equation}\label{p_E_def}
p_d(Q,X_t):=\langle \psi_t,\psi_t\rangle-|\langle \psi_t,\Psi_0(X_t)\rangle|^2\, ,
\end{equation}
using the initialization $\psi(X_t)=\Psi_0(X_t)$ for  $X_t\in I$, where the inflow domain is given by 
\begin{equation}\label{I_plane}
I=\{X:(X-Y)\cdot P=0\}
\end{equation}
for a given point  $Q:=(Y,P)\in\rset^{6N}$ { and the ground state $\Psi_0$ is defined in \eqref{eq:ground}.}
By neglecting the small term $-G/(2M) \Delta(\psi/G)$ in the transport equation \eqref{wkb_ekv1} and using
the simplification with the potential $\lambda_0(X)=0$ and
\begin{equation}\label{tilde_v_def}
\check V(X)=  \left[\begin{array}{cc}
X & \delta\\
\delta & -X\\
\end{array}\right]\COMMA
\end{equation}
 we obtain the Landau-Zener model. The Landau Zener model was constructed to model and explain the dynamic transitions from the ground state to an excited state when the electron potential surfaces cross
or nearly cross (with a minimal distance $\delta$) and the eigenstates change rapidly (near $X=0$ for small values of $\delta$): the eigenvalues of $\check V$ in \eqref{tilde_v_def} are $\lambda_\pm(X)= \pm \sqrt{X^2+\delta^2}$
and the eigenvectors 
\[
\Psi_\pm = \frac{1}{\sqrt{\delta^2+\big(\lambda_\pm(X)-X\big)^2}}
\left[\begin{array}{c}\delta\\ 
{\lambda_\pm(X)-X}\\
\end{array}\right]\, .
\]

\begin{remark}\label{Dpsi}{\rm
In the one dimensional avoided crossing case \eqref{tilde_v_def}  the electron eigenvectors satisfy
$\Psi_+(X)= f(X/\delta)$ for a smooth function $f:\rset\rightarrow\rset^2$
with
\[
\begin{split}
\lim_{X\rightarrow\infty} \Psi_+(X)&=\left[ \begin{array}{c} 1 \\ 0\end{array}\right]\, ,\\
\lim_{X\rightarrow-\infty} \Psi_+(X)&=\left[ \begin{array}{c} 0 \\ 1\end{array}\right]\, .
\end{split}
\]
We see that 
\[
\begin{split}
\|\partial_X \Psi_+\|_{L^\infty(\rset)} &=\mathcal O(\delta^{-1})\, ,\\
\end{split}
\]
for $\delta\ne 0$. 
The other eigenvector $\Psi_-$ satisfies the same bounds.
Therefore assumption $(iii)$ in Theorem \ref{thm:potential} holds for this avoided crossing with any $\delta\ne 0$.}
\end{remark}

\subsubsection{Ehrenfest dynamics}

As far as computations are concerned, the system~\eqref{wkb_ekv1} is still
very demanding due to the Laplacian on the right hand side. A further
simplification mainly introduced for computational reasons (and used in the
computational examples of this paper) is the so-called
Ehrenfest dynamics
\begin{equation}\label{eq:ehrenfest}
\begin{split}
\dot X_t & = P_t\COMMA\\
\dot P_t &= -\nabla\lambda_0(X_t) -\frac{\langle \psi_t,\nabla\check V(X_t)\psi_t\rangle}{\langle \psi,\psi\rangle}\COMMA\\
iM^{-1/2}\dot \psi_t &= \check V(X_t)\psi_t\COMMA
\end{split}
\end{equation}
which is an approximate WKB solution,
by neglecting the small term $-G/(2M) \Delta(\psi/G)$
in the transport equation and replacing the Eikonal equation with 
$|P|^2/2 + \lambda_0(X) + \langle \psi,\check V\psi\rangle/\langle\psi,\psi\rangle=E$, as in \REF{old_paper}.  If we write the wave function in its real and imaginary parts $\psi=\psi^r+i\psi^i$, and use initial data that satisfy $\langle \psi_0,\psi_0\rangle =2M^{-1/2}$, then the Ehrenfest dynamics \eqref{eq:ehrenfest} is a Hamiltonian system with the Hamiltonian
\begin{equation}\label{eq:H_E}
H_E=|P|^2/2 +\lambda_0(X) +M^{1/2}\langle \psi,\check V(X)\psi\rangle/2\, ,
\end{equation}
using the primal variables
$X$ and $\psi^r$, the dual variables $P$ and $\psi^i$. If the probability, $\langle \psi^\perp,\psi^\perp\rangle$, to be in the excited state is small (using the projection in \eqref{p_e_def1})
the Ehrenfest Hamiltonian is $\mathcal O(\langle \psi^\perp,\psi^\perp\rangle)$ close to the Born-Oppenheimer Hamiltonian
$|P|^2/2 +\lambda_0(X)$.

\subsection{An estimate of the probability to be in excited states for matrix-valued potentials} \label{sec_excite}
This section presents a formal stability study of a
perturbed eigenvalue problem that
provides an approximation for the probability to be in excited states, $p_{ex}$.
The two ingredients are first to determine the perturbation 
as a dynamic transition problem, related to the Landau-Zener model,
and thereafter to  use the stability analysis of a matrix eigenvalue problem
to identify the  probability to be in the excited state with the  squared norm of the change of the eigenvector.
The dynamic transition problem can be approximated numerically using Ehrenfest dynamics,
as described in Sections \ref{zener_ehren} and \ref{sec_num}.
{Our formulation of the  perturbed eigenvalue problem is related to the transformation from local WKB solutions,
which relate to the dynamic transition probability $p_d$, to a global solution of the Schr\"odinger equation.
In dimension one we view a local WKB solution as a solution to the left of the hyperplane $I$ (which is a point in dimension one) and another WKB solution to the right of the hyperplane. The continuity condition to match these two WKB solutions to the right and left
to one global differentiable solution  forms the eigenvalue problem we will study. The eigenvalue problem is
consequently not formulated directly for the unbounded Schr\"odinger operator. It  is instead an eigenvalue problem with bounded solution operators obtained from WKB solutions. A reason we use this form of the perturbation analysis
is that small dynamic transition probabilities $p_d$ from Landau-Zener like dynamics can be viewed
as small regular perturbations while the corresponding perturbations of the potential in the Schr\"odinger equation is not small and not regular. 
 The construction is extended to several dimensions in Section \ref{sec_multiD}.}

\subsubsection{The perturbed eigenvalue problem in one space dimension}\label{sec_3}
Consider first the scalar Schr\"odinger eigenvalue problem in one space dimension.
Assume we have solutions $\Phi$ in the domain $X>0$
and $\Phi$ in $X<0$, with right and left limits {$\lim_{X\rightarrow 0+} \Phi(X)=:\Phi^r, \ \lim_{X\rightarrow 0-} \Phi(X)=:\Phi^\ell,$
and $\lim_{X\rightarrow 0+} \Phi'(X)=:\Phi'^{r}, \ \lim_{X\rightarrow 0-} \Phi'(X)=:\Phi'^{\ell}$.}
In the one dimensional case, the WKB method for an eigenvalue problem typically gives a caustic to the left, with phase $\theta_\ell$,
and a caustic to the right,  with phase $\theta_r$, (see Section \ref{sec:wkb}, \eqref{airy_stat_fas} and \eqref{caustic_expansion}) so that
\begin{equation}\label{wkb_cos}
\begin{split}
\Phi(X)&\simeq \phi_r(X)\cos\big(M^{1/2}\theta_r(X)-\pi/4\big), \mbox{ for $ X>0$,}\\
\Phi(X)&\simeq \phi_\ell(X)\cos\big(M^{1/2}\theta_\ell(X)-\pi/4\big), \mbox{ for $- X>0$,}\\
\end{split}
\end{equation}
for smooth functions $\theta_{r,\ell}$ and $\phi_{r,\ell}$.
The continuity condition in order to have a global solution
\[
\begin{split}
\Phi^r&=\Phi^\ell\, ,\\
\Phi'^r &=\Phi'^\ell\, ,\\
\end{split}
\]
can with the notation $\Phi'^{r}=:R^{r} \Phi^{r}$ (and similarly for the limit to the left)
be written
\begin{equation}\label{egenkont}
\left[
\begin{array}{cc}
1 & -1\\
R^r & -R^\ell\\
\end{array}\right] \left[\begin{array}{c} \Phi^r\\ \Phi^\ell\end{array}\right] =0
\end{equation}
which implies that the $2\times 2$ matrix  
\[
\left[
\begin{array}{cc}
1 & -1\\
R^r& -R^\ell\\
\end{array}
\right]
\]
must be singular, i.e., $R^r=R^\ell$.
For the WKB solutions the derivative satisfies
\[
\begin{split}
R^r&=-M^{1/2}\theta_r'(0)\tan\big(M^{1/2}\theta_r(0)-\pi/4\big) + \mathcal O(M^{0}),\\
R^\ell&=-M^{1/2}\theta_\ell'(0)\tan\big(M^{1/2}\theta_\ell(0)-\pi/4\big) + \mathcal O(M^{0}).\\
\end{split}
\] 
The two eigenvalues
of the matrix 
 \[
 \left[\begin{array}{cc}
1 & -1\\
R^r & -R^\ell\\
\end{array}\right] \, 
\]
satisfy 
 \begin{equation}\label{egenvardet}
\mu= \frac{1-R^\ell}{2} \pm \sqrt{(\frac{1-R^\ell}{2})^2 +R^\ell -R^r}\PERIOD
\end{equation}

\begin{figure}[htbp]
  \includegraphics[height=5cm]{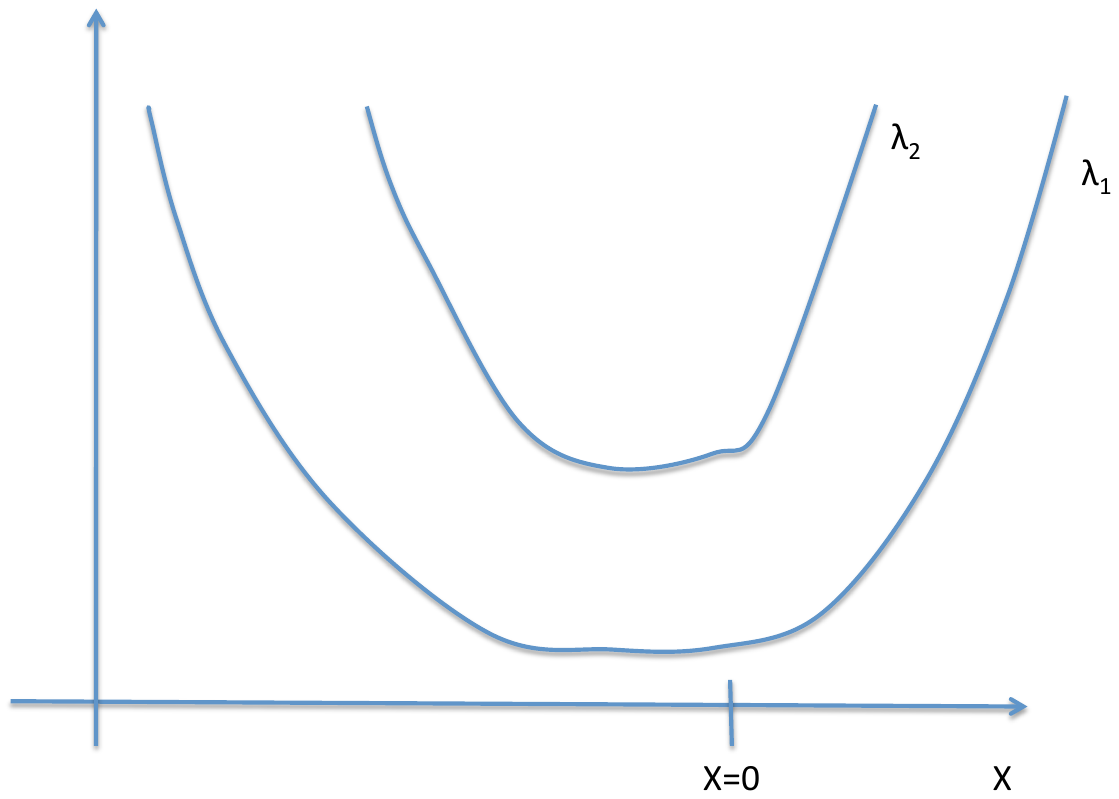}
  \includegraphics[height=5cm]{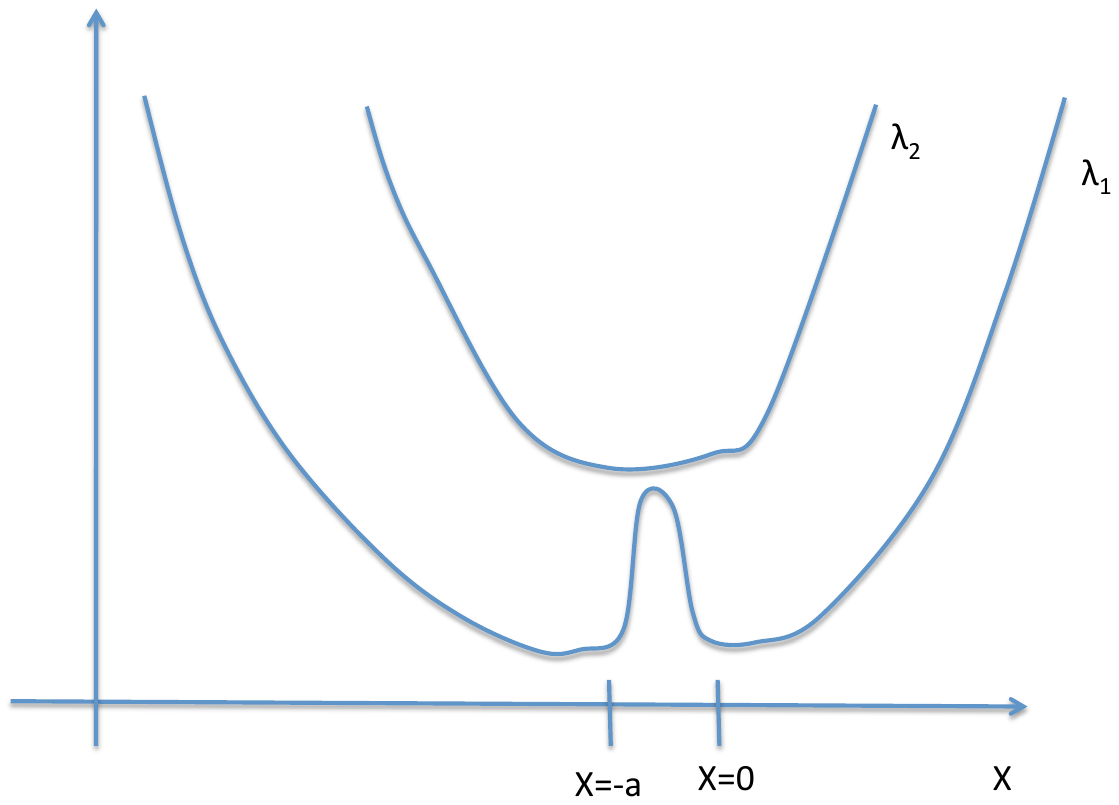}
  \caption{The eigenvalues $\lambda_1$ and $\lambda_2$ of the potential $V$ in the unperturbed (left) and perturbed (right) cases as a function of position $X$. For  $X\notin[-a,0]$ the potential $V(X)$ is diagonal in both the unperturbed and perturbed cases while for $X\in[-a,0]$ the potential $V(X)$ is diagonal  in the unperturbed case and non diagonal in the perturbed case. }
  \label{lamfig}
\end{figure}

Consider the
two decoupled scalar Schr\"odinger eigenvalue problems 
\begin{equation}\label{diagonal_sch}
\begin{split}
\big(-\frac{1}{2M}\Delta + \lambda_1(X)\big)\Phi_1(X) & =E\Phi_1(X)\COMMA \mbox{ and }\\
\big(-\frac{1}{2M}\Delta + \lambda_2(X)\big)\Phi_2(X) & =E\Phi_2(X)\, .\\
\end{split}
\end{equation}
The continuity condition corresponding to \eqref{egenkont}  for \eqref{diagonal_sch}
will be used as our unperturbed eigenvalue problem,
and in one space dimension this continuity  condition becomes the unperturbed problem
\begin{equation}\label{R_first}
\left[
\begin{array}{cc}
\left[\begin{array}{cc}
1 & -1\\
R_1^r & -R^\ell_1\\
\end{array}\right] & 0\\
0& \left[\begin{array}{cc}
1 & -1\\
R_2^r & -R^\ell_2\\
\end{array}\right]\\
\end{array}\right] \left[\begin{array}{c} \Phi_1^r\\ \Phi_1^\ell\\
\Phi_2^r\\ \Phi_2^\ell\end{array}\right] =0\, .
\end{equation}
We consider a special case of the  Schr\"odinger eigenvalue equation \eqref{schrodinger_stat} in one space dimension
where the potential  is diagonal
outside the interval $[-a,0]$, i.e. 
$V(X)=\left[\begin{array}{cc}
\lambda_1(X) & 0 \\
0 & \lambda_2(X)\\
\end{array}\right], \ X\notin [-a,0],$ and view it
 as a perturbation of the diagonal Schr\"odinger equation \eqref{diagonal_sch},
 as illustrated in Figure \ref{lamfig}.
 To have a diagonal potential in the unperturbed case is not necessary. What is important
 is that the electron eigenvalues are well separated so that the transition from ground state to 
 excited states is negligible.

The WKB-method \eqref{wkb_met} shows that 
\begin{equation}\label{bada}
\Phi^\pm(X):=\psi_\pm(X)  G^{-1}_\pm(X)e^{\pm iM^{1/2}\theta(X)}
\end{equation}
are local solutions to the Schr\"odinger equation \eqref{schrodinger_stat}, with positive and negative wave speeds $\pm\theta'(X)$. When the path $X_t$ passes through a domain
where the electron eigenvalues are close, the wave function $\psi_+$ (or $\psi_-$) determined by the transport equation \eqref{schrod_first} yields
the probability $p_d$ to transit from the ground state to an excited state, 
as in the Landau-Zener model \eqref{eq:LZ-model}.
For a WKB-solution \eqref{bada}, with $\Phi=\Phi^+$ (or $\Phi=\Phi^-$), we can define the transition operator $\tilde S^+$ by
\begin{equation}\label{splusminus}
\Phi^+(0)= \tilde S^+\Phi^+(-a)=
\left[\begin{array}{cc} \tilde S^+_{11} &\tilde S^+_{12}\\
\tilde S^+_{21}& \tilde S^+_{22}\end{array}
\right]\left[\begin{array}{c}
\Phi_1^+(-a)\\ \Phi_2^+(-a)\end{array}\right]\, ,
\end{equation}
(and similarly for $\Phi^-$)
{by solving the eikonal and transport equations \eqref{wkb_ekv1} with data given on $X=-a$. The idea is that
given the point $X=-a$ and the momentum $P$ in this point the two first equations determine the path $(X(t),P(t))$.
The third Schr\"odinger like transport equation in \eqref{wkb_ekv1} then yields a linear solution operator along this path determining $\psi$ at $X=0$ from $\psi$ at $X=-a$. Let us now determine $\tilde S$ more precisely.} 
We want to determine the global effect of a small transit probability $|\tilde S_{21}|^2$
viewed as the perturbation to the eigenvalue problem \eqref{R_first} { with  $\tilde S^+\Phi^+(-a)+\tilde S^-\Phi^-(-a)$ replacing $\Phi^\ell$.} To simplify the 
perturbation analysis we  factorize the transition%
\[
\left[\begin{array}{cc} \tilde S_{11} &\tilde S_{12}\\
\tilde S_{21}& \tilde S_{22}\end{array}
\right]
=
\left[\begin{array}{cc} 1 &S_{12}\\
S_{21}& 1\end{array}
\right]
\left[\begin{array}{cc} \tilde S_{11} &0\\
0& \tilde S_{22}\end{array}
\right]
\]
into a transition between the states
$S:=\left[\begin{array}{cc} 1 &S_{12}\\
S_{21}& 1\end{array}
\right]
$
and a diagonal matrix $\left[\begin{array}{cc} \tilde S_{11} &0\\
0& \tilde S_{22}\end{array}
\right]$, which does not contribute to the transition between states.
Since we are only interested in estimates of the small components $\tilde S_{12}=S_{12}\tilde S_{22}$ and $\tilde S_{21}=\tilde S_{11}S_{21}$, we consider only perturbations generated by the
transition matrix $S$. We note that the transition element $\tilde S_{21}$ by \eqref{splusminus} can be written
as the $2$-component of $\Phi(0)$ if $\Phi(-a)=(1,0)$.
Therefore \eqref{bada} shows that the transition element 
\begin{equation}\label{transit_element}
S_{21}^+=\frac{\tilde S^+_{21}}{\tilde S^+_{11}}
=\frac{ G_+(0)e^{iM^{1/2}\theta(0)} \psi_{2+}(0)}{ G_+(0)e^{iM^{1/2}\theta(0)}\psi_{1+}(0)}
=\frac{\psi_{2+}(0)}{\psi_{1+}(0)}\quad \mbox{where $|\psi_{1+}(-a)|=1$}
\end{equation}
can be determined by the WKB-amplitude function $\psi_+$ and similarly we have the transition element $S_{21}^-={\psi_{2-}(0)}/{\psi_{1-}(0)}$. The dynamic transition probability $p_d$, defined by Ehrenfest dynamics  in \eqref{p_E_def}, measures in this case the amplitude squared in the excited state at $X=0$  for a wave starting in $X=-a$ in the ground state. Consequently we have, for $X_t=0$ and $Q$ corresponding to the $+$ wave speed in \eqref{bada}, that 
\[
p_d(Q,0)=\langle \psi_{+}(0),\psi_+(0)\rangle-|\langle \psi_+(0),\Psi_0(0)\rangle|^2
=|\psi_{1+}(0)|^2 + |\psi_{2+}(0)|^2 - |\psi_{1+}(0)|^2=|\psi_{2+}(0)|^2\]
so that by \eqref{transit_element}
 $S^\pm_{21}=\mathcal O(p_d^{1/2})$, as $p_d\rightarrow 0$ (namely as the spectral gap becomes large).
 We assume for simplicity that the transition matrix $S=S^\pm$ is translation invariant, so that
$\Phi_\pm'(0)=S^\pm\Phi_\pm'(-a)$. The total perturbation becomes the sum
$S\Phi:= S^+\Phi^+ + S^-\Phi^-$ of the perturbations of the two WKB-solutions $\Phi^\pm$, for the splitting $\Phi=\Phi^++\Phi^-$.
Hence, we consider perturbations
\[
\Phi^\pm(0)=S^\pm\Phi^\pm(-a)= \left[\begin{array}{cc} 1- \mathcal O(p_{d}) &\mathcal O(p^{1/2}_{d})\\
\mathcal O(p^{1/2}_{d}) &1-\mathcal O(p_{d}) \\
\end{array}
\right]\left[\begin{array}{c}
\Phi^\pm_{1}(-a)\\ \Phi^\pm_{2}(-a)\end{array}\right]\, 
\]
with the transition from the state $1$ to state $2$ determined by  the matrix component $S_{21}^\pm$, which is of the (small) order $p_d^{1/2}$, where
$p_d$, defined in \eqref{p_E_def}, is related to the Landau-Zener probability.
A reason to decompose the solution into WKB solutions is that their transition is determined by the amplitude functions $\psi_\pm$ with approximately conserved norm $|\psi_\pm(t)|$. The approximate conservation follows from 
the conservation $\frac{d}{dt}|\psi_\pm(t)|^2=0$ for the Ehrenfest dynamics \eqref{eq:ehrenfest}
and the fact that
the transport equation \eqref{schrod_first}, determining $\psi_\pm$, becomes the Ehrenfest dynamics in the limit $M\rightarrow\infty$.

The perturbed eigenvalue condition becomes
\begin{equation}\label{eigenpert}
(A+\beta)r_\beta=0\COMMA
\end{equation}
where 
\begin{equation}\label{S_sum}
\begin{split}
 A &:= \left[
\begin{array}{cc}
\left[\begin{array}{cc}
1 & -1\\
R_1^r & -R_1^\ell\\
\end{array}\right] & 0\\
0& \left[\begin{array}{cc}
1 & -1\\
R_2^r & -R_2^\ell\\
\end{array}\right]\\
\end{array}\right]  =: \left[\begin{array}{cc} A_1 & 0\\ 0 &A_2\\ \end{array}\right]\, ,\\
r_\beta &=\left[\begin{array}{c} \Phi_1(0+)\\ \Phi_1(-a)\\
\Phi_2(0+)\\ \Phi_2(-a)\end{array}\right]\, , \\
A+\beta &=\left[\begin{array}{cccc}
1 &- S_{11} & 0 & -S_{12}\\
R_1^r & -S_{11}R_1^\ell & 0 & - S_{12}R_1^\ell\\
0 & -S_{21} & 1 & -S_{22}\\
0 & -S_{21}R_2^\ell & R_2^r & -S_{22}R_2^\ell\\
\end{array}
\right]\, ,\\
S\Phi(-a)&=S^+\Phi^+(-a)+S^-\Phi^-(-a), \quad  \beta \in \mathbb{C}^{4\times4},\\
\end{split}
\end{equation}
and we want to determine  the change in the eigenvector $|r_\beta-r_0|^2$ (which measures the probability to be in the excited state)
 from the perturbed eigenvalue problem \eqref{eigenpert}. 
For this we use the result in the following lemma, whose derivation follows \REF{golub}. 

\begin{lemma}\label{lem:eigenproblemdiff}
Assume that the matrix $A\in\mathbb{C}^{n\times n}$ has $n$ distinct
eigenvalues, $\mu_j$, $0\leq j \leq n-1$. Let $l_j$ and $r_j$ be 
left and and right
eigenvectors of $A$, i.e.\ $l_jA=\mu_j l_j$, $Ar_j=\mu_j r_j$,
satisfying 
\begin{equation*}
l_j\cdot r_k=\begin{cases}
1, & \text{ if } j=k,\\
0, & \text{ if } j\neq k.
\end{cases}
\end{equation*}
For
sufficiently small  matrices $\beta\in\mathbb{C}^{n \times n}$
the eigenvectors and eigenvalues of the perturbation matrix $A+\beta$
are differentiable functions of $\beta$. If $r_\beta$ and $\mu_\beta$ denote the
eigenvector and eigenvalue to $A+\beta$ that 
equal the eigenpair $(r_0,\mu_0)$ to the matrix $A$ for $\beta=0$,
we have that
\begin{equation}\label{eq:rbeta}
l_j\cdot(r_\beta-r_0) = -\frac{l_j\cdot \beta r_0}{\mu_j-\mu_0}+o(|\beta|).
\end{equation}
and
\begin{equation}\label{eq:mubeta}
\mu_\beta-\mu_0 = l_0\cdot \beta r_0 + o(|\beta|).
\end{equation}
\end{lemma}
\begin{proof}
That the eigenvectors and eigenvalues of $A+\beta$ are differentiable
functions of $\beta$ in a neighborhood of $\beta=0$ follows directly
by differentiation of the relations
\begin{equation*}
(A+\beta)r=\mu r,\quad l(A+\beta)=\mu l.
\end{equation*}

Let $B:=|\beta|^{-1}\beta$, where $|\beta|$ is a matrix norm, e.g.,\ the Euclidean operator norm. 
Let $r(\gamma)$ and $\mu(\gamma)$ be the perturbed normalized eigenvector and
eigenvalue corresponding to $r_0$ and $\mu_0$ given by $r(0)=r_0$ and
$\mu(0)=\mu_0$, and
\begin{equation}\label{eq:eigengamma}
(A+\gamma B)r(\gamma)=\mu(\gamma)r(\gamma).
\end{equation}
Differentiation of \eqref{eq:eigengamma} at $\gamma=0$ gives
\begin{equation}\label{eq:eigendiff}
Ar'(0)+Br(0)=\mu'(0)r(0)+\mu(0)r'(0).
\end{equation}
Using that $r(0)=r_0$ and $\mu(0)=\mu_0$ and multiplying by $l_j$ from
the left, we get
\begin{equation*}
(\mu_j-\mu_0)l_j\cdot r'(0) = l_j\cdot Ar'(0)-\mu_0l_j\cdot r'(0) = -
  l_j\cdot B r_0 + \mu'(0) l_j\cdot r_0 = - l_j\cdot B r_0, \quad \text{for $j\neq 0$},
\end{equation*}
by the orthogonality between left and right eigenvectors.
The Taylor expansion 
$r(|\beta|)=r(0)+r'(0)|\beta| +o(|\beta|)$ gives \eqref{eq:rbeta}.

To prove \eqref{eq:mubeta} we take the scalar product with $l_0$ from
the left in \eqref{eq:eigendiff}, which gives $\mu'(0)=l_0\cdot B
r_0$. The Taylor expansion $\mu(|\beta|)=\mu(0)+\mu'(0)|\beta| +
o(|\beta|)$ gives \eqref{eq:mubeta}.
\end{proof}

 The combination of \eqref{eq:rbeta} in Lemma \ref{lem:eigenproblemdiff} and the assumtion that $r_\beta-r_0=\mathcal O(1)$, as $|\beta|\rightarrow 0+$, also when $|\mu_j-\mu_0|$ is small, implies that
\begin{equation}\label{scatt_mu}
\begin{split}
\ell_j\cdot (r_\beta-r_0) &=
\left\{
\begin{array}{cl}
- \frac{\ell_j\cdot \beta r_0}{\mu_j-\mu_0 }+ o(|\beta|)\COMMA & \mbox{ if } |\mu_j-\mu_0|^{-1}=o(|\beta|^{-1})\COMMA\\
\mathcal O(1)\COMMA & \mbox{ otherwise }\COMMA\\
\end{array}\right. \\
&= - \frac{\ell_j\cdot \beta r_0}{\mu_j-\mu_0 + c|\ell_j\cdot \beta r_0|\, \mbox{sign}(\mu_j-\mu_0)}+ 
o(|\beta|)  \\
&[\mbox{ and by the definition } (\mu_j-\mu_0)^\sharp:=\mu_j-\mu_0 + c|\ell_j\cdot \beta r_0|\, \mbox{sign}(\mu_j-\mu_0)]\\
&=: - \frac{\ell_j\cdot \beta r_0}{(\mu_j-\mu_0)^\sharp}+o(|\beta|) \COMMA
\end{split}
\end{equation}
 for some positive constant $c$, which determines the perturbation orthogonal to $r_0$.
 We will also use that this representation holds separately 
 for each perturbation $S^+$ and $S^-$. 
We do not
determine the perturbation in the $r_0$ direction, which does
 not contribute to the transition between states. 
We denote by $r_\beta^\perp- r_0$ the 
projection of the perturbation $r_\beta - r_0$ on the hyperplane
orthogonal to $r_0$.
Since  this plane is spanned by the left eigenvectors
$l_1,\ldots,l_{n-1}$, we have by \eqref{eq:rbeta} in Lemma
\ref{lem:eigenproblemdiff} that 
\begin{equation*}
r_\beta^\perp - r_0 = -\sum_{l_j\cdot r_0 \neq 0} \frac{l_j\cdot \beta
r_0}{\mu_j-\mu_0} + o(|\beta|).
\end{equation*}

{
Returning to the perturbed problem \eqref{eigenpert}, we recall from 
equation~\eqref{egenvardet} that the unperturbed eigenvalues are given by 
\[
\begin{bmatrix}
\mu_0\\
\mu_1\\
\mu_2\\
\mu_3
\end{bmatrix} 
= 
\begin{bmatrix}
\frac{1-R_1^\ell}{2} - \sqrt{\left(\frac{1-R_1^\ell}{2}\right)^2 +R_1^\ell -R_1^r}\\
\frac{1-R_1^\ell}{2} + \sqrt{\left(\frac{1-R_1^\ell}{2}\right)^2 +R_1^\ell -R_1^r}\\
\frac{1-R_2^\ell}{2} - \sqrt{\left(\frac{1-R_2^\ell}{2}\right)^2 +R_2^\ell -R_2^r}\\
\frac{1-R_2^\ell}{2} + \sqrt{\left(\frac{1-R_2^\ell}{2}\right)^2 +R_2^\ell -R_2^r}
\end{bmatrix},
\]
and the corresponding unperturbed eigenvectors are
\[
\ell_j = \begin{cases} \frac{(1, -(1-\mu_j)/R_1^r, 0, 0)}{(1-(1-\mu_j)^2/R_1^r)}  & \text{when } j=0,1,\\
                      \frac{(0, 0, 1, -(1-\mu_j)/R_2^r )}{(1-(1-\mu_j)^2/R_2^r)}  & \text{when } j =2,3.
\end{cases}
\quad \text{and } \quad 
r_j =\begin{cases} (1, 1-\mu_j, 0, 0) & \text{when }  j=0,1,\\
(0, 0, 1, 1-\mu_j ) & \text{when } j =2,3.
\end{cases}
\]}
Lemma \ref{lem:eigenproblemdiff} yields 
\begin{equation}\label{egen_p_E}
\mu_\beta-\mu_0=\ell_0\cdot \beta r_0 + o(|\beta|)
= (1-S_{11})\frac{1-(1-\mu_0)^2\frac{R_1^\ell}{R_1^r}}{1-(1-\mu_0)^2\frac{1}{R_1^r}} +o(\sqrt{p_d})
= o(\sqrt{p_d})\, .
\end{equation}
To have the perturbed eigenvalue equal to zero
 means that we start with an unperturbed eigenvalue $\mu_0=o(\sqrt{p_d})$ in order to obtain
the perturbed eigenvalue  $\mu_\beta=0$, satisfying $(A+\beta)r_\beta=0$.
Using the
unperturbed  eigenvector $r_0=(1,1+o(\sqrt{p_d}),0,0)$, with $\mu_0=o(\sqrt{p_d})$, yields a perturbation in the $\Phi_2$ component
from zero (in the third component of $r_0$) to
\begin{equation}\label{S_definition}
S_{21} \sum_{j=2,3} \frac{1-(1-\mu_j)^2\frac{R_2^\ell}{R_2^r}}{(\mu_j-\mu_0)^\sharp(1-(1-\mu_j)^2\frac{1}{R_2^r})} +o(\sqrt{p_d})
=\mathcal O(\sqrt{p_d})
\end{equation}
and its squared absolute value measures the  probability to be in the excited state, i.e.  the probability density $|\Phi_2|^2$, which integrated yields $p_{ex}$ defined in \eqref{p_e_def1}.

\begin{figure}[htbp] %
        \centering
        \subfigure[The probability $p_{ex}$]{\label{pex1D:subfig3}
		\includegraphics[width=0.8\textwidth]{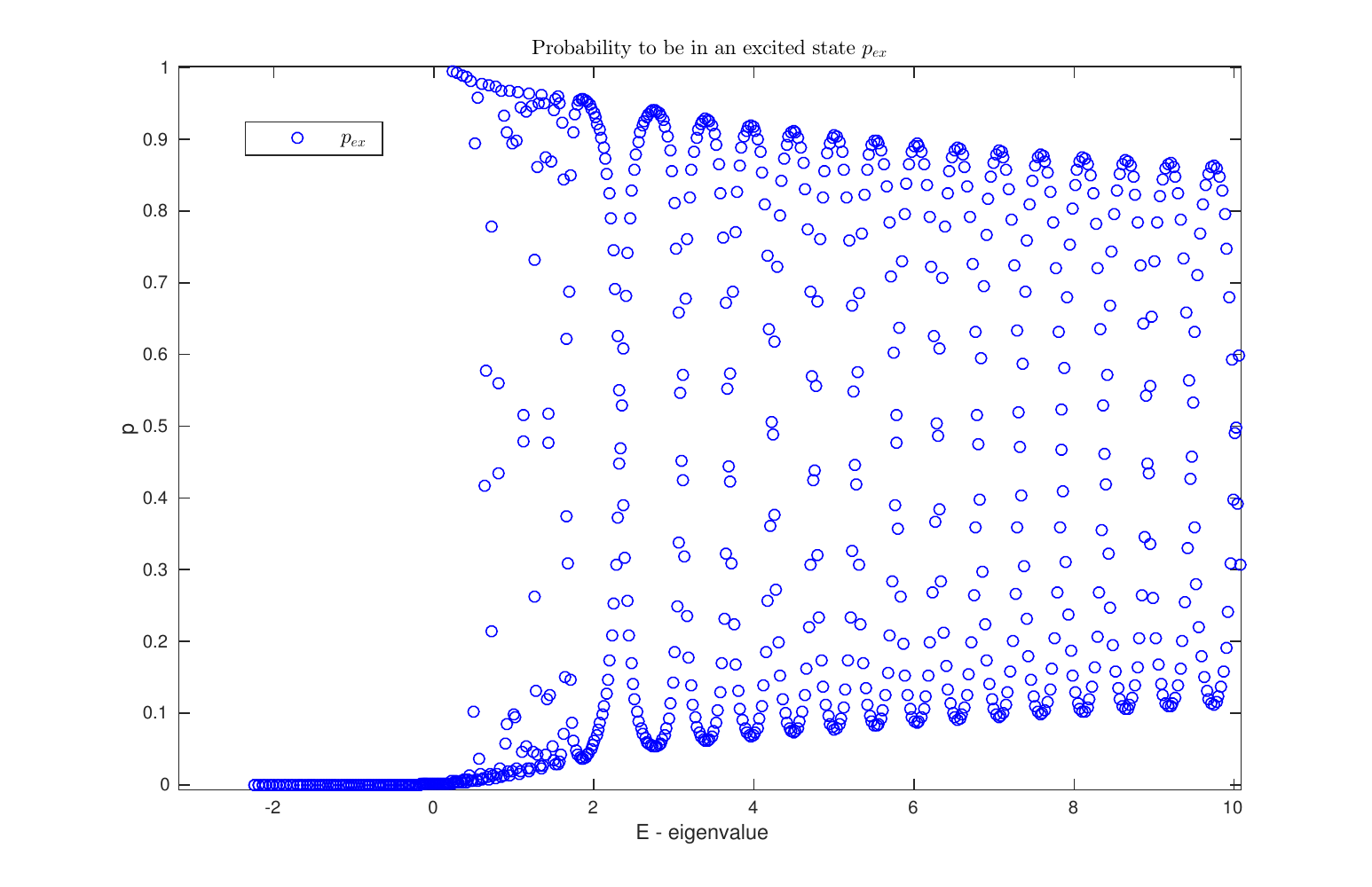}}

        \subfigure[The probability $p_{ex}$ (left panel) to be in an excited state compared to the estimate $\hat p_{ex}$ (right panel)]{\label{pex1D:subfig2}
		\includegraphics[width=0.8\textwidth]{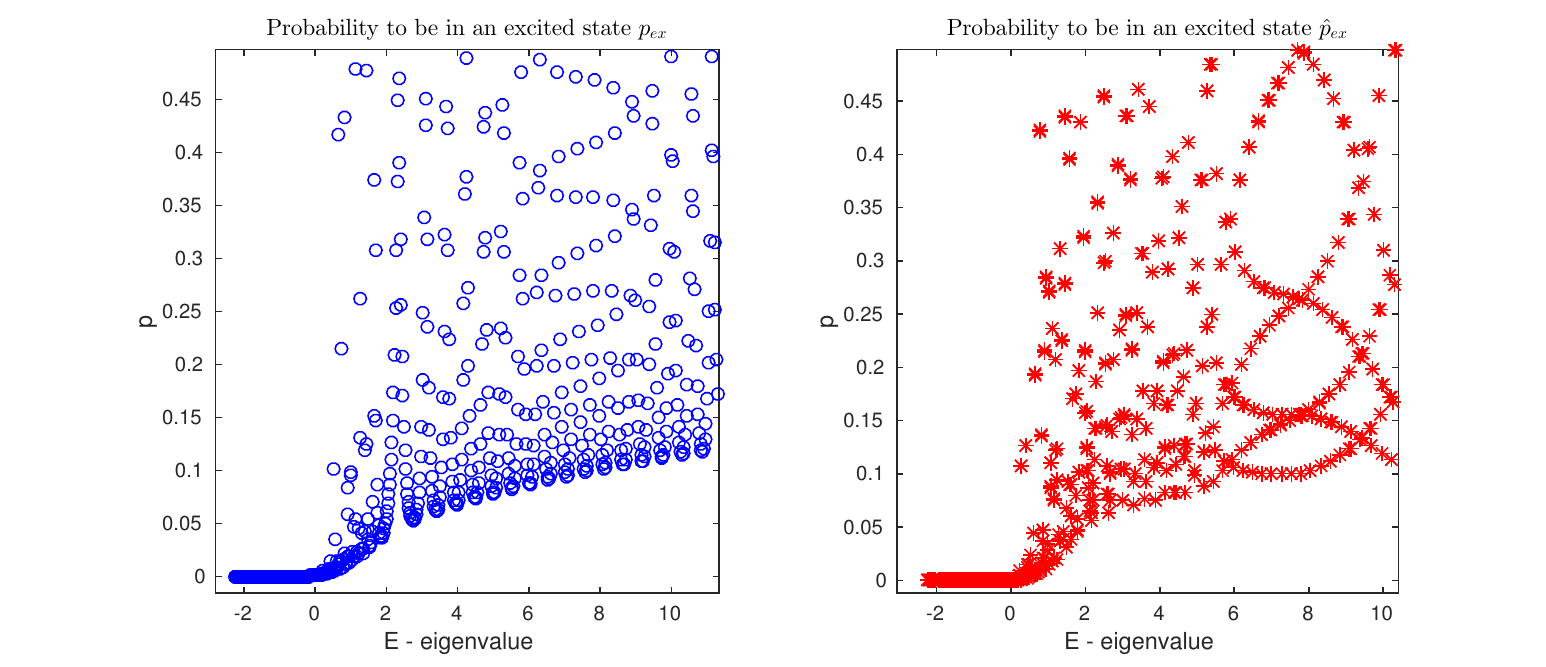}}
	
        \subfigure[The probability $p_{ex}$ ('o') to be in an excited state compared to the estimate $\hat p_{ex}$ ('*')]{\label{pex1D:subfig1}
		\includegraphics[width=0.8\textwidth]{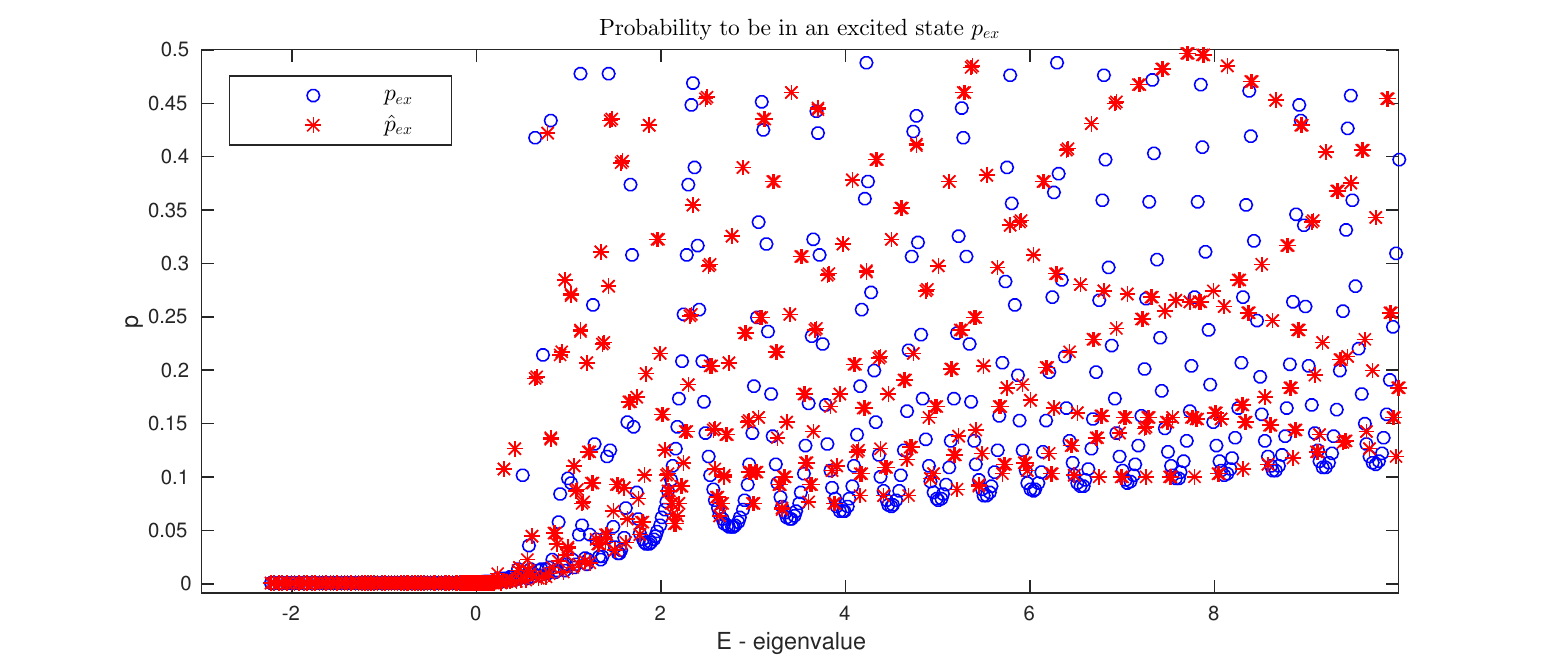}}

	\caption{Probabilities to be in the excited state,~$p_{ex}$, using 
  numerical approximation of \eqref{p_e_def1}, precisely defined in  
\eqref{eq:Pr-excited}, and the estimations,~$\hat p_{ex}$, using 
  the formula~\eqref{eq:pe-approx}, obtained by solving the discrete one-dimensional
  Schr\"odinger eigenvalue problem~\eqref{eq:discrete-schrod-1d} with~$M=1000$, $\delta = \sqrt{\frac{3}{5}} M^{-0.2}\approx 0.1946$,
  $C=0.09$, and mesh size $h =0.0001$.}
	\label{fig:pex-pexdetail1D}	
\end{figure}

 \begin{remark}\label{rem:resonance}{\rm
 Having a zero eigenvalue of  $A_1$  in \eqref{S_sum} means that the scalar Schr\"odinger eigenvalue problem
 \[
 \big(-\frac{1}{2M}\Delta + \lambda_1(X)\big)\Phi_1(X)  =E\Phi_1(X)
 \] 
 has a solution in the whole domain.
 We see from \eqref{S_definition} that when also $A_2$ has a zero eigenvalue $\mu_2=0$ we obtain a resonance in the sense that a small perturbation $\beta$ (where $|S_{21}|=\mathcal O(|\beta|)=\mathcal O(p_d^{1/2})$) yields
 a large change in the eigenvector $r_\beta=\Phi$, meaning that the probability to be in the excited state (which by \eqref{S_definition} is of the order $|\beta|^2/|\mu_2|^2$)  can be
 of order one even if the perturbation $\beta$ is tiny. Resonance means that the two scalar Schr\"odinger eigenvalue problems have a coinciding eigenvalue $E$, see Figure \ref{fig:pex-pexdetail1D}. An approximation of this 
 probability  is 
 \begin{equation}\label{pert_b}
 |\beta|^2/|\mu_2^\sharp|^2\approx p_d/(|\mu_2-\mu_0|^2+ c^2p_d). 
 \end{equation}
 }
 \end{remark}
 
 \subsubsection{Numerical test of the perturbation analysis}
A method to numerically test the validity of the perturbation analysis leading to the estimate \eqref{S_definition}, for the one dimensional model,
is to compare the probability to be in excited states $p_{ex}$ in %
\eqref{p_e_def1} (determined from numerical approximation of the 
Schr\"odinger equation \eqref{eq:discrete-schrod-1d}) with the expression \eqref{pert_b} derived from  \eqref{S_definition} in the perturbation analysis. If $p_{ex}$ shows some similarity with \eqref{pert_b} the
perturbation analysis is in some sense justified.
The difference $\mu_2-\mu_0$, of the eigenvalues   in \eqref{pert_b}, can be obtained from numerical approximation of a difference $E_{2,k'}-E_{1,k}$, of the eigenvalues 
of the decoupled system \eqref{diagonal_sch},
$  \big(\lambda_i(X)-\frac{1}{2}M^{-1}\frac{\partial^2}{\partial X^2}\big)\Phi_{i,k}(X)  = E_{i,k}\Phi_{i,k}(X)$, $i=1,2$,
using the same centered difference approximation as in \eqref{eq:discrete-schrod-1d}. Let $\mathcal E_i=\{E_{i,k}\ |\ k=1,2,3,\ldots\}$ be the
set of all $E_i$ eigenvalues, for $i=1,2$.
By the approximation 
\[|\mu_2-\mu_0|^2\approx C'|E_{1,k}-E_{2,k'} |^2 
\approx C^{-1} |E_{1,k}-E_{2,k'} |^2/|E_{1,k} -E_{1,k+1}|^2,\] 
where $k$ and $k'$ is chosen such that $E_{1,k}$ and $E_{2,k'}$ are the closest to $E$  and 
using two  constants $C$ and $C'= C^{-1}/|E_{1,k}-E_{1,k+1}|^2$,
we obtain the following estimate of the probability 
to be in the excited state
\begin{equation}\label{eq:pe-approx}
  \hat p_{ex} := \frac{p_{d}}{C^{-1}|E_1^l-E_2^l|^2/ |E_1^l-E_1^m|^2 + p_{d}}\COMMA
\end{equation}
with $C$~being a constant, $E$~being an eigenvalue of the discrete two-state Schr\"odinger 
equation corresponding to \eqref{eq:pe-approx} and
\[
\begin{split}
E_1^l &:=\underset{E_1\in\mathcal E_1}{\operatorname{argmin}} |E_1-E|,\quad
E_1^m := \underset{ E_1\in\mathcal E_1\,  \&\, E_1\ne E^l_1}{\operatorname{argmin}} |E_1-E^l_1|,\\
 E_2^l &:= \underset{E_2\in\mathcal E_2}{\operatorname{argmin}} |E_2-E_1^l|,\quad
p_{d} := \exp(-\pi\delta^2\sqrt M / \sqrt{2(E-\lambda_-(0))})\, .
\end{split}
\]

\begin{figure}[htbp] %
        \centering
        \subfigure[]{\label{1D:subfig1}
                \includegraphics[width=0.45\textwidth]{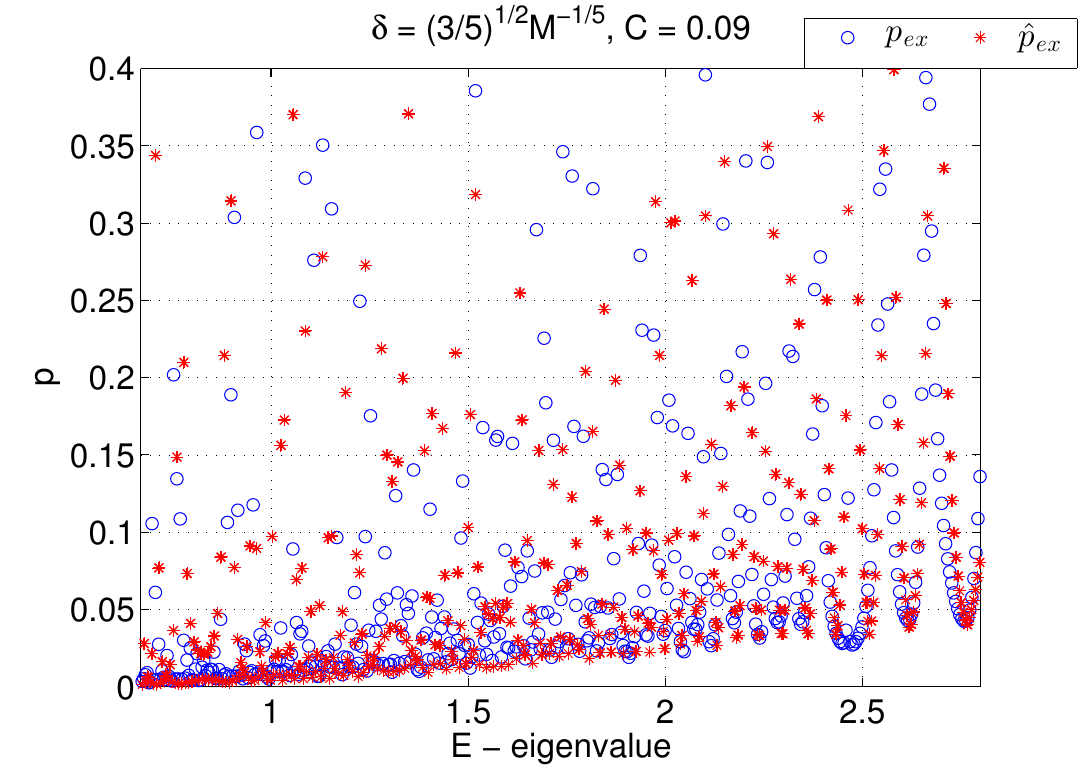}}
          \qquad
        \subfigure[]{ \label{1D:subfig2}
                \includegraphics[width=0.45\textwidth]{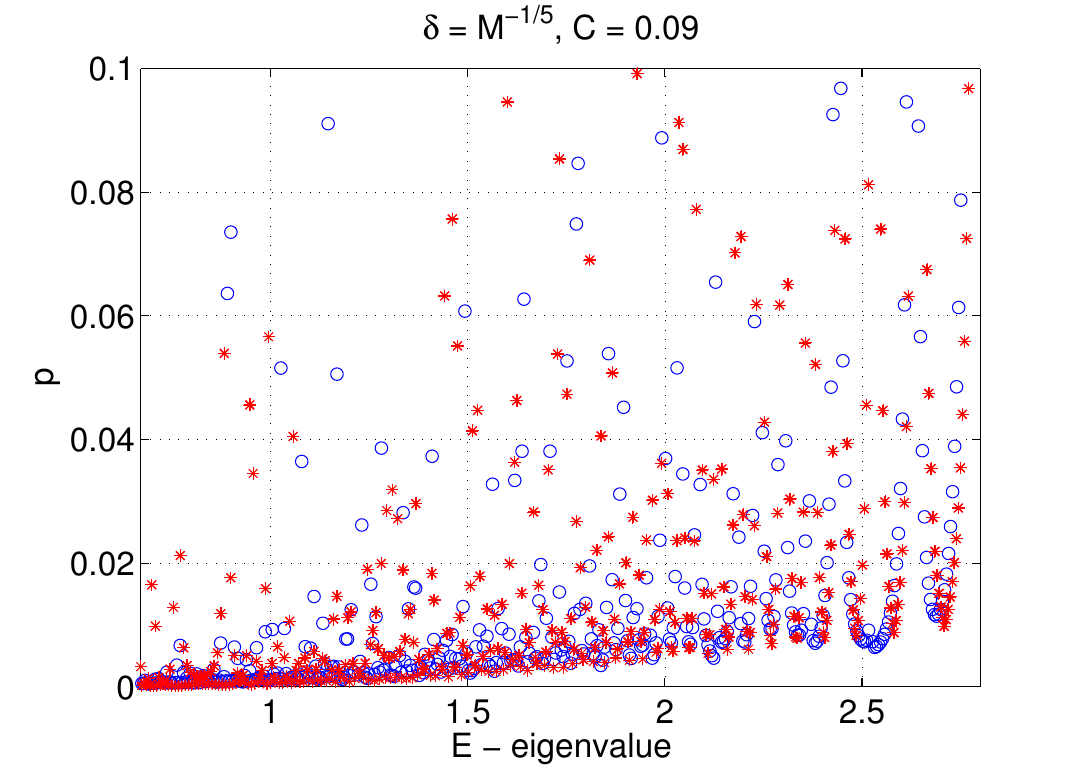}}

        \subfigure[]{ \label{1D:subfig3}
                \includegraphics[width=0.45\textwidth]{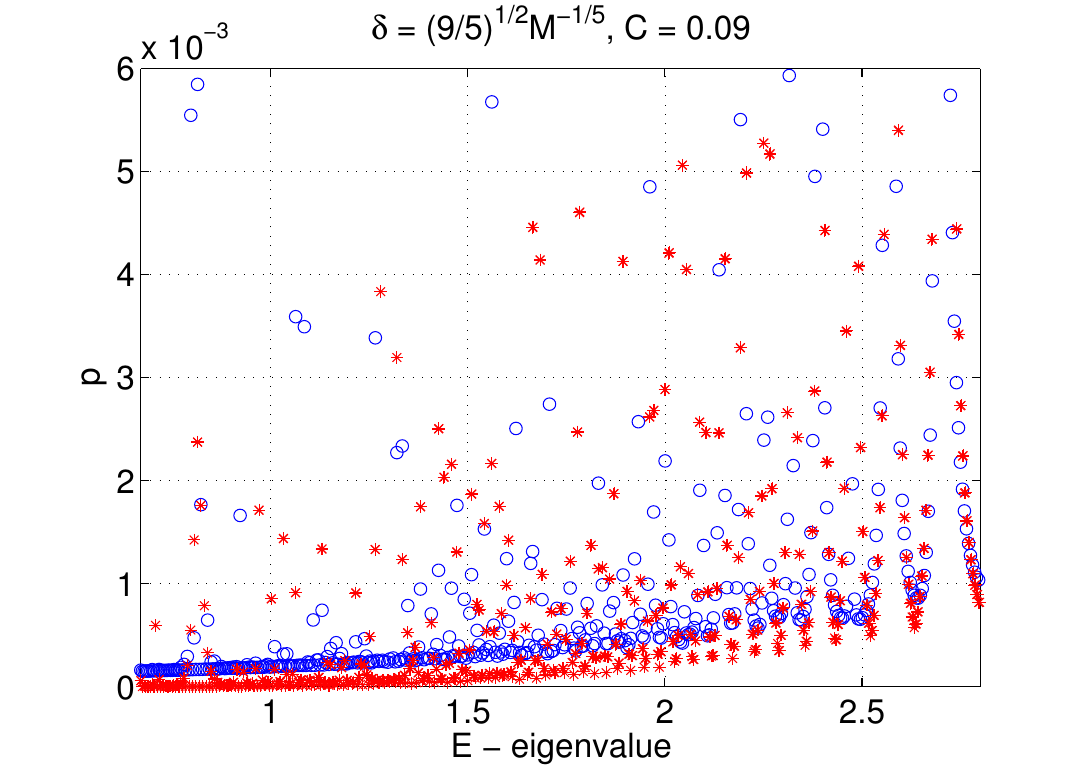}}
        \qquad
        \subfigure[]{ \label{1D:subfig4}
                \includegraphics[width=0.45\textwidth]{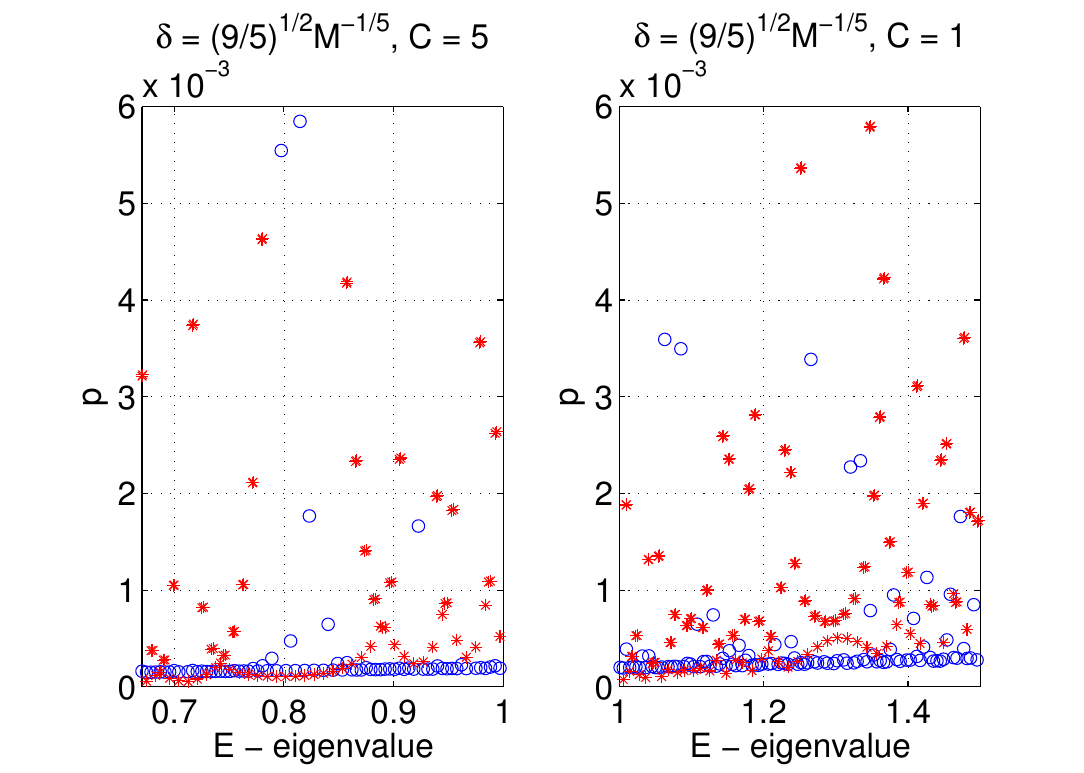}}

  \caption{Probabilities to be in the excited state,~$p_{ex}$, using 
  numerical approximation of \eqref{p_e_def1}, precisely defined in  
\eqref{eq:Pr-excited}, and the estimations,~$\hat p_{ex}$, using 
  the formula~\eqref{eq:pe-approx},
obtained by solving the discrete one-dimensional
  Schr\"odinger eigenvalue problem~\eqref{eq:discrete-schrod-1d} with~$M=20000$.}
  \label{fig:computed-estimated}
\end{figure}

\subsubsection*{Conclusions.}
Figures~\ref{fig:pex-pexdetail1D} and~\ref{fig:computed-estimated} %
show that the estimation of the probabilities,~$\hat p_{ex}$,
to be in the excited state, obtained using the 
formula~\eqref{eq:pe-approx}, 
and %
the numerically computed 
probabilities,~$p_{ex}$, to be in the excited state, obtained using 
the formula~\eqref{eq:Pr-excited} from the solution of the discrete one dimensional 
Schr\"odinger eigenvalue problem~\eqref{eq:discrete-schrod-1d},
have  similar qualitative behavior in the following two aspects: 
(1) the minimal value of $p_{ex}$ and $\hat p_{ex}$ are similar, and 
(2) when $E$ is close to a resonance, seen for $\hat p_{ex}$ in \eqref{eq:pe-approx} as $|E_1^l - E_2^l | \ll | E_1^l-E_1^m |$, 
both $p_{ex}$ and $\hat p_{ex}$ become larger. 
By also adjusting the constant $C$ in \eqref{eq:pe-approx} the behavior can also be quantitatively similar, see Figure \ref{fig:pex-pexdetail1D} and \ref{fig:computed-estimated}.
That is the perturbation analysis
seems to capture this resonance phenomenon, at least qualitatively. 

\subsubsection{The perturbed eigenvalue problem in the multi-dimensional case}\label{sec_multiD}
The main idea to extend the one dimensional perturbed eigenvalue problem to multiple dimensions is
to write the Schr\"odinger wave function as a phase-space integral of highly oscillatory
functions 
using the Fourier-Bros-Iagolnitzer transform (FBI-transform). { 
This integral can be viewed as an integral over all hyperplanes $I$  in $\rset^{3N}$ defined in \eqref{I_plane}.
Therefore 
each  FBI-mode can be initial data
for a WKB solution  \eqref{wkb_ekv1}
to the right of a hyperplane $I$  in $\rset^{3N}$ and a WKB solution to the left of the hyperplane. The condition to have a differentiable WKB solution
across the hyperplane then uses the continuity condition \eqref{eigenpert}, as the one dimensional setting,
for each point in the hyperplane.
}
The direction of the oscillations in the FBI-mode generates through the WKB-method  a molecular dynamics path and an amplitude function.
The amplitude function solves, along the path, a transport equation that is well approximated
by Ehrenfest dynamics { \eqref{eq:ehrenfest}} and can consequently be computed similarly as the molecular dynamics paths but with smaller time steps, see Sections \ref{zener_ehren} and \ref{sec_num}.
 The integral of WKB-solutions over all hyperplanes contribute
through their amplitude functions to the transition, $S$.
We will now explain this extension to multiple dimensions more precisely
and then in the next section show how to computationally approximate the probability to be in excited states $p_{ex}$.

We will use the FBI-transform. %
We could use the standard semi-classical Fourier transform as the initial data for the WKB-method but
the FBI-transform has the advantage of giving the high frequency content locally in $X$, i.e. microlocally,  which yields a more
accurate WKB Ansatz. The important property of the FBI transform
\begin{equation}\label{fbi_def}
T\varphi(Y,P):=\underbrace{(2^{1/3}\pi)^{-3N/4}M^{9N/8}}_{=:\alpha_M}\int_{\rset^{3N}} e^{iM^{1/2}(Y-X)\cdot P - |Y-X|^2M^{1/2}/2} \varphi(X) dX
\end{equation}
is the identity 
\begin{equation}\label{fbi_iso}
\varphi=T^*T\varphi
\end{equation}
for $\varphi\in L^2(\rset^{3N})$, where the adjoint operator $T^*$ is defined by
\[
T^*\phi(X)= \alpha_M\int_{\rset^{6N}} 
e^{-iM^{1/2}(Y-X)\cdot P - |Y-X|^2M^{1/2}/2} \phi(Y,P) dY dP
\]
e.g.\ for all Schwartz functions $\phi$ on $\rset^{6N}$, see \REF{martinez}.
Therefore the integral representation
\[
\begin{split}
\Phi(X) &=T^*T\Phi(X)\\
&=\alpha_M
\int_{\rset^{6N}} 
e^{-iM^{1/2}(Y-X)\cdot P - |Y-X|^2M^{1/2}/2} T\Phi(Y,P) dY dP\\
\end{split}
\]
yields suitable boundary data
on  hyperplanes (e.g. $X_1=0$)
for the WKB-method. %
We will, for  each point $(Y,P)=:Q$ (in the classically allowed region), use a WKB-function \eqref{wkb_form}.
In the case of caustics this single WKB-function is replaced by a finite sum of WKB-functions
based on  the Legendre transform $\theta^*(P)$ of $\theta(X)$, as explained in Section \ref{sec:caustics} and \ref{stat_phase_sec},
\begin{equation}\label{wkb_decomp}
\Phi_Q(X)\sim\sum_{\{P^*: \nabla_{P}\theta_Q^*(P^*)= X\}} \phi_{Q}(X;P^*)e^{iM^{1/2}\theta_{Q}(X;P^*)}\, .
\end{equation}
This WKB solution $\Phi_Q$ solves the Schr\"odinger eigenvalue problem \eqref{schrodinger_stat} in the two domains to the left and right of the  hyperplane, with  the boundary condition  constructed so that the WKB-mode $\Phi_Q(X)$ is equal to 
the FBI-mode \[\alpha_M %
e^{-iM^{1/2}(Y-X)\cdot P - |Y-X|^2M^{1/2}/2} T\Phi(Y,P)\]
for $X-Y$ in the hyperplane $L_P:=\{q\in \rset^{3N}\, | \, q\cdot P=0\}$ orthogonal to $P$:
\begin{equation}\label{wkb_data}
\begin{split}
q\cdot \nabla_{}\theta_Q(X)&= q\cdot P=0\COMMA\quad \mbox{for $q\in L_P$ and $X-Y\in L_P$}\\
\phi_Q( X;P^*) &= \alpha_M T\Phi(Q)e^{-|X-Y|^2M^{1/2}/2}\COMMA\quad \mbox{for $X-Y\in L_P$}\, .
\end{split}
\end{equation}
The construction of the WKB-solutions in the case of caustics 
uses asymptotic convergence as $M\rightarrow\infty$, described by the asymptotically equal sign $\sim$, see \REF{zworski,maslov} and \eqref{caustic_expansion}.
The obtained decomposition \begin{equation}\label{phi_Q_ekv}
\Phi(X)=%
\int_{\rset^{6N}}\Phi_Q(X)
dQ\end{equation}
determines, for each $Q$, %
the corresponding operators $R^{r,\ell}$.
The total perturbation comes from transitions from all $\phi_Q$ along the paths:
 for each $Q=(Y,P)$ %
we let the function $\phi_Q(X_t)$ solve the transport equation \eqref{schrod_first} along
the  WKB-path $\{X_t\}_{t=0}^\infty$
and initialize $\phi_Q(X_\tau)$ to the ground state 
$\Psi_0(X_\tau)$ 
each time $X_\tau$ is in the hyperplane 
$(X_\tau-Y)\cdot P=0$, generated by $Q$.
The perturbed wave function is then as in \eqref{transit_element}
determined by
\[
\langle S(Q,X),\Psi_i(X)\rangle =\frac{\langle \phi_Q(X),\Psi_i(X)\rangle}{ \langle \phi_Q(X),\Psi_0(X)\rangle}\, ,
\]
which depends on $Q$,
by relating $\psi_{1+}$ in \eqref{transit_element} with $\langle \phi_Q, \Psi_0\rangle$
and $\psi_{2+}$ with $\langle\phi_Q,\Psi_1\rangle$ in the case of two states. %
Integrating over all $Q$ yields 
the total perturbation, as in \eqref{S_sum},
\[
\int_{\rset^{6N}} S(Q,X)\Phi_Q(X) dQ 
\] at the point $X$.
The corresponding perturbed eigenvalue problem \eqref{eigenpert} for the wave function $\Phi_Q$
implies, as in  \eqref{S_definition}, the following larger
change in the excited state
\begin{equation}\label{phi_exp}
\begin{split}
\Phi^\perp(X)&=\Phi(X)-\langle \Psi_0(X),\Phi(X)\rangle\Psi_0(X)\\
&=\int_{\rset^{6N}}
\underbrace{S_{21}(Q,X)}_{=p_d^{1/2}(Q,X)} \Phi_Q(X) \sum_{j=2,3} \frac{1-(1-\mu_j)^2\frac{R_2^\ell}{R_2^r}}{(\mu_j-\mu_0)^\sharp(1-(1-\mu_j)^2\frac{1}{R_2^r})}dQ \\
&\qquad 
+\mathcal O(|\int_{\rset^{6N}}
S_{21}(Q,X)\Phi_Q(X) d Q|^2)\end{split}
\end{equation}
at the point $X$, where  $\mu_j$ and $R^{r,\ell}_2$ depend on $Q$. 
The index $1$ in $S$ now corresponds to the component in the ground state
and the index $2$ to be in an excited state, i.e.,
\[
\Phi(X)=\gamma_1(X)\Psi_0(X) + \gamma_2(X)\Psi^\perp(X)\COMMA
\]
and $\gamma_j(X) \in \mathbb C$ with the orthogonality conditions
\[
\begin{split}
\langle \Psi_0,\Psi^{\perp}\rangle & =0, \quad
\langle \Psi_0,\Psi_0\rangle =1,\quad 
\langle \Psi^\perp,\Psi^{\perp}\rangle  =1\, .\\
\end{split}
\]

The eigenvalues of the unperturbed problem \eqref{egenvardet} vary highly %
due the oscillatory functions $R^{r,\ell}$, with the phase $M^{1/2}\theta^\pm(X)$. 
Figure \ref{fig:computed-estimated} illustrates a consequence of this variation and indicates that mean values of excitation probabilities corresponding to eigenvalues in a neighborhood 
are more stable than individual excitations probabilities. 
Let us therefore consider a local average of the probability to be in the excited state based on  a
simple model of the smallest eigenvalue
$\mu_2%
\sim \tan (M^{1/2}(\theta_2^+-\theta_2^\ell))$ %
(in the excited component) %
where we assume $|\mu_2|$ has a bounded density, $\rho_Q$, defined on $\rset$. %
In the case of one space dimension we can fix a point $X\in\rset$ (in the classically allowed region) where $\mu_2$ is evaluated, that cuts the spatial domain into two. The average then corresponds to an ensemble of Schr\"odinger eigenvalues localized around $E$, related to a local mean value of excitation probabilities in Figure   \ref{fig:computed-estimated}. We denote this local average by $\mathbb A$.
In higher dimensions we can similarly fix a co-dimension one cutting surface and evaluate $\mu_2$ on the surface.
Let us now study the average excitation probability in this model.
Using that $(\mu_2-\mu_0)^\sharp$ yields transition probabilities bounded by one for each $Q$ as in \eqref{scatt_mu}
and that outside a compact set (corresponding to the classically allowed region) the factor $|S_{21}(Q,X)\Phi_Q(X)|$ is negligible small, we have
\[
\begin{split}
\mathbb A[\langle \Phi^\perp(X),\Phi^\perp(X)\rangle ]
&\le C \mathbb A\big[|\int_{\rset^{6N}}
\frac{|S_{21}(Q,X) \Phi_Q(X)|}{|(\mu_2-\mu_0)^\sharp| } dQ|^2\big]\\
&\le C \mathbb A\big[\int_{\rset^{6N}}\frac{|S_{21}(Q,X) \Phi_Q(X)|^2}{|\mu_2-\mu_0|^2 + c^2|S_{21}(Q,X)|^2} d Q\big] \\
&= C\int_{\rset^{6N}} \int_\rset \frac{|S_{21}(Q,X) \Phi_Q(X)|^2}
{|\mu_2-\mu_0|^2 + c^2|S_{21}(Q,X)|^2}\rho_Q(\mu_2) d\mu_j d Q\\
&= C\int_{\rset^{6N}} |S_{21}(Q,X)| |\Phi_Q(X)|^2\int_\rset\frac{1}
{\frac{|\mu_2-\mu_0|^2}{ c^2|S_{21}(Q,X) |^2}+1}\rho_Q(\mu_2) 
d\frac{\mu_2}{c^2|S_{21}(Q,X)|} dQ\\
&\le C \int_{\rset^{6N}}{|S_{21}(Q,X)| |\Phi(X)|^2}  \|\rho_Q\|_{L^\infty(\rset)} dQ\  \int_\rset \frac{dx}{x^2+1}
\\
&=\mathcal O\big(
\int_{\rset^{6N}}{|S_{21}(Q,X)| |\Phi(X)|^2} dQ\big)
\, .\\
\end{split}
\]
We conclude that in a model of local ensemble averages of eigenvalues $E$, assuming that  
for each $X$ the eigenvalue difference $|\mu_2-\mu_0|$ has a bounded density, then the average
probability to be in the excited state has the bound
\[
\mathbb A[p_{ex}]=\mathcal O(\int_{\rset^{3N}}\int_{\rset^{6N}}
{|S_{21}(Q,X)||\Phi_Q(X)|^2} %
d Q d X)\, ,
\]
which, by \eqref{phi_exp} can be written as the average of the square root of the dynamic transition probability $p_d$
\begin{equation}\label{prop2}
\begin{split}
\mathbb A[p_{ex}] &=\mathcal O(\int_{\rset^{3N}}\int_{\rset^{6N}}
p_d^{1/2}(Q,X)|\Phi_Q(X)|^2 dQd X)\, .\\
\end{split}
\end{equation}

\subsubsection{Ergodic computation of probabilities to be in excited states}\label{sec_erg_p_e}
\begin{figure}[hbp]
\centering
\includegraphics[width=0.8\textwidth]{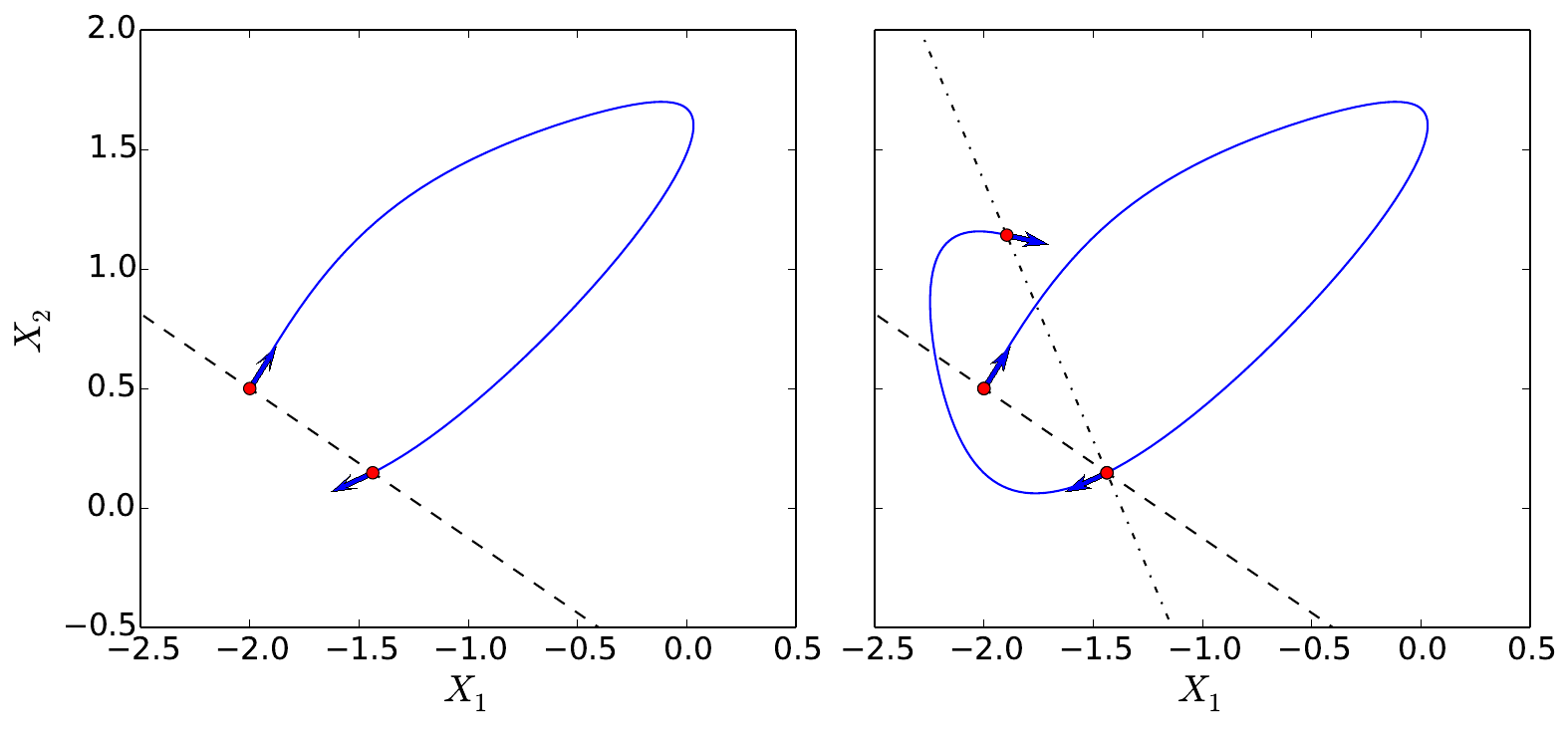}
   \caption{Illustration of the
       sampling of hyperplanes and initialization times described in
       part 1 of Algorithm~\ref{alg:compute-pe-md}. 
       Left figure: the 2D Born-Oppenheimer
       molecular dynamics trajectory
       $X_t$ (blue line) with initial conditions $X_{T_0}=(-2,0.5)$
       and $\dot X_{T_0} =P_{T_0}$ is simulated until 
 the time $t=T_1$ when it crosses the hyperplane $P_{T_0}^\bot$ (dashed
 line). Right figure: a new hyperplane $P_{T_1}^\bot$ (dash-dotted line) 
is constructed and the trajectory of $X_t$ is simulated further until
the time $t=T_2$ when it crosses that hyperplane. The sampling
procedure is iterated until the final time is reached. }
\label{fig:init_times}
\end{figure}

The FBI-transform satisfies
\[
  \|T\Phi\|_{L^2(\rset^{6N})}=\|\Phi\|_{L^2(\rset^{3N})}=1
\]
and $T\Phi$ concentrates on the phase-space set $H_0(Y,P)=|P|^2/2+\lambda_0(Y)=E$ in the limit
as $M\rightarrow\infty$ and $p_{ex}\rightarrow 0+$, see \REF{martinez}. When the dynamics is ergodic, the phase-space
measure is in addition uniform on the set $E<H_0(Q)<E+\gamma$ as $\gamma\rightarrow 0+$, cf. \REF{reed_simon}. The WKB functions $\Phi_Q(X)$ behave similarly as $T\Phi$ locally, since their initial conditions are given by $T\Phi$ and
they solve the Schr\"odinger equation. We may write $S_{21}(Q,X)=:S_{21}(Q;X,P)$ where
$P$ is the momentum for the path at the position $X$ that started in the plane given by $Q=(Y',P')$, since this $P$ is a function of $Q$ and $X$.
Therefore we approximate $p_{ex}$ by the following time average of $p_d^{1/2}$
\begin{equation}\label{ergod_p_scatt}
\begin{split}
p_{ex}\approx \lim_{\gamma\rightarrow 0+} \int_{\rset^{3N}} \frac{\int_{E<H_0(Q)<E+\gamma} |S_{21}(Q;X,P)| dQ}{
\int_{E<H_0(Q)<E+\gamma} dQ} dX
& =
\lim_{\gamma\rightarrow 0+} \int_{\rset^{3N}} \frac{\int_{E<H_0(Q)<E+\gamma} |S_{21}(X,P;Q)| dQ}{
\int_{E<H_0(Q)<E+\gamma} dQ} dX\\
& =\lim_{T\rightarrow\infty} \int_{\rset^{3N}} \int_{0}^T \underbrace{|S_{21}(X,P;X_t,P_t)|}_{=p_d^{1/2}(X;X_t)} \frac{dt}{T} dX\COMMA
\end{split}
\end{equation}
where we used that the transition probability $S_{21}(Q;X,P)=S_{21}(X,P;Q)$ is symmetric by reversing time $t$ to $-t$.
In fact, along a given path $\{ X_t\}_{t=0}^T$ 
also the momentum $P'$, at position $Y'$, is determined by the position $Y'$ and $X_t$ so that 
$p_d^{1/2}(Y';X_t)=S_{21}(Y',P';X_t,P_t)$ is well defined.
The relation \eqref{ergod_p_scatt} means that we  sample square roots, $p_d^{1/2}(X;X_t)$, of transition probabilities along the molecular dynamics path, normalized for each passage through a hyperplane, and then take the average of all  hyperplanes. In our computations we sample hyperplanes by means of the phase space trajectory of ground state Born-Oppenheimer molecular dynamics, cf.~Algorithm~\ref{alg:compute-pe-md} and Figure~\ref{fig:init_times}.

\begin{algorithm}[hbp]
\caption{Approximations of the  probabilities, $p_{ex}$,  to be in excited states using molecular dynamics in 2D based on \eqref{ergod_p_scatt}}
\label{alg:compute-pe-md} 
\begin{algorithmic}
\STATE {\bf Input: } Energy $E$; potential functions $\VOPER$; mass $M$; time $T$;  
   initial position $X_0$, initial momentum $P_0$.
\STATE {\bf Output: } Approximated  probability $\hat{\hat{p}}_{ex}$ to be in excited states.

\medskip

\STATE {\bf 1.} Sampling of initialization times $\{T_n\}_{n}$ and hyperplane
coordinates $\{ \overline X_{T_n}, \overline P_{T_n} \}_n$:

\STATE Set $t\gets 0, n\gets0, T_n\gets 0$ and $(\overline X_{T_n}, \overline P_{T_n}) 
\gets(X_0, P_0)$ and define
the hyperplane $\overline P_{T_n}^\bot :=\{ X \in \R^2 | X \cdot \overline P_{T_n} =0\}$.

\WHILE{$T_n<T$}

\STATE Simulate the ground state Born-Oppenheimer molecular
dynamics for $(X_t,P_t)$ until $X_t$ crosses the plane $\overline P_{T_n}^\bot$.

\STATE At the crossing time, set $n\gets n+1, T_n \gets t$ and $(\overline X_{T_n}, \overline P_{T_n}) \gets (X_t, P_t)$ and define
the hyperplane $\overline P_{T_n}^\bot :=\{X \in \R^2 | X \cdot \overline P_{T_n} =0\}$.

\ENDWHILE

\medskip

\STATE {\bf 2.} Solve the Ehrenfest molecular dynamics using the St\"ormer-Verlet method with the entries in $\{T_k,\overline X_{T_k},\overline P_{T_k}\}_{k=0}^n$ as input parameters and compute $\hat{\hat{p}}_{ex}$ :
\STATE $temp \gets 0$
\FOR{$k = 0$ to $ n-1$}
  \STATE $t_1 \gets T_k$
  \STATE $t_2 \gets T_{k+1}$
  \STATE $X_{t_1} \gets \overline X_{T_k}$
  \STATE $P_{t_1} \gets \overline P_{T_k}$
  \STATE $\psi_{t_1} = \Psi_-(X_{t_1})$
  \STATE $\{X_t, P_t, \psi_t\} \gets \mbox{Ehrenfest dynamics path obtained from time $t_1$ to $t_2$ with initial data $X_{t_1}, P_{t_1}, \psi_{t_1}$}$
  \STATE $temp \gets temp + \int_{t_1}^{t_2} |\langle \psi_t, \Psi_+(X_t) \rangle| dt$ 
\ENDFOR
\STATE $\hat{\hat{p}}_{ex} \gets temp/T_n$
\end{algorithmic}
\end{algorithm}

\clearpage

\section{Numerical examples}\label{sec_num}
The purpose of this section is to present two simple model problems,
where the Schr\"odinger eigenvalue solution can be studied computationally
and compared to
the molecular dynamics approximation.
In the following subsections we 
 \begin{itemize}
 \item show that the dynamic transition probability $p_d$ in \eqref{p_E_def} can be determined by numerical solution of Ehrenfest dynamics,
 \item verify assumption \eqref{g_bar} on exponential convergence rate in finite time, %
 \item compare numerical approximations of observables from
 the Schr\"odinger equation and Born-Oppenheimer molecular dynamics,
 \item  compare approximations of the probability to be in excited states $p_{ex}$ from the Schr\"odinger equation with the molecular dynamics approximation $\hat{\hat p}_{ex}$ from Algorithm 1.
 \end{itemize}
\subsection{Model 1: A one dimensional problem}\label{sec:1d-example}
We consider the one dimensional, time-independent Schr\"odinger equation~\eqref{schrodinger_stat} 
with the heavy-particle coordinate~$X \in \mathbb{R}$, two electron states $J=2$
and the Hamiltonian operator,~$\HOPER(X)$, defined by
\begin{equation}\label{eq:schrod-hamiltonian}
\HOPER(X) := \VOPER(X)-\frac{1}{2}M^{-1}\frac{\partial^2}{\partial X^2}
\end{equation}
with the potential operator $\VOPER$ defined by the matrix
\begin{equation}
\VOPER(X) := \left[\begin{array}{cc}
                               X+r(X) & \delta \\
                               \delta & -X+r(X)
                        \end{array}\right]\COMMA \label{eq:potential-1d}
\end{equation}
where the parameter $\delta$ is a non-negative constant and the function~$r:\mathbb R \rightarrow \mathbb R$ is given by
$$r(X) := \left\{\begin{array}{ll}
                                (a_l-X)^2 & \mbox{if $X<a_l$}\COMMA \\
                               (X-a_r)^2 & \mbox{if $X>a_r$}\COMMA \\
			0 & \mbox{otherwise}.
                        \end{array}\right.$$
For each $X$ the potential $\VOPER$ defines the eigenvalue problem~$\VOPER(X)\Psi_\pm(X) = \lambda_\pm(X)\Psi_\pm(X)\COMMA$ 
with the two eigenvalues~$\lambda_\pm(X)\in \mathbb{R}$ and two eigenvectors~$\Psi_\pm(X)\in \mathbb{R}^2$. 
The eigenvalues are given by~$\lambda_\pm (X) = r(X) \pm \sqrt{X^2+\delta^2}$,
and $\Psi_-(X)$ and $\Psi_+(X)$ respectively denote the ground and the excited state vector.

\begin{figure}[hbp]
  \centering
  \includegraphics[height=5cm]{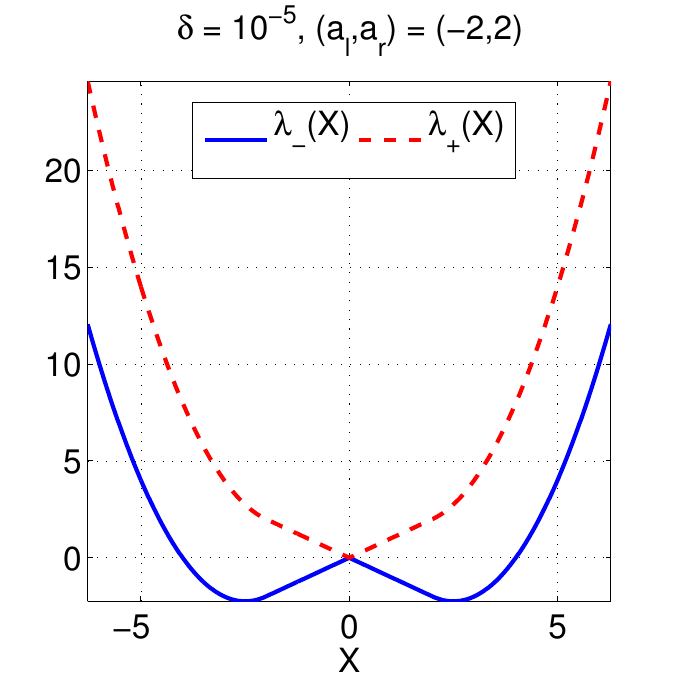}
  \hspace{0.5cm}
  \includegraphics[height=5cm]{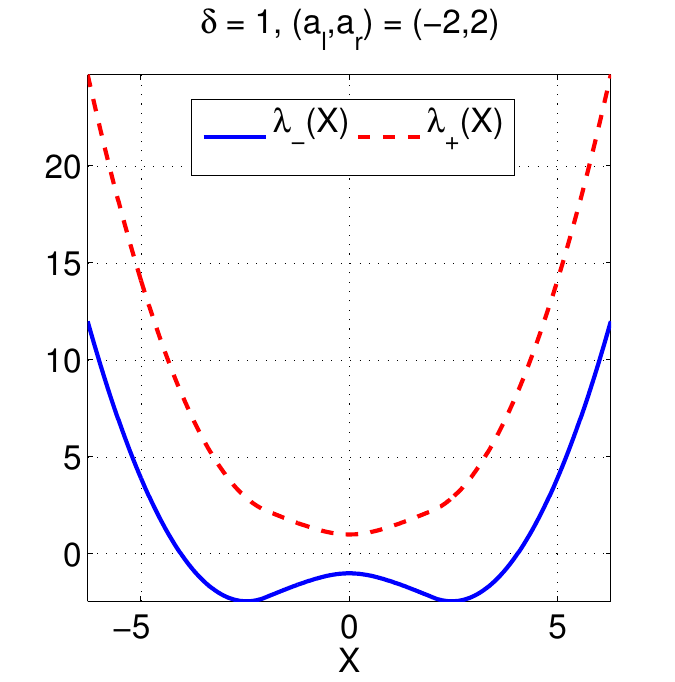}
  \caption{Eigenvalues of the potential matrix~\eqref{eq:potential-1d} for 
     small and large values of  the parameter $\delta$ corresponding 
     to the small and large minimum gaps between the eigenvalues,
     respectively.}
  \label{fig:eigenvalues-1d}
\end{figure}

Choosing~$\delta=0$ gives a conical intersection at $X=0$, and a positive value of~$\delta$
gives a minimum gap~$2\delta$ between~$\lambda_-(X)$ and~$\lambda_+(X)$ at~$X=0$.
A small value of $\delta$ corresponds to a large  probability to be in excited states and a large
value of~$\delta$ corresponds to a small  probability to be in excited states. Figure~\ref{fig:eigenvalues-1d} 
illustrates examples of small and large spectral gaps between~$\lambda_+$ and~$\lambda_-$ 
with two different values of~$\delta$, respectively. 
\subsubsection{The 1D
  discretized Schr\"odinger equation} %
We use a central difference method to solve the one dimensional two-state 
Schr\"odinger equation~\eqref{schrodinger_stat} in the domain $X\in(-2\pi,2\pi)$ with
$(a_l,a_r)=(-2,3)$, mesh-size $h = 4\pi/[10M^{3/4}]$ and partition $X^j = -2\pi+jh, \ j=0,1,2,\ldots, [10M^{3/4}]$.  
The approximation of the eigenvalue problem with the eigenvector components $\Phi_h^j\simeq \Phi(X^j)$ and corresponding eigenvalue $E_h\in \rset$ becomes
\begin{equation}\label{eq:discrete-schrod-1d}
 -\frac{1}{2M}\frac{\Phi_h^{j-1}-2\Phi_h^{j}+\Phi_h^{j+1}}{h^2} +\VOPER(X^j)\Phi_h^j  = E_h\Phi_h^j \COMMA
 \mbox{ for } j=1,2,\dots, [10M^{3/4}]-1,
\end{equation}
using homogeneous Dirichlet boundary conditions $\Phi_h^0=\Phi_h^{10M^{3/4}}=0$.

The approximation of the probability $p_{ex}$ to be in the excited state is given by
\begin{equation}\label{eq:Pr-excited}
 p_{ex} = \frac{ \sum_{j=0}^{[10M^{3/4}] }\langle\Phi_h^j,\Psi_+(X^j)\rangle^2 }{\sum_{j=0}^{[10M^{3/4}] } \langle
 \Phi_h^j,\Phi_h^j\rangle } \, .
\end{equation}
\subsubsection{Approximation of the dynamic transition probability $p_d$}
We show in this computational example that the dynamic transition probability, $p_{d}$, defined in~\eqref{p_E_def}, 
can be obtained from Ehrenfest molecular dynamics simulations~\eqref{eq:ehrenfest} using the following formula
\begin{equation}\label{eq:scattering-probability}
p_d(t)=\left|\frac{\langle\Psi_+(X_t),\psi_t\rangle}{\langle\Psi_+(X_t),\Psi_+(X_t)\rangle}\right|^2\COMMA
\end{equation}
with $t\in\mathbb{R}$ denoting the time. We observe that numerical experiments illustrate 
that the transition probability $p_d(t)$ 
approximates, as time tends to infinity,  the Landau-Zener probability $p_{LZ}$ as defined in~\eqref{eq:p_LZ},
given the Landau-Zener model~\eqref{eq:LZ-model}, see Figure~\ref{fig:md-1d}.

We approximate the solution of the transport equation~\eqref{schrod_first} by the Ehrenfest 
molecular dynamics simulations \eqref{eq:ehrenfest} with $\lambda_0=\lambda_-$ (equal to the smallest eigenvalue of $V$ in\eqref{eq:potential-1d}) in 
the Hamiltonian~\eqref{eq:H_E}. For this computational 
example, we choose $E=1$, $(a_l,a_r)=(-2,2)$ and the initial conditions $t=0$, $X_0 = -4.0$, 
$\psi_0 = \Psi_-(X_0)$, and $P_0 = \sqrt{2\left(E-\lambda_-(X_0)\right)}$.
We approximate the solution of the molecular dynamics using the St\"ormer-Verlet method, \REF{verlet-method}, based on the two symplectic Euler steps, 
\begin{equation}\label{eq:verlet-1d}
\begin{array}{rcl}
 \psi^i_{k+\frac{1}{2}} & = & \psi^i_k-\frac{\Delta t}{2}\sqrt{M}\check\VOPER(X_k)\psi^r_k \COMMA \\
 P_{k+\frac{1}{2}} &=& P_k-\frac{\Delta t}{2}\left(\lambda'_-(X_k) +\left\langle\psi^r_k, \check\VOPER'(X_k)\psi^r_k\right\rangle +  \left\langle \psi^i_{k+\frac{1}{2}}, \check\VOPER'(X_k) \psi^i_{k+\frac{1}{2}}\right\rangle\right) \COMMA\\
 X_{k+1} & = & X_k+ \Delta tP_{k+\frac{1}{2}} \COMMA \\ 
 \psi^r_{k+1} &=& \psi^r_k+\frac{\Delta t}{2}\sqrt{M} \left(\check\VOPER(X_k) + \check\VOPER(X_{k+1})\right) \psi^i_{k+\frac{1}{2}} \COMMA \\
 \psi^i_{k+1} &=& \psi^i_{k+\frac{1}{2}}-\frac{\Delta t}{2}\sqrt{M} \check\VOPER(X_{k+1})\psi^r_{k+1} \COMMA \\
 P_{k+1} & = &  P_{k+\frac{1}{2}} - \frac{\Delta t}{2}\left(\lambda'_-(X_{k+1}) + \left\langle \psi^r_{k+1},\check \VOPER'(X_{k+1})\psi^r_{k+1}\right\rangle + \left\langle\psi^i_{k+\frac{1}{2}},\check \VOPER'(X_{k+1})\psi^i_{k+\frac{1}{2}}\right\rangle\right) .
\end{array}
\end{equation}

\begin{figure}[hbp]
  \centering
  \includegraphics[width=8cm]{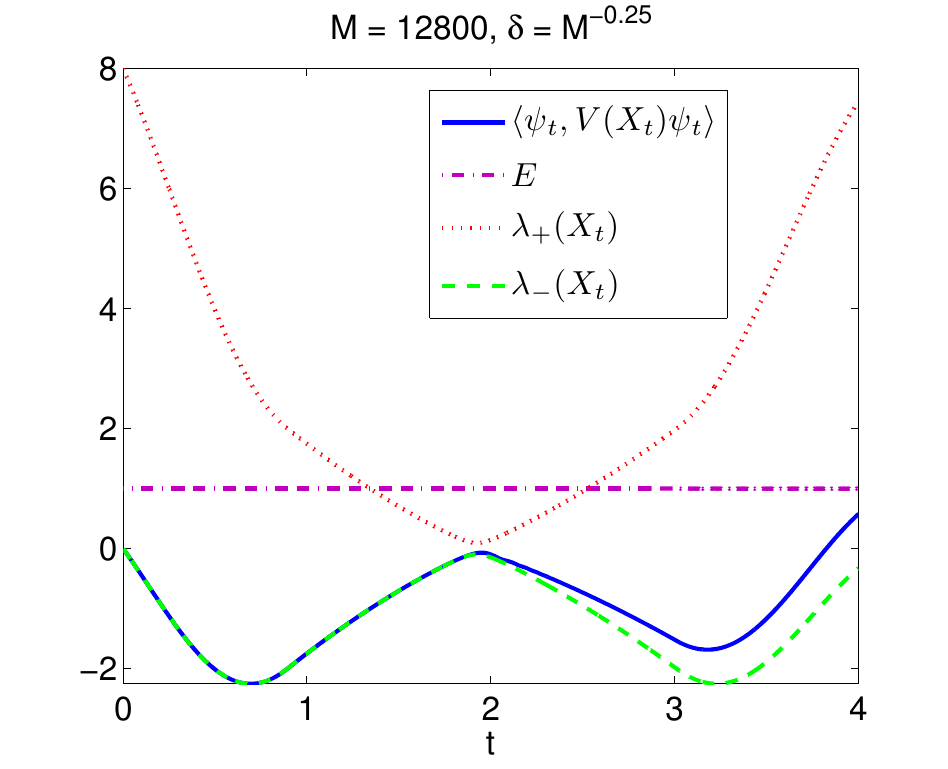}
  \caption{The potential energy function, $t\mapsto \langle \psi_t, V(X_t)\psi_t\rangle$, 
  approximated from the Ehrenfest molecular dynamics  simulations~\eqref{eq:ehrenfest}, 
  deviates away from the eigenvalues after the avoided crossing, if the gap is small and the mass is sufficiently small.} 
  \label{fig:md-1d-pot}
\end{figure}
\begin{figure}[hbp]
 \centering
 \subfigure{\label{md-1D:subfig1}
   \includegraphics[width=0.32\textwidth,height=0.32\textwidth]{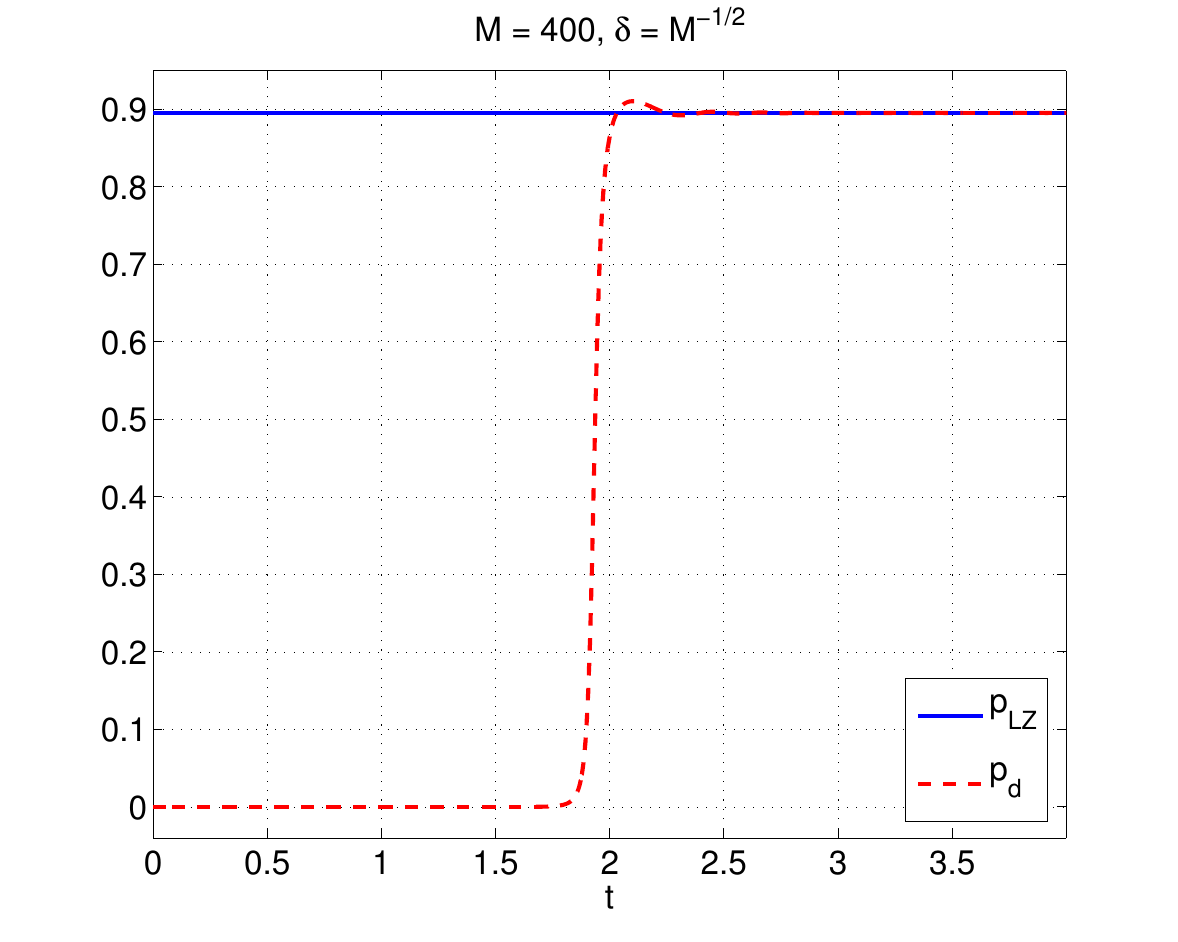}
 }
 \qquad\quad
 \subfigure{\label{md-1D:subfig2}
   \includegraphics[width=0.32\textwidth,height=0.32\textwidth]{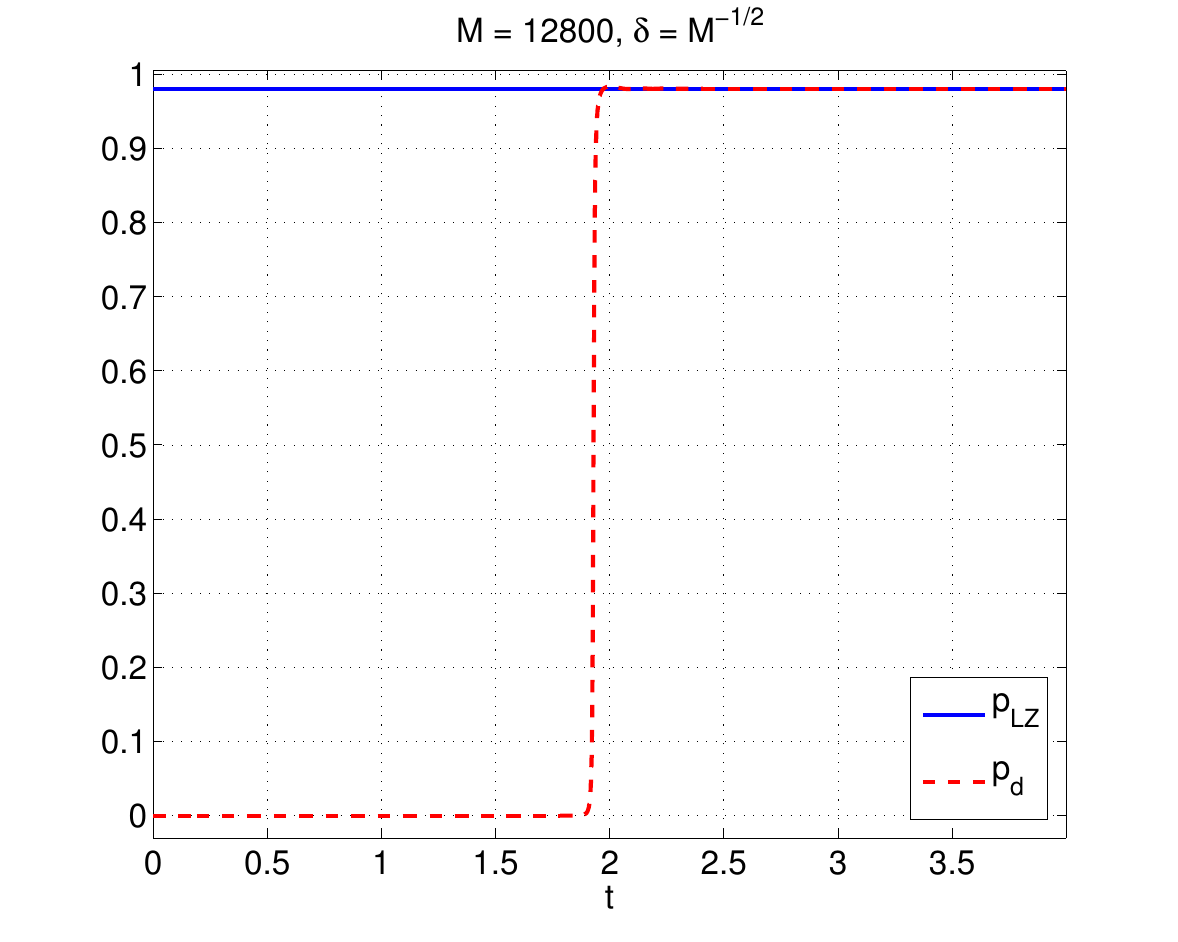}
 }

 \subfigure{\label{md-1D:subfig3}
   \includegraphics[width=0.32\textwidth,height=0.32\textwidth]{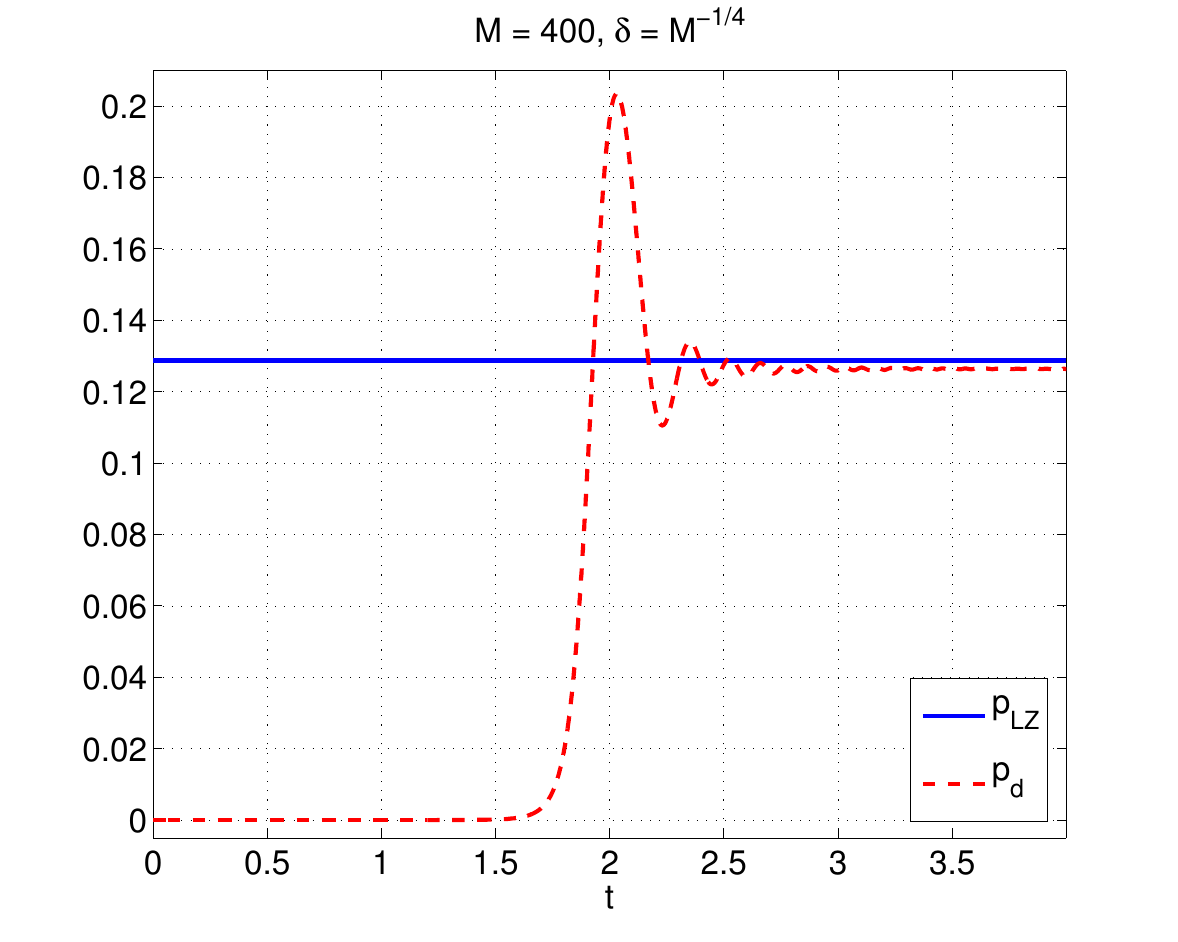}
 }
 \qquad\quad
 \subfigure{\label{md-1D:subfig4}
   \includegraphics[width=0.32\textwidth,height=0.32\textwidth]{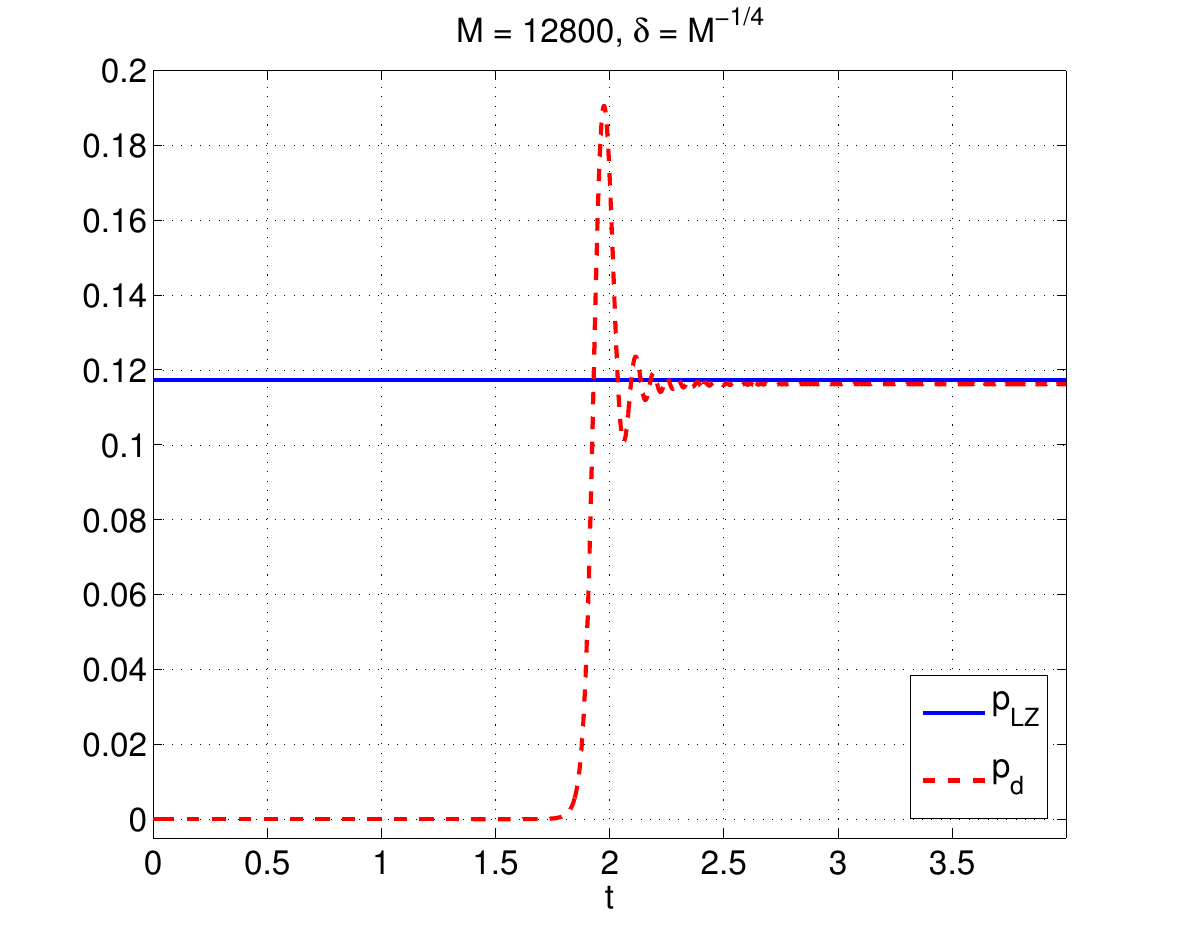}
 }
 \subfigure{\label{md-1D:subfig5}
   \includegraphics[width=0.3\textwidth,height=0.3\textwidth]{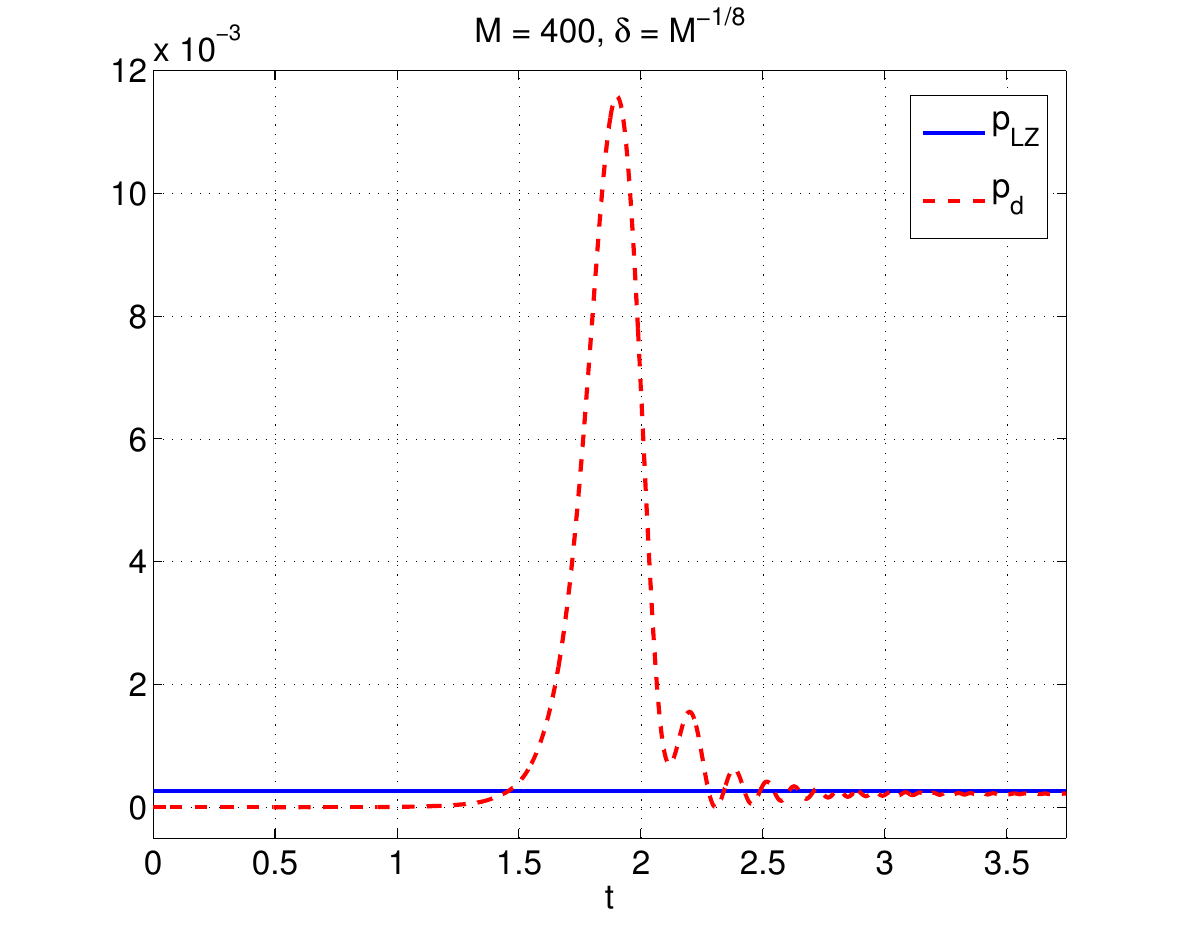}
 }
 \qquad\quad
 \subfigure{\label{md-1D:subfig6}
   \includegraphics[width=0.3\textwidth,height=0.3\textwidth]{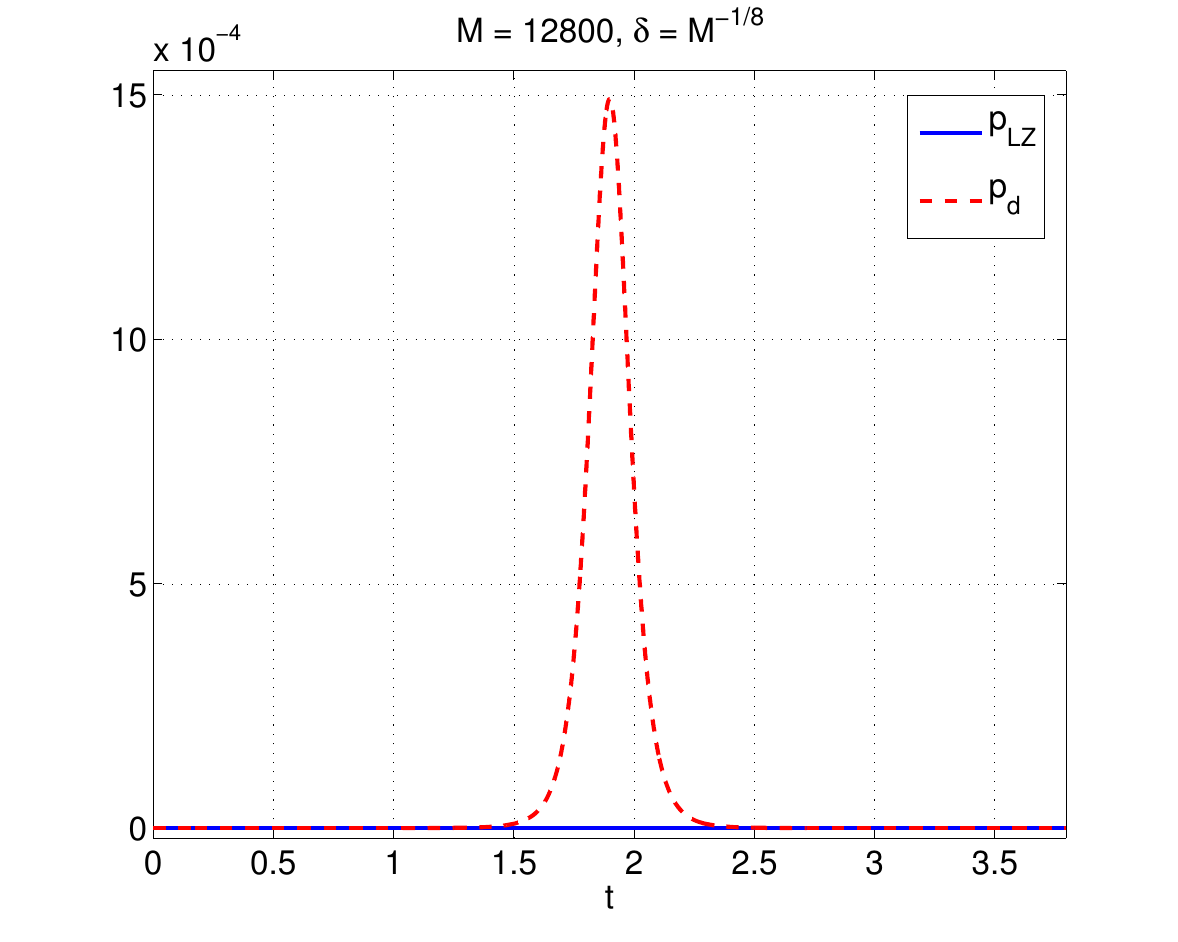}
 } 
  \caption{Estimation of the Landau-Zener probability, $p_{LZ}$, from the Ehrenfest
  molecular dynamics simulations~\eqref{eq:ehrenfest}. The estimation,~$p_d(t)$, 
  overshoots near the avoided crossing and eventually stabilizes around~$p_{LZ}$.
  The overshoot period and the relative magnitude are more prominent for the cases of smaller
  Landau-Zener probabilities.}
  \label{fig:md-1d}
\end{figure}
\subsubsection*{Conclusions.}
Figure~\ref{fig:md-1d-pot} illustrates that the potential energy function $t\mapsto \langle \psi_t, V(X_t)\psi_t\rangle$,
obtained from the Ehrenfest molecular dynamics simulations is equal to  the ground 
state energy,~$\lambda_-(X_t)$, until the dynamics reaches the avoided conical 
intersection. Then it deviates from the ground state energy and continues in 
between the ground and excited state energies. 

Figure~\ref{fig:md-1d} shows that the transition probability, $p_d(t)$, obtained 
using the formula~\eqref{eq:scattering-probability} based on Ehrenfest molecular 
dynamics simulations, remains zero until it reaches the avoided conical intersection and 
then it oscillates near the avoided conical intersection region and approaches~$p_{LZ}$ 
asymptotically as~$t\rightarrow\infty$.

\subsection{Model 2: A two dimensional problem}\label{sec:2d-example}
\begin{figure}[htbp]
  \centering
  \includegraphics[width=0.4\textwidth]{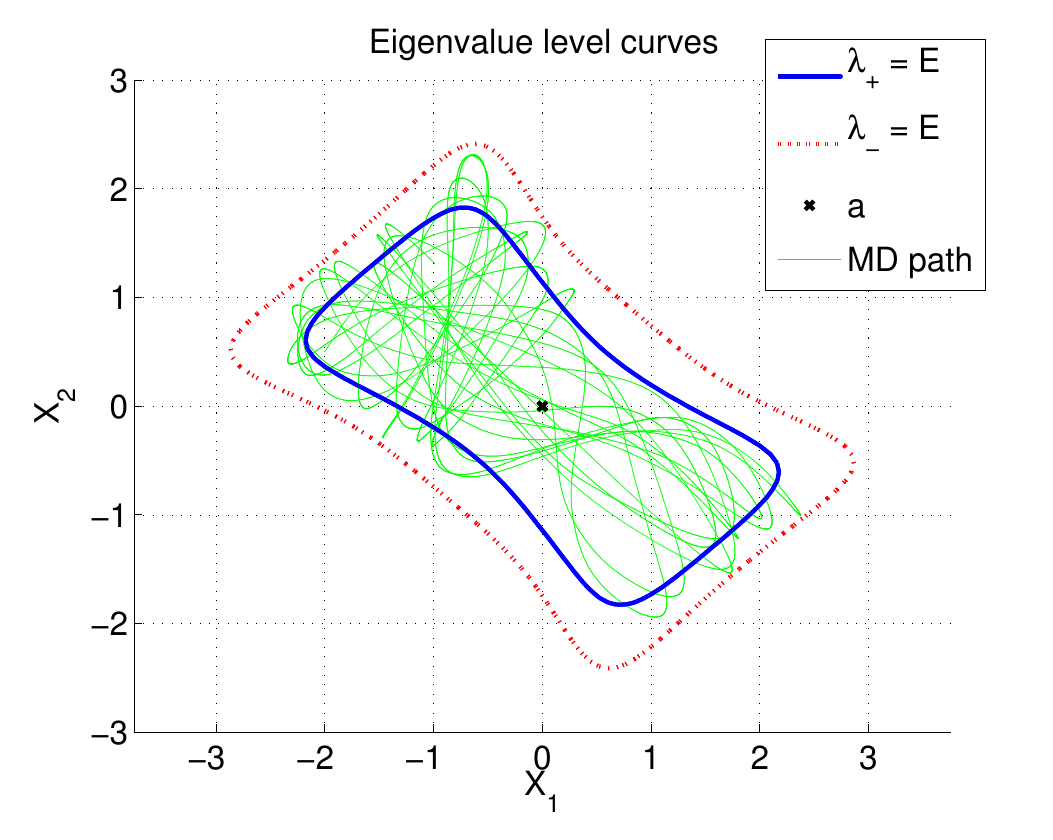}
  \includegraphics[width=0.4\textwidth]{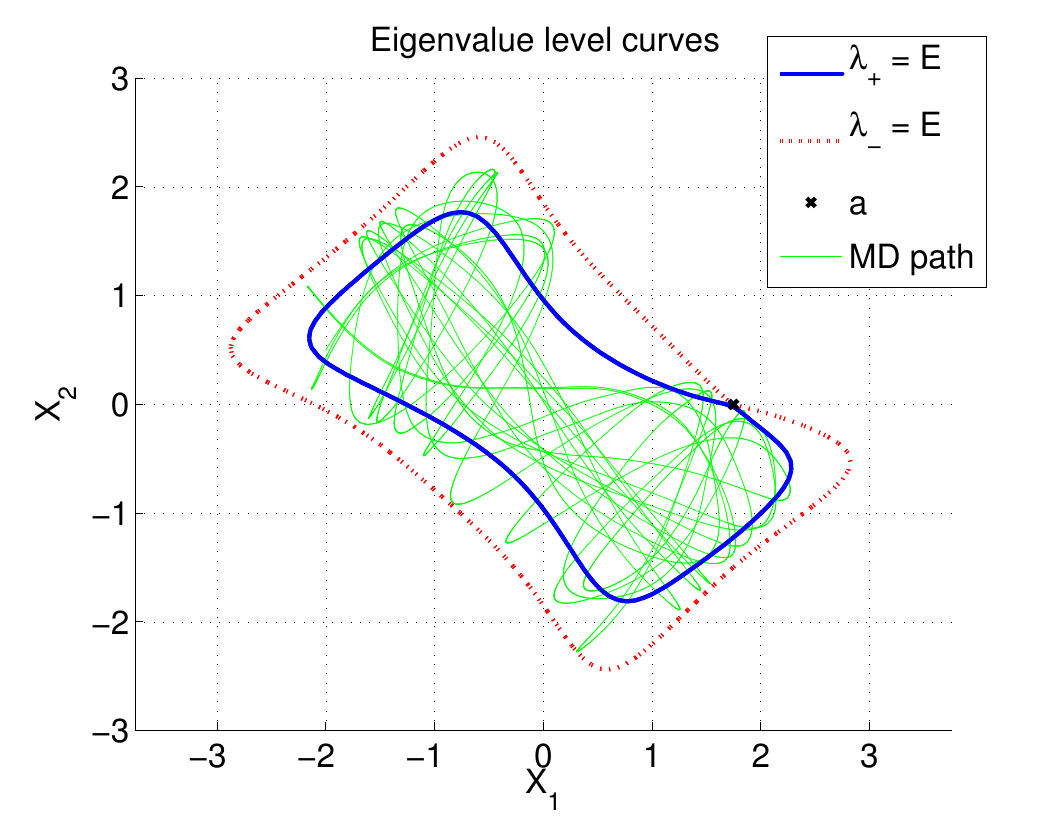} \\
  \includegraphics[width=0.4\textwidth]{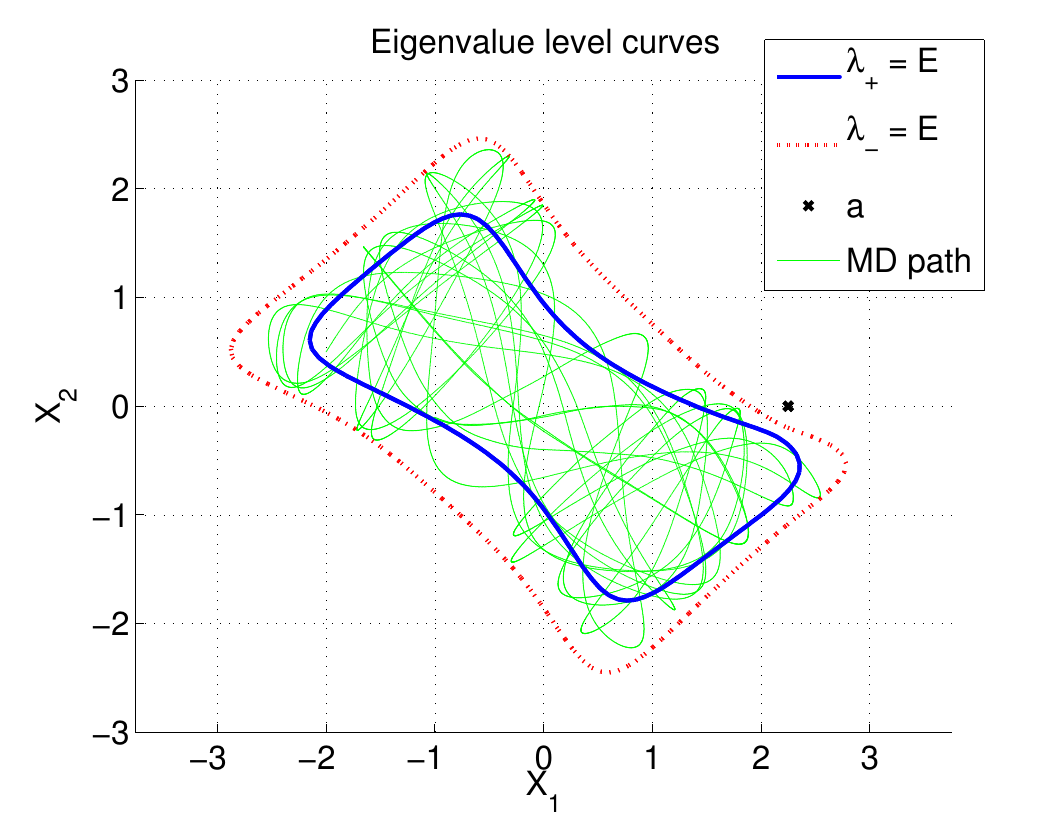}
  \includegraphics[width=0.4\textwidth]{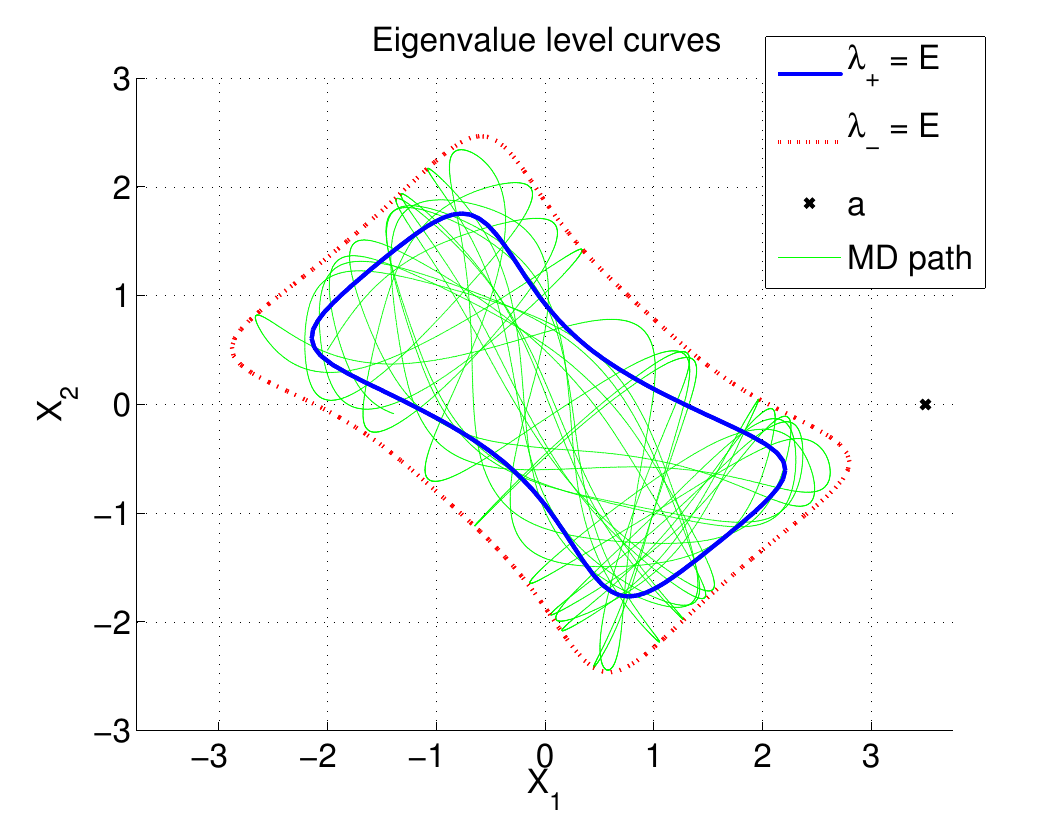}
  \caption{Level curves of the eigenvalues $\lambda_\pm$ of the potential~\eqref{eq:2d-potential}
  at the energy $E=1.5$ with conical intersection points $a=(0,0),(1.75,0),(2.25,0),$ and $(3.5,0),$ respectively.
  We illustrate examples of the Ehrenfest molecular dynamics~\eqref{eq:ehrenfest} paths computed
  using St\"ormer-Verlet method \eqref{eq:verlet-1d} with mass $M=100,$ time $t\in[0,100],$
  time steps $\Delta t = 1/(32\sqrt M),$ and initial data $X_0=[-2,0.5], P_0 = [1, \sqrt{2(E-\lambda_-(X_0))-1}],$ 
  and $\psi_0=\Psi_-(X_0)$.}   
  
  \label{fig:eig-level-curves-2d}
\end{figure}

We  consider the two dimensional, time-independent Schr\"odinger equation~\eqref{schrodinger_stat} 
with the heavy-particle coordinate~$X = (X_1,X_2) \in \mathbb R^2$ and 
two electron states $J=2$. The Hamiltonian is
$$\HOPER(X) = V(X) -\frac{1}{2}M^{-1}\Delta_X \COMMA$$
with
\begin{equation}\label{eq:2d-potential}
\VOPER(X) := \Big(\underbrace{\frac{1}{2}(X_1^2+\alpha X_2^2) + \beta\sin(X_1X_2)}_{=: \lambda_s(X)}\Big)I + \eta\left[\begin{array}{cc}
                               v_{1}(X) & v_{2}(X)\\
                               v_{2}(X) & -v_{1}(X)
                        \end{array}\right],
\end{equation}
where $I$ is the $2 \times 2$ identity matrix, the functions $v_{1},v_{2}:\mathbb{R}^2\rightarrow\mathbb{R}$ 
are given by (A) or (B) below and $\alpha=\sqrt 2$, $\beta = 2$ and $\eta = 1/2$. 
As in Example 1, we have for each $X$ the eigenvalue problem, $\VOPER(X)\Psi_\pm(X) = \lambda_\pm(X)\Psi_\pm(X)$, 
with the eigenvalues $\lambda_\pm(X) = \lambda_s(X) \pm \eta\sqrt{(v_1(X))^2 + (v_2(X))^2}$ and
the ground state and excited state eigenvectors $\Psi_-(X)$ and $\Psi_+(X)$, respectively.
We choose the energy, $E=1.5$.
In this example we study the following two cases of interactions between the two 
potential surfaces~$\lambda_+$ and~$\lambda_-$:
\begin{enumerate}
 \item[(A)] {\it A line intersection.}
 We choose $v_1(X) = \arctan(X_1/\eta)$ and $v_2(X) = \delta/\eta$, 
 where $\delta$ is a non-negative constant that defines the minimum distance 
 between the potential surfaces. In this example, the potential surfaces intersect 
 each other (for $\delta = 0$) or have the minimum distance between each other 
 at the line $X_1=0$.
 Here, choosing a small value for $\delta$ corresponds to a large  
 probability to be in excited states and a large value for $\delta$ will give a small  probability to be in excited states. 
 The parameter $\delta$ in this two dimensional example is analogous to the 
 parameter~$\delta$ in the one dimensional example given in Section~\ref{sec:1d-example}.
 \item[(B)] {\it A conical intersection.}
 We choose $v_1(X)=\arctan((X_1-a_1)/\eta)$ and $v_2(X)=\arctan((X_2-a_2)/\eta)$ 
 with $a = (a_1,a_2)\in\mathbb{R}^2$ being a chosen point in the two dimensional space. 
 If $a$ is chosen such that $\lambda_\pm(a)$ are smaller than the energy, $E$, 
 we have a conical intersection between the potential surfaces at the point~$a$ 
 in the { classically allowed} region $R := \{X: \lambda_-(X) \le E\}$, otherwise 
 we have a positive gap between the potential surfaces in the domain $R$. 
 Here, choosing $a$ in the origin gives a larger absolute momentum, $|P|$, at 
 the conical intersection whereas choosing $a$ far from the origin will yield
 a smaller $|P|$, which yield larger and smaller probabilities to be in excited states, 
 respectively. Figure~\ref{fig:eig-level-curves-2d} shows level curves of the 
 eigenvalues, along with example molecular dynamics paths, for varying conical intersection
 points.
\end{enumerate}
\subsubsection{The 2D discretized Schr\"odinger equation}
We use a standard finite difference method to discretize the two dimensional Schr\"odinger 
eigenvalue problem~\eqref{schrodinger_stat} with mesh-size $h$  
in both the $X_1$ and $X_2$ directions in the computational domain $\Omega = [-4,4]^2$. 
The unknown eigenvalue components $ \Phi_h^{j,k}\simeq \Phi(X^{j,k})$,  with the nodal point
 $X^{j,k}=(-4+jh, -4+kh)$  for $j,k=1,2,\dots, 8/h-1$, solve the discrete
eigenvalue problem
\begin{equation}\label{eq:discrete-schrod-2d}
-\frac{1}{2M}\left(\frac{\Phi_h^{j-1,k}-2\Phi_h^{j,k}+\Phi_h^{j+1,k}}{h^2} + \frac{\Phi_h^{j,k-1}-2\Phi_h^{j,k}+\Phi_h^{j,k+1}}{h^2} \right) + \VOPER(X^{j,k})\Phi_h^{j,k}= E_h\Phi^{j,k}_h\, %
\end{equation}
with homogeneous Dirichlet boundary conditions  and $h = 1/(4\sqrt{M})$. 
We use the solution of the discrete Schr\"odinger eigenvalue problem to determine the approximate 
Schr\"odinger observables
\begin{equation}\label{eq:schrod-obs}
  g_{\SCH} = \frac{\sum_{j,k} g(X^{j,k})\rho(X^{j,k})}{\sum_{j,k}\rho(X^{j,k})}\COMMA
\end{equation}
with $g(X):\mathbb R^2 \rightarrow \mathbb R$ and $\rho(X)=|\Phi_h(X)|^2$.

The value of the molecular dynamics observable, $g_{\MD}$, %
is given by 
\begin{equation*}
  g_{\MD} := \lim_{\delta \rightarrow 0+} \frac{\int_{E<H_0(X,P)<E+\delta} g(X,P)dXdP}{\int_{E<H_0(X,P)<E+\delta} dXdP}\COMMA
\end{equation*}
which, for the two dimensional case, can be written\footnote{For a fixed $X\in\rset^{d}$, $P-$integration using spherical coordinates in the shell $E<|P|^2/2+\lambda_0(X)<E+\delta$ yields the $P-$ dependence 
$[|P|^{d}]_{\sqrt{E-\lambda_0(X)}}^{\sqrt{E-\lambda_0(X)+\delta}}$ which differentiated with respect to $\delta$ gives
\eqref{eq:g-MD} for $d=2$ and the additional factor $( E-\lambda_0(X))^{d/2-1}$ for $d\ne 2$. } as 
\begin{equation}\label{eq:g-MD}
  g_{\MD} = \frac{\int_{H_0(X,0)\le E} g(X)dX}{\int_{H_0(X,0)\le E} dX}\COMMA
\end{equation}
when $g$ depends only on $X$.

\subsubsection{Verification of the ergodic rate condition \eqref{g_bar}.}
\begin{figure}[htbp]
\centering
  \includegraphics[width=0.4\textwidth,height=0.4\textwidth]{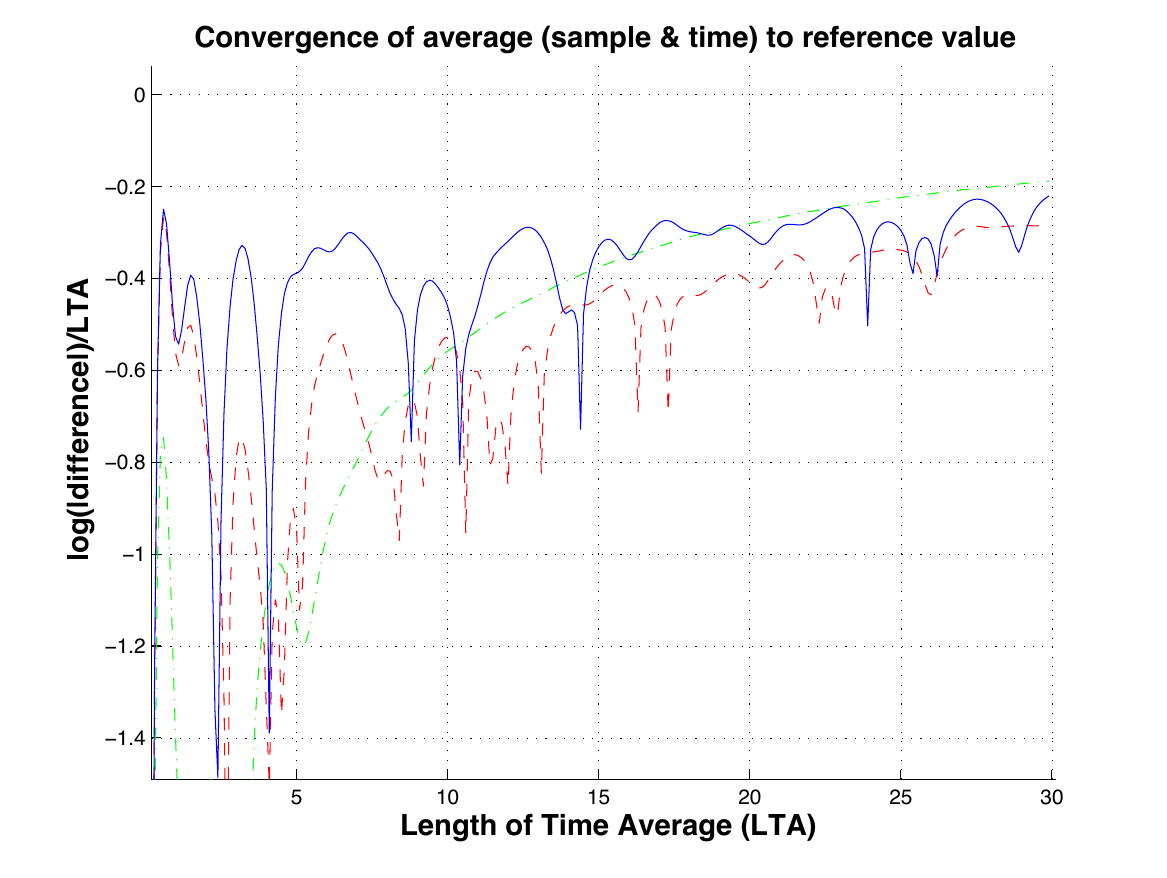}
   \caption{
The figure shows the function 
$
T \mapsto
T^{-1} \log \Big(T^{-1}\int_{T}^{2T}\big(\bar{\mathbb E}[g\circ S_{t0}(Z_0)] -g_{MD}(Z_0) \big)dt\Big)\, ,$
which approximates the
decay rate $2\gamma$ in assumption \eqref{g_bar},   for the projected dynamics \eqref{sde_proj} with
the potential $\lambda(X)=X_1^2/2+X_2^2/\sqrt 2+ 2\sin(X_1X_2)$ and energy $E=H(Z_0)=1.5$.
The empirical mean is 
 $\bar{\mathbb E}[g\circ S_{t0}(Z_0)] := \sum_{n=1}^N \sin(\bar X_1(t;n)\bar X_2(t;n))/N$ and the number of stochastic paths are $N=2\times 10^5$. The initial point is $Z_0=(0,0, \sqrt{1.5},\sqrt{1.5})$ and $(\bar X(t;n),\bar P(t;n))$ is the phase space point for sample path $n$ at time $t$  determined by the numerical method \eqref{proj_nm} for time steps $\Delta t=0.01$. The blue solid curve is for diffusion coefficient
  $\epsilon=10^{-4}$, the red dashed for $\epsilon=10^{-2}$ and the green dashed-dot for $\epsilon=0.25$.
Here the approximation $g_{MD}(Z_0)=-0.4388$  is determined by an accurate Monte Carlo simulation in the microcanonical ensemble.
 }
\label{fig313}
\end{figure}

Figure \ref{fig313}  determines the convergence rate $\gamma$ in assumption \eqref{g_bar} 
based on numerical approximation of the ergodic projected stochastic dynamics \eqref{sde_proj}. The   numerical method is given by the steps $\bar Z(n\Delta t)\mapsto \bar Z\big((n+1)\Delta t\big)$ satisfying
 \begin{equation}\label{proj_nm}
 \begin{split}
 (X,P)&=\bar Z(n\Delta t)\\
 X_*&=X+P\Delta t \\
 P_*&=P-\nabla\lambda(X_*)\Delta t \\
 Z^1&=(X_*,P_*)\\
 Z^2&=Z^1+\sqrt\epsilon\mathbb P(Z^1)\Delta W\\
 Z^3&=Z^2+\mu \nabla H(Z^2)/|\nabla H(Z^2)|\\
 \bar Z\big((n+1)\Delta t\big)&=Z^3\, ,\\
 \end{split}
  \end{equation}
  where $\mu\in \rset$ is chosen so that $H(Z^2+\mu \nabla H(Z^2)/|\nabla H(Z^2)|)=E$.
  Here $\Delta W=\sqrt{\Delta t}(\xi_1,\xi_2,\xi_3,\xi_4)$ and $\xi_i$ are all independent taking values $\pm 1$ with probability $1/2$. The temperature parameter is $\tau=1/2$. This scheme with phase-space path 
  $\{\bar Z(n\Delta t)=\big(\bar X(n\Delta t),\bar P(n\Delta t)\big)\ | \ n=0,1,2\ldots \}$ satisfies the convergence assumptions in \cite{tony_faou}\footnote{Related to examples (4.15) and (4.7) in \cite{tony_faou}.} so that $\mathbb E[g\big(\bar X(T),\bar P(T)\big)- g\big(X(T),P(T)\big)]=\mathcal O(\Delta t)$, uniformly in $T$.
Figure \ref{fig313} shows that in this case $2\gamma>0.2$. The maximal Lyapunov constant is determined numerically to $\hat C=0.35$ and the convergence rate exponent $\gamma/(\hat C+\gamma)$ in Theorem \ref{thm:potential} becomes $0.1/0.45\approx 0.2$.

\subsubsection{Convergence of observables in 2D}\label{sec:conv-obs-2D}
\begin{figure}[htbp]
 \centering
 \subfigure{\label{obs-schrod-spin-2D-M200:subfig1}
  \includegraphics[width=0.35\textwidth,height=0.35\textwidth]{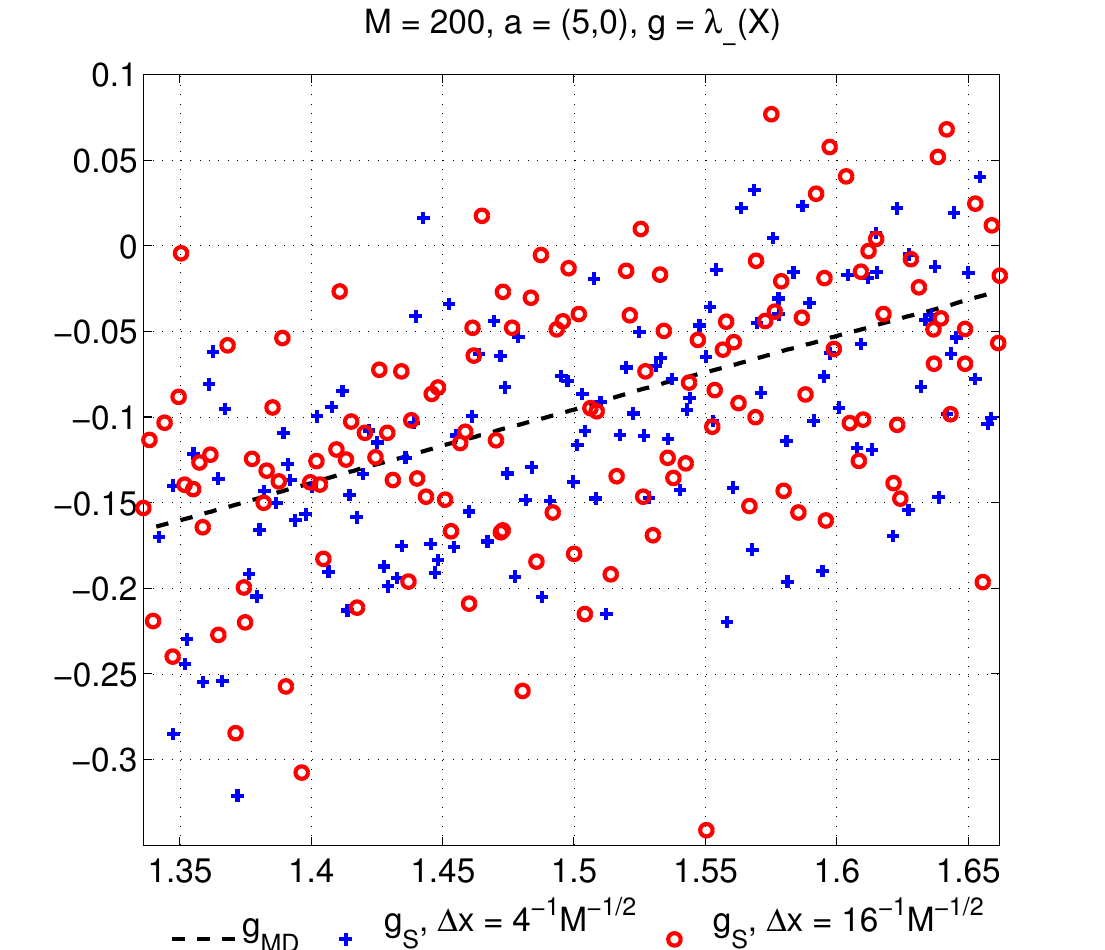}
}
\qquad
\subfigure{\label{obs-schrod-spin-2D-M200:subfig2}
  \includegraphics[width=0.35\textwidth,height=0.35\textwidth]{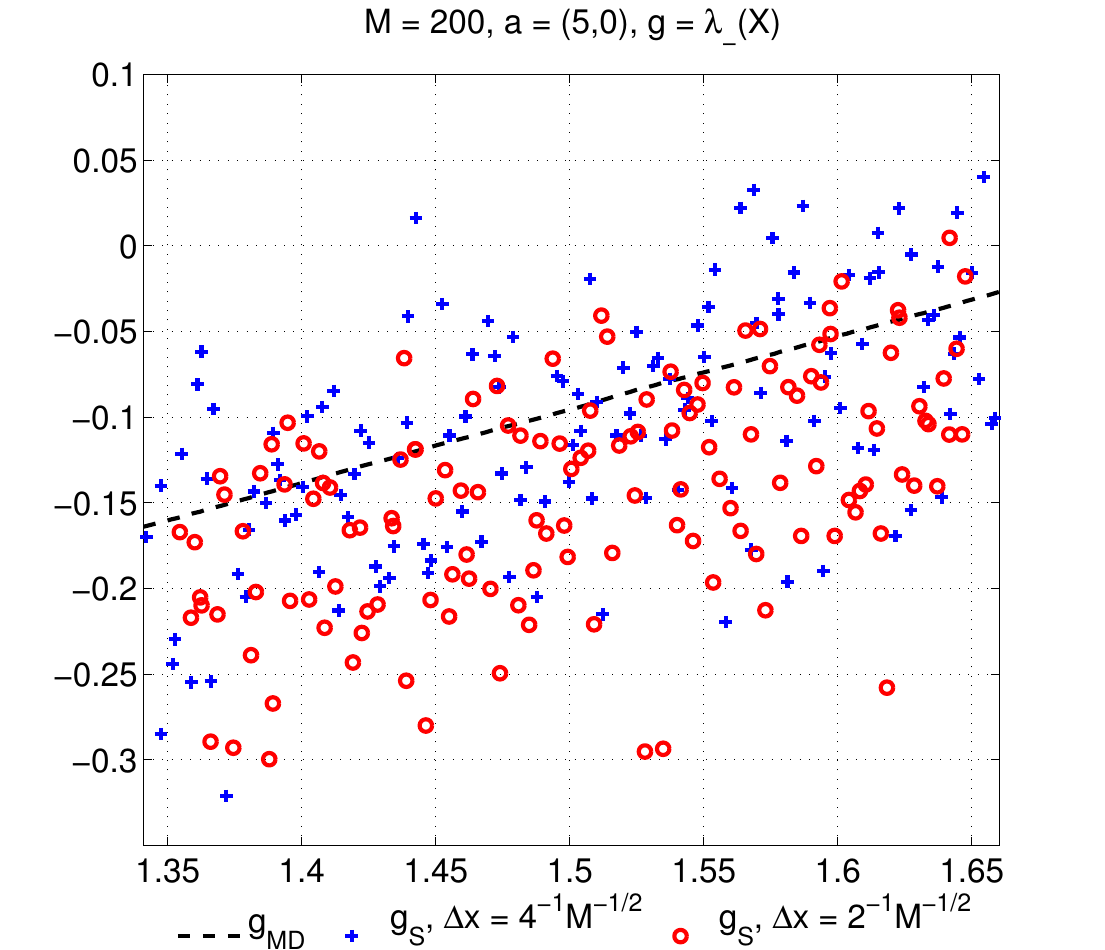}
}

 \subfigure{\label{obs-schrod-spin-2D-M200:subfig3}
  \includegraphics[width=0.35\textwidth,height=0.35\textwidth]{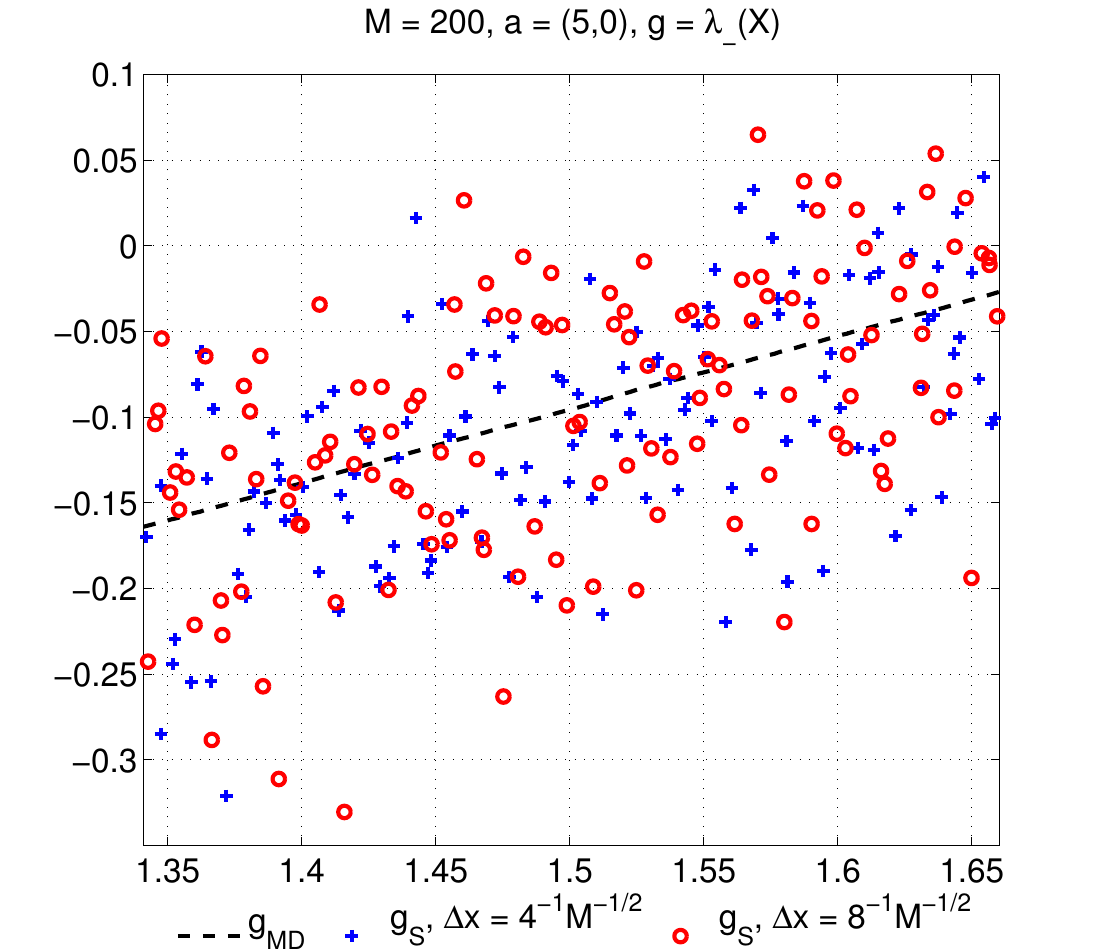}
 }
 \qquad
 \subfigure{\label{obs-schrod-spin-2D-M200:subfig4}
  \includegraphics[width=0.35\textwidth,height=0.35\textwidth]{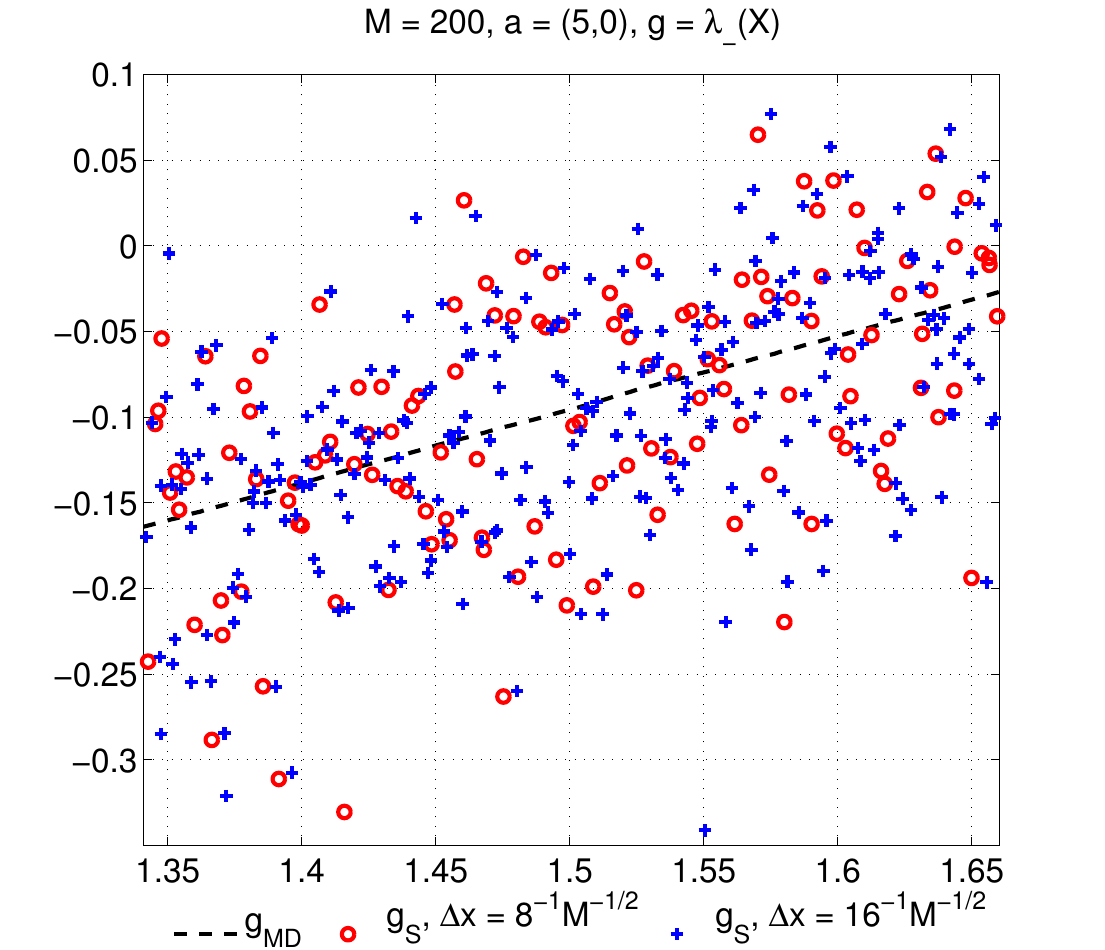}
 }
  \caption{Schr\"odinger observables,~$g_{\SCH}$, as a function of Schr\"odinger 
  eigenvalues, $E_h$, for mass $M=200$, and the conical intersection case with~$a=(5,0)$ 
  outside the {classically allowed} region, compared with the corresponding molecular 
  dynamics observables,~$g_{\MD}$, which we compute using Monte Carlo integrations 
  based on the formula~\eqref{eq:g-MD}. The plots show that the solutions obtained with
  the mesh size $h = 1/(4\sqrt{M})$, $h = 1/(8\sqrt{M})$ and $h = 1/(16\sqrt{M})$ are comparable, 
  whereas the solution obtained with the mesh size $h = 1/(2\sqrt{M})$ appears less accurate.}
\label{fig:obs-schrod-spin-2D-M200}
\end{figure}

\begin{figure}[htbp]
 \centering
 \subfigure{\label{obs-schrod-spin-2D-M3200:subfig1}
  \includegraphics[width=0.35\textwidth,height=0.35\textwidth]{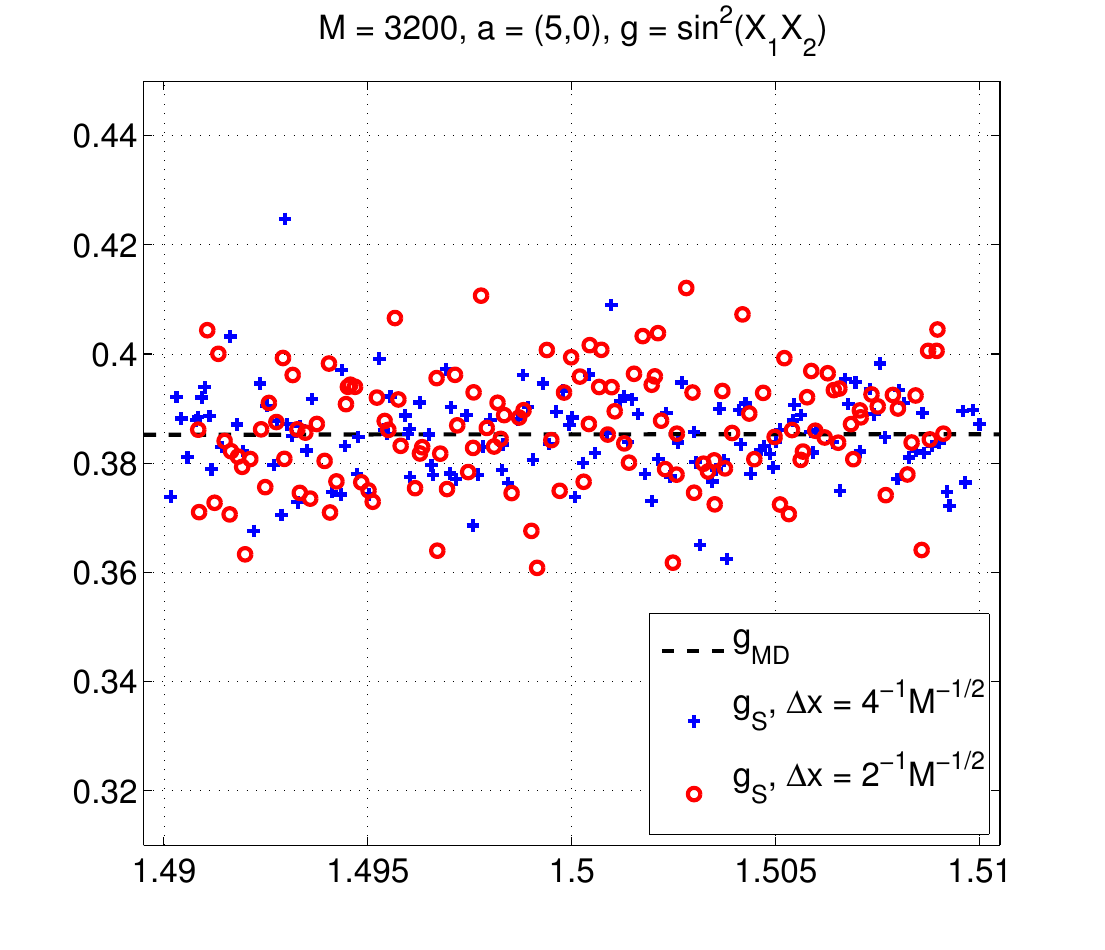}
}
\qquad
\subfigure{\label{obs-schrod-spin-2D-M3200:subfig2}
  \includegraphics[width=0.35\textwidth,height=0.35\textwidth]{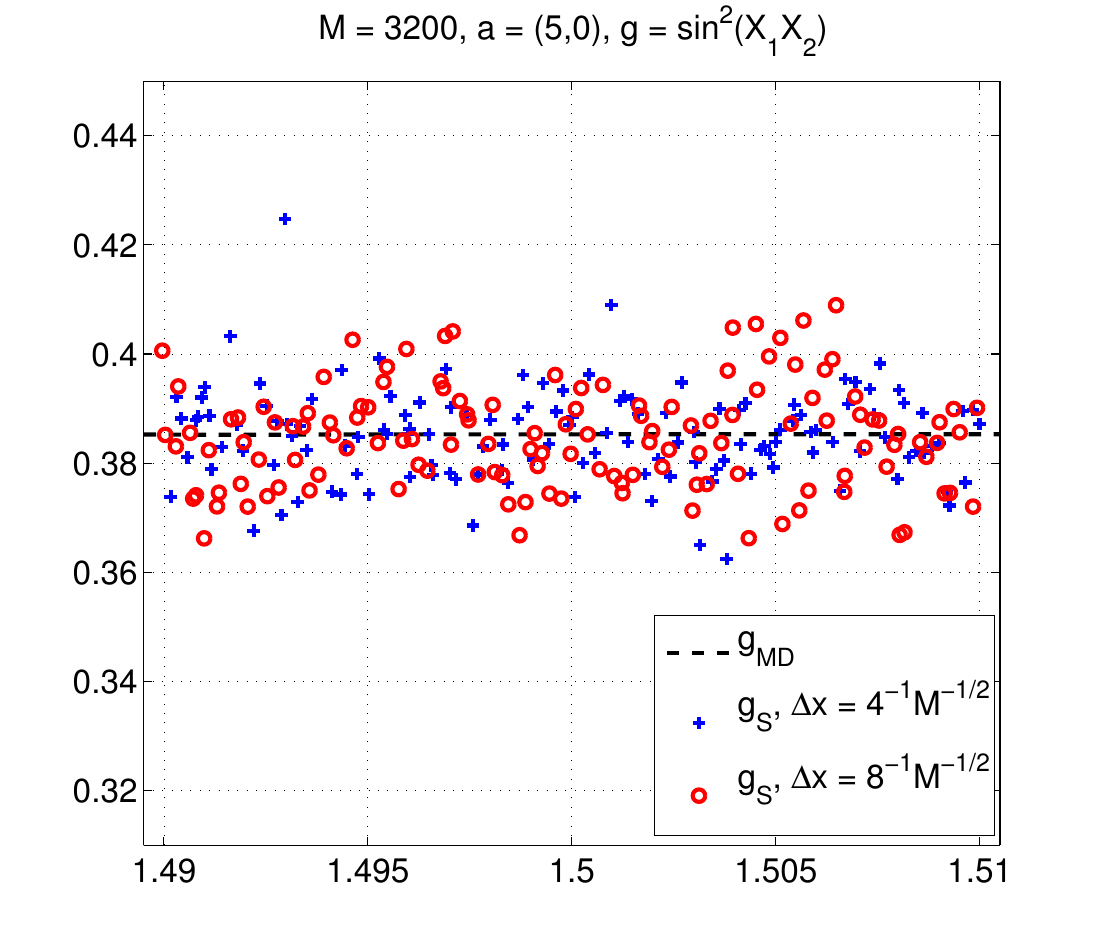}
}

 \subfigure{\label{obs-schrod-spin-2D-M3200:subfig3}
  \includegraphics[width=0.35\textwidth,height=0.35\textwidth]{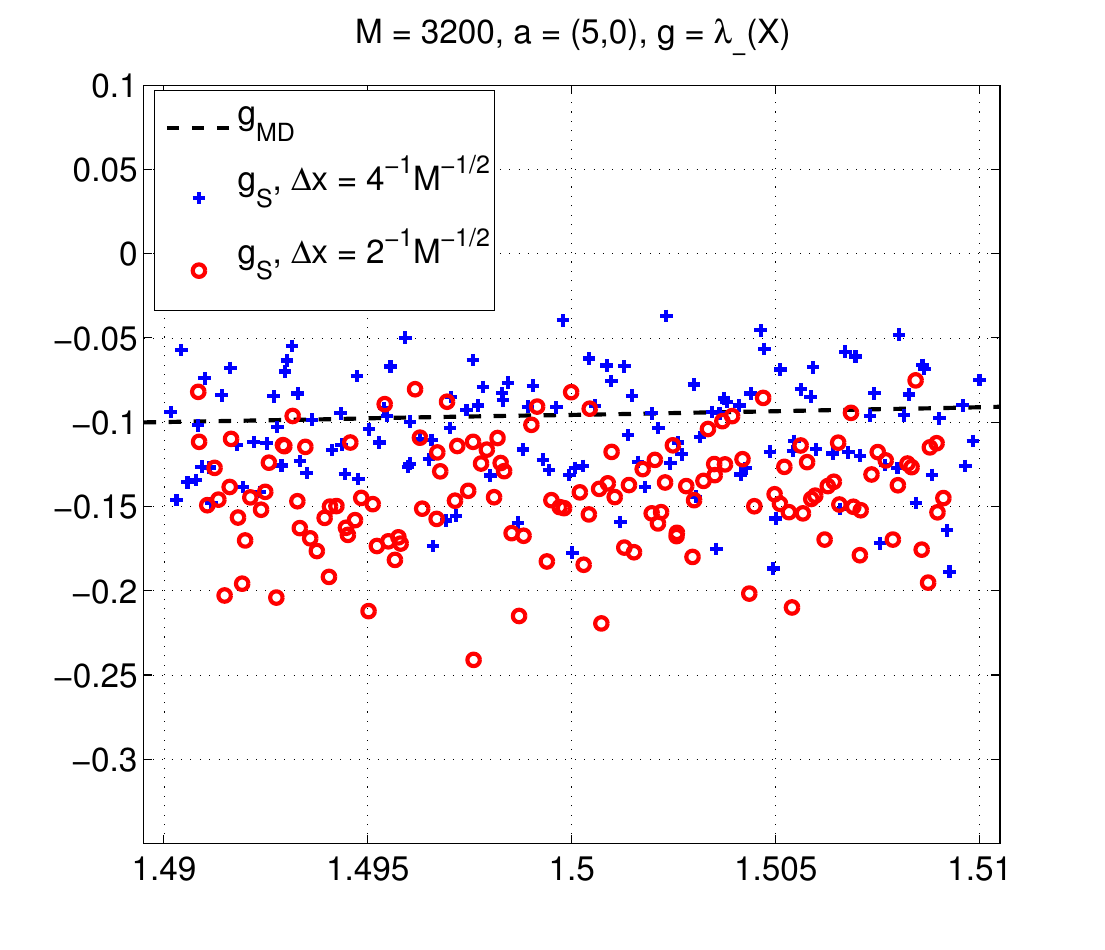}
 }
 \qquad
 \subfigure{\label{obs-schrod-spin-2D-M3200:subfig4}
  \includegraphics[width=0.35\textwidth,height=0.35\textwidth]{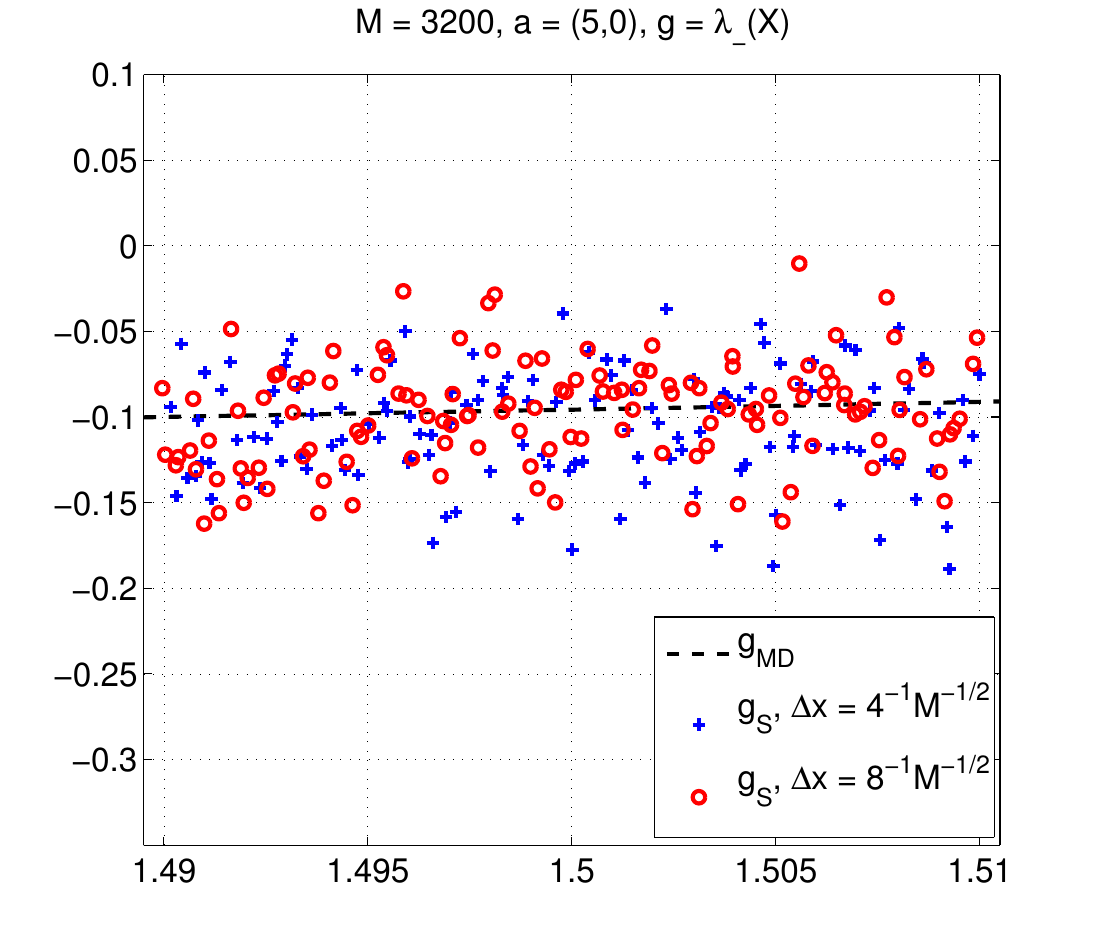}
 }
  \caption{Schr\"odinger observables,~$g_{\SCH}$, as a function of Schr\"odinger 
  eigenvalues, $E_h$, for mass $M=3200$, and the conical intersection case with~$a=(5,0)$ 
  outside the {classically allowed} region, compared with the corresponding molecular 
  dynamics observables,~$g_{\MD}$, which we compute using Monte Carlo integrations 
  based on the formula~\eqref{eq:g-MD}. The plots show that the solutions obtained with
  the mesh size $h = 1/(4\sqrt{M})$ and $h = 1/(8\sqrt{M})$ are comparable, whereas the solutions
  obtained with the mesh size $h = 1/(2\sqrt{M})$ appear less accurate.}
\label{fig:obs-schrod-spin-2D-M3200}
\end{figure}

\begin{figure}[htbp]
 \centering
\includegraphics[width=0.7\textwidth]{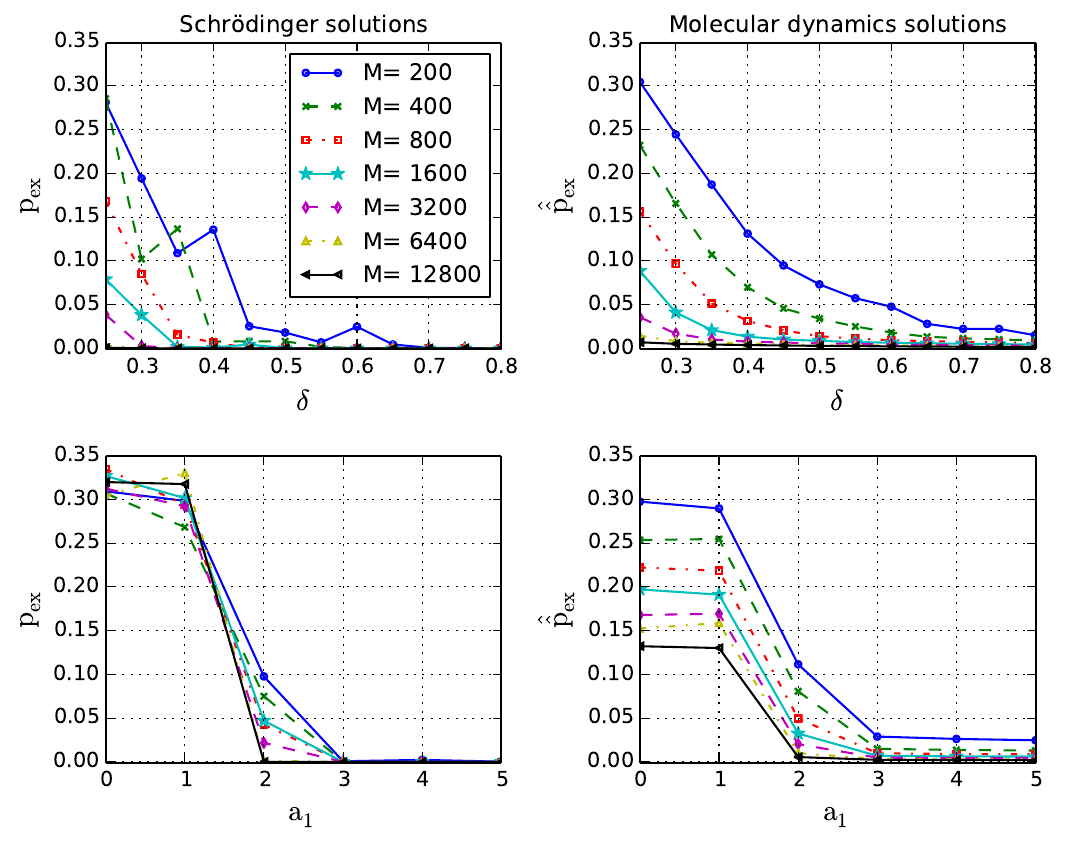}
  \caption{Plots showing $p_{ex}$ estimated from the
      solution of the discrete Schr\"odinger eigenvalue problems (left
      column) and the corresponding $\hat{\hat{p}}_{ex}$ estimated by the
      molecular dynamics Algorithm~\ref{alg:compute-pe-md} (right
      column) for the line intersection cases (first row) and conical
      intersection cases (second row). The Schr\"odinger solution is
    computed by averaging the solutions for $20$ eigenvalues around
    $E=1.5$ for the line intersection and by averaging the solutions
    for $500$ eigenvalues around $E=1.5$ for the conical
    intersection. }
\label{fig:pe-Mpow0p5}
\end{figure}

Figures~\ref{fig:obs-schrod-spin-2D-M200} and~\ref{fig:obs-schrod-spin-2D-M3200} 
show that the Schr\"odinger observables,~$g_{\SCH}$,
are close to the ergodic molecular dynamics observables,~$g_{\MD}$.
We can compare the numerical results in Figures~\ref{fig:obs-schrod-spin-2D-M200} and~\ref{fig:obs-schrod-spin-2D-M3200}
with Theorem~\ref{ergod_sats} for the case when the probabilities 
to be in the excited state,~$p_{ex}$,
are small, as in Figure~\ref{fig:pe-Mpow0p5} for
$a=(5,0)$ in the conical intersection case. We see in Figures~\ref{fig:obs-schrod-spin-2D-M200} and~\ref{fig:obs-schrod-spin-2D-M3200}
that the standard deviation of $|g_{\SCH}-g_{\MD}|$ for $M=3200$ is roughly $1/2$ times of the 
standard deviation for $M=200$, i.e.~the numerical results agree roughly with an $\mathcal O(M^{-1/4})$ estimate in
the theorem, see~Table~\ref{tab:obs-conv}.

\begin{table}
\centering
\begin{tabular}{ c | c | c }
& $g=\sin^2(X_1X_2)$ & $g=\lambda_-(X)$ \\ \hline
$e_\infty$ &   -0.3079 &  -0.3159 \\
$e_W$ &   -0.3052 &  -0.2506 \\
$e_1$ &   -0.2744 &  -0.2950 \\
$e_2$ &   -0.2671 &  -0.2890 \\
\end{tabular}
\caption{Numerical results for the order of convergence of the observables in 2D as explained in Section~\ref{sec:conv-obs-2D}.
The table shows $\varrho$ where the order of convergence is $\mathcal O(M^\varrho)$. The convergence
rate is computed as $\varrho = (\log e({M_1})-\log e({M_2}))/(\log M_1-\log M_2),$ with $e({M_i})$ being the norm of error corresponding to the
mass $M_i,$ where $M_1=200,$ and $M_2=3200$. We use the norms $e_\infty = \max_i |(g_{\SCH}-g_{\MD})(E_i)|,$ $e_W = |\sum_{i} (g_{\SCH}-g_{\MD})(E_i)|/N_E,$
$e_1 = \sum_{i}|(g_{\SCH}-g_{\MD})(E_i)|/N_E,$ and $e_2 = \sqrt{\sum_i((g_{\SCH}-g_{\MD})(E_i))^2/N_E},$ with
$N_E$ being the number of eigenvalues $E_i$.}
\label{tab:obs-conv}
\end{table}

\begin{figure}[htbp]
\centering
\includegraphics[width=0.4\textwidth]{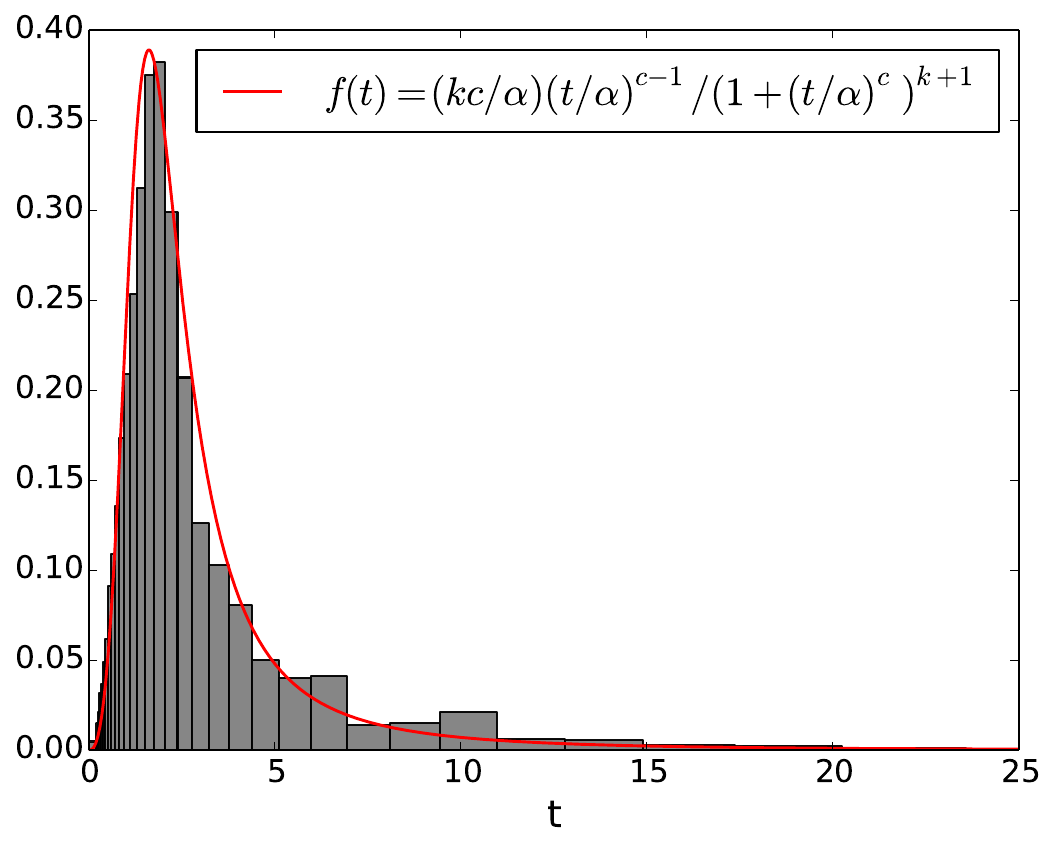}
   \caption{Normalized histogram of 
       the inter-arrival initialization times obtained by Algorithm~\ref{alg:compute-pe-md}
       when applied to the conical
       intersection problem with $a=(5,0)$ and $M=12800$, 
       and, for comparison, the probability density function of the Burr XII distribution
       with parameter values $\alpha = 1.65$, $c=3.6$ and $k=0.5$. We do not know why this distribution fits reasonably well to the data.}
\label{fig:hitting_times}
\end{figure}

\subsubsection{Approximation of the probability to be in excited states}
We approximate the probability to be in the excited state, $p_{ex}$,  from the Sch\"odinger equation
by the approximation, $\hat{\hat p}_{ex}$, from Algorithm 1, using 
Ehrenfest molecular dynamics simulations based on the  St\"ormer-Verlet method
as in~\eqref{eq:verlet-1d}.
We choose $t_0=0$, $P_0 = [1, \sqrt{2(E-\lambda_-(X_0))-1}]$,
$\psi_0 =  \Psi_-(X_0)$, $X_0 = [2, -0.5]$ (line intersection) and $X_0=[-2,0.5]$ (conical intersection). 
\begin{figure}[htbp]
 \centering
  \subfigure[Line intersection]{\label{pe-2D:subfig1}
    \includegraphics[width=0.35\textwidth,height=0.35\textwidth]{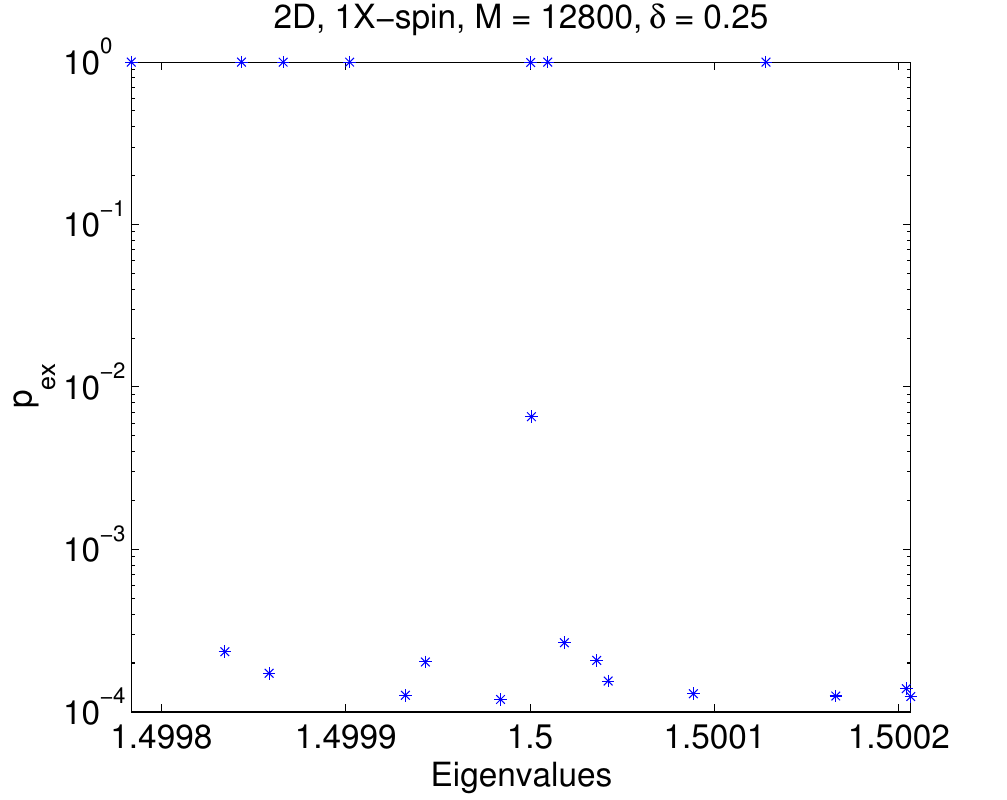} %
  }
  \qquad
  \subfigure[Line intersection]{\label{pe-2D:subfig2}
    \includegraphics[width=0.35\textwidth,height=0.35\textwidth]{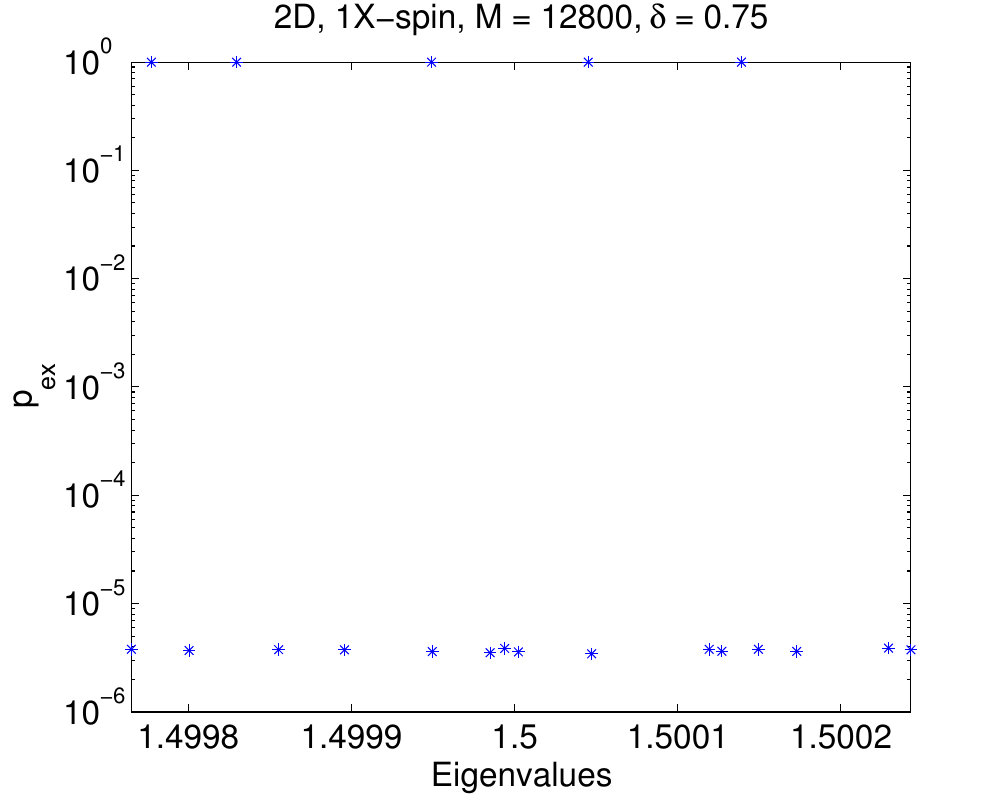} %
  }

  \subfigure[Conical intersection]{\label{pe-2D:subfig3}
    \includegraphics[width=0.35\textwidth,height=0.35\textwidth]{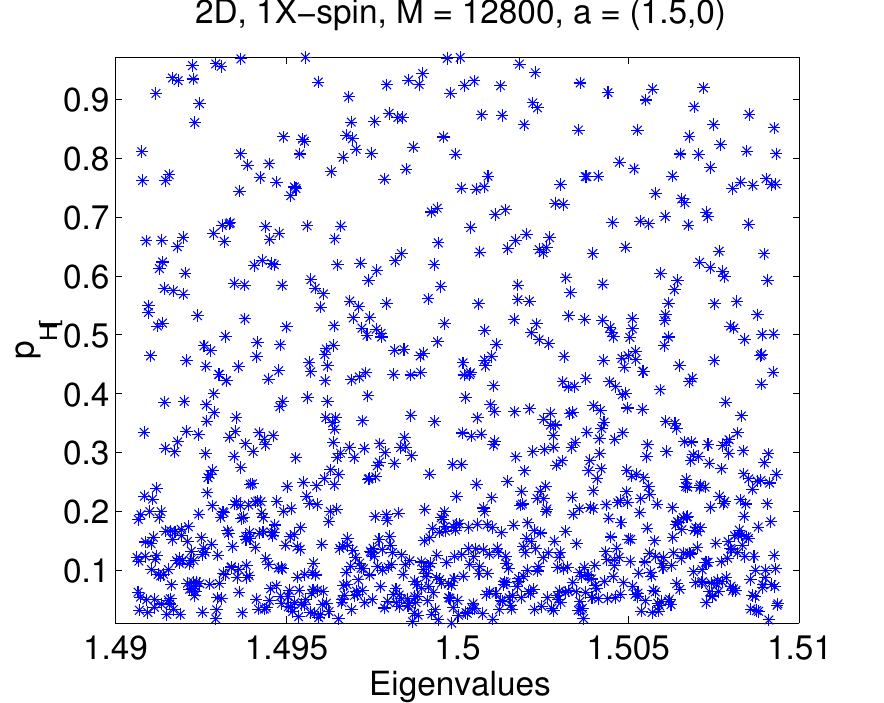}  %
  }
  \qquad
  \subfigure[Conical intersection]{\label{pe-2D:subfig4}
    \includegraphics[width=0.35\textwidth,height=0.35\textwidth]{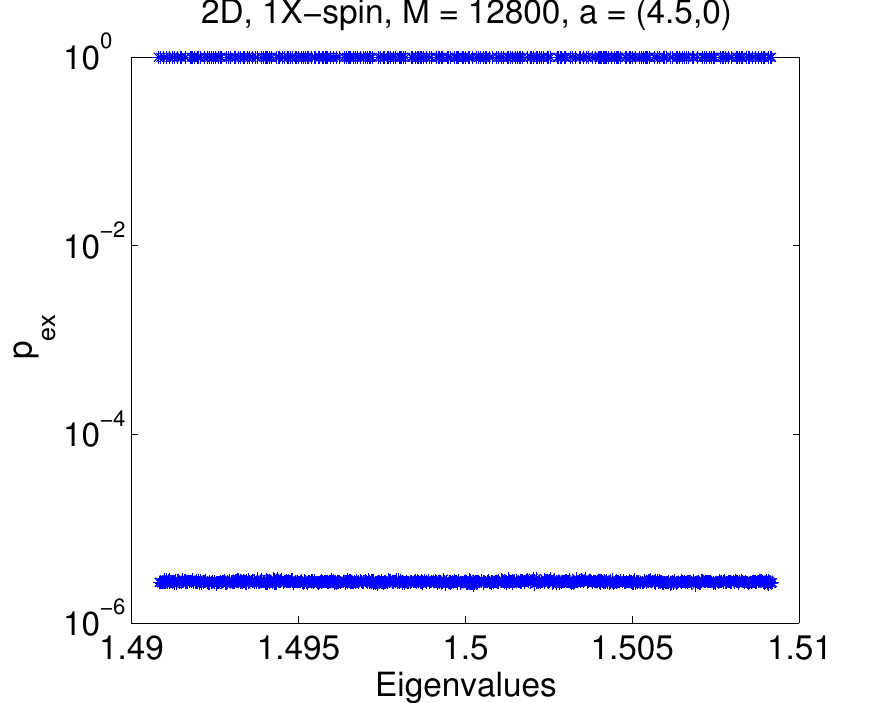}  %
  }
  \caption{Probabilities, $p_{ex}$, to be in the excited state
  as a function of the eigenvalues, $E_h$, 
  computed from the solution of the two dimensional Schr\"odinger equations 
  and based on the formula~\eqref{eq:Pr-excited}. The first row corresponds 
  to the line intersection for two different values of $\delta$; the second row corresponds 
  to the conical intersection case for the intersection point $a=(1.5,0)$ inside and 
  $a=(4.5,0)$ outside the {classically allowed} region. We note that the plots contain
  the probabilities  to be in the excited state, $p_{ex}$, for $20$ and $1000$ 
  eigenvalues for line and conical intersection cases, respectively.}
\label{fig:pe-2D}
\end{figure}

\begin{figure}[htbp]
 \centering
 \includegraphics[width=0.7\textwidth]{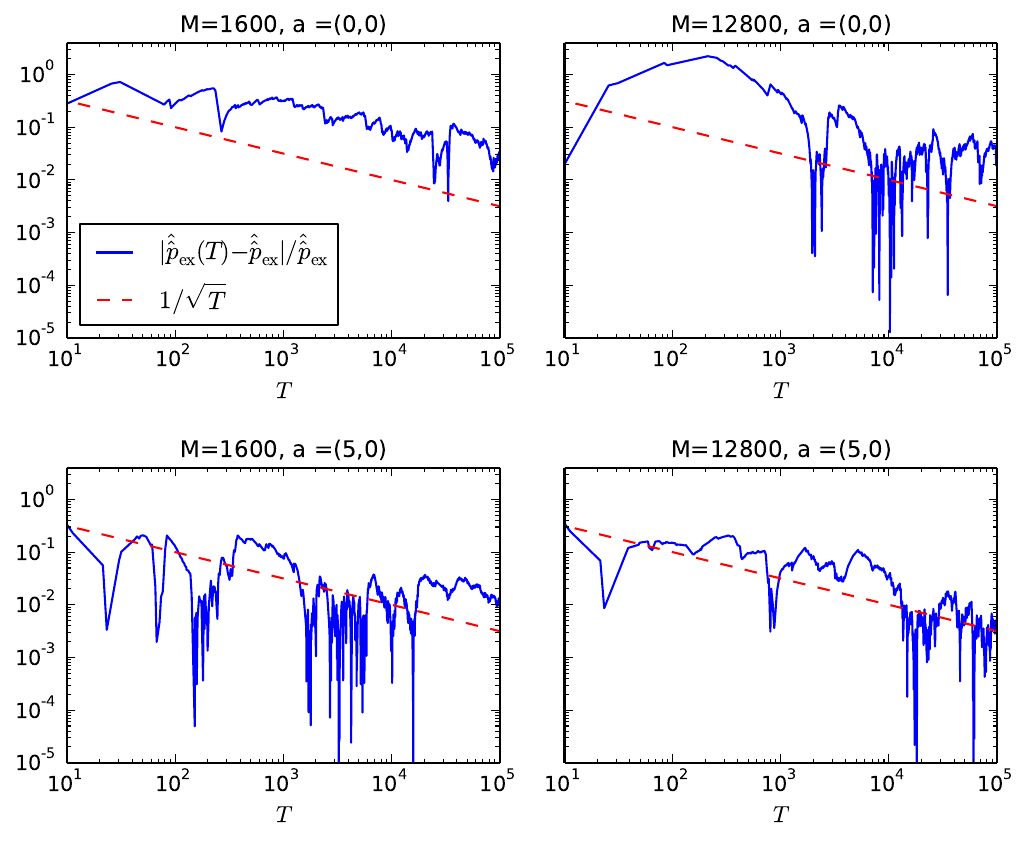}
  \caption{Convergence plots for $\hat{\hat p}_{ex}$ obtained for 
      the hyperplane sampling and molecular dynamics simulation methods in 
      Algorithm~\ref{alg:compute-pe-md} applied to the conical intersection 
      problem.}
  \label{fig:pe-convT-error-y1-x-plane}
\end{figure}

\subsubsection*{Conclusions.}

Figure~\ref{fig:pe-2D} shows the probabilities,~$p_{ex}$, to be in the
excited state determined from the discrete Schr\"odinger eigenvalue
problem~\eqref{eq:discrete-schrod-2d} based on the
formula~\eqref{eq:Pr-excited}. We note that for large gaps, in the
right column, the probabilities to be in the excited state, $p_{ex}$,
are either small or close to one, as expected, while for small gaps, 
in the left column, $p_{ex}$ takes values in the whole interval $[0,1]$. 
We also observe the lack of the structure present in the
analogous Figure~\ref{fig:computed-estimated} for the one dimensional
case, due to not having the resonances ordered with the eigenvalues as
in the one dimensional case discussed in Remark~\ref{rem:resonance}.

  Figure~\ref{fig:pe-Mpow0p5} shows $p_{ex}$ determined from the
  discrete Schr\"odinger eigenvalue problem and from molecular
  dynamics simulations for both line and conical intersection cases
  computed by Algorithm~\ref{alg:compute-pe-md}. The shapes of the
  respective methods' solutions are similar, especially when the
  probability to be in the exited state, $p_{ex}$, is small.
  Figure~\ref{fig:hitting_times} presents an empirical probability
  density of inter-arrival hitting times of hyperplanes,
  cf.~Algorithm~\ref{alg:compute-pe-md}. The empirical
  density is well approximated by a heavy tailed Burr XII
  probability density, which indicates that long
  inter-arrival hitting times are non-negligible.
  In Figure~\ref{fig:pe-convT-error-y1-x-plane}, we observe an
  expected decay of~$\mathcal O(T^{-1/2})$ for the relative error
  $|\hat{\hat p}_{ex}(T)-\hat{\hat p}_{ex}|/\hat{\hat p}_{ex}$ when studying the conical
  intersection problem with $a=(5,0)$. In the more unstable parameter
  setting of $a=(0,0)$, however, the error decay is slower with an
  observed rate of convergence slightly slower than the expected
  $-1/2$. This might be due to a longer correlation length in the
  molecular dynamics and hyperplane sampling of
  Algorithm~\ref{alg:compute-pe-md} when the parameter values of $a$
  are close to $(0,0)$.

\begin{remark} \label{rem_HD}[Extension to higher dimension] {\rm
The formation of a conical intersection for a $2\times 2$ symmetric
matrix $V$ requires the eigenvalues to be equal, which means that the
off diagonal element is zero and the diagonal elements are equal at an
intersection point.  To have a conical intersection point is therefore
generic in dimension two and in higher dimensions the intersection is
typically a co-dimension two set. An implementation of the algorithm
for computing $\hat{\hat p}_{ex}$ in higher dimensions and for potentials
used in chemistry is future work. The relevant geometry in the
multi-dimensional case is the distance to this co-dimension two set,
which makes the link to two dimensional case studied here.}
\end{remark}

\clearpage
\section{Regularity and Weyl calculus}\label{sec:weyl}
This section establishes four lemmas where the last three are used in Theorem \ref{thm:potential} and  \ref{ergod_sats}.
The first lemma proves the  
bound \eqref{hatc_first} on derivatives of the solution $u(t,Z)=\mathbb E[g\circ S_{Tt}(Z)]$ to \eqref{kb},
which verifies that  assumption (v)   in Theorem \ref{thm:potential} holds for a case with a near crossing of eigenvalues.
The second lemma derives $L^2$-bounds on Weyl quantized operators obtained from mollified symbols,
the third and fourth lemma estimate remainder terms related to the Weyl quantization and establishes the Moyal expansion
for the setting used in the proofs of Theorem \ref{thm:potential} and Theorem \ref{ergod_sats}.

We will use equation \eqref{ito_eq} written as the Ito equation
\begin{equation}\label{zekvation}
dZ_t=a(Z_t)dt + b(Z_t)d\tilde W_t\, 
\end{equation}
for the molecular dynamics path $Z:[0,T]\times\Omega\rightarrow \rset^{6N}$, where $\Omega$ is the sample space
and assume that the Jacobians $a'$ and $b'$ and the higher order derivatives $\partial_Z^{\alpha}a(Z)\in \rset^{(6N)^{(|\alpha|+1)}}$, with $|\alpha|=n$ for $n=1,2,3,\ldots $, satisfy the $\ell^2$ bounds
\begin{equation}\label{assump_aprim}
\begin{split}
\|a'\|_{2,\infty}+\|b'\|_{2,\infty}^2&=\mathcal O(\delta^{-1})\, ,\\
\|\partial_Z^{\alpha}a\|_{2,\infty}+\|\partial_Z^{\alpha}b\|_{2,\infty}^2&=\mathcal O(\delta^{-|\alpha|})\, ,\\
\end{split}
\end{equation}
uniformly in $N$ (but not in $|\alpha|$).
Here we use the $\ell^2$ (Frobenius) norm: assume that $|\alpha|=n$ then
$\|\partial_Z^\alpha a\|_2=\sqrt{\sum_{|\alpha|=n}\sum_{i=1}^{6N}|\partial_Z^\alpha a_i|^2}$ and the notation
\[
\|a(Z_t)\|_{2,\infty}:=\sup_{\omega\in \Omega} \|a\big(Z_t(\omega)\big)\|_2\, .
\]

\begin{lemma}\label{lemma_time}
Assume \eqref{assump_aprim} and that the near crossing is weak, with
non vanishing velocity $P$ across the near crossing, namely that %
\begin{equation}\label{cdeupp}
\begin{split}
\max_{X\in\rset^{3N}}\|\lambda_\eta''(X)\|_2&\le c_\delta\delta^{-1}\, ,\\
 \int_0^T\|\lambda_\eta''(X_t)\|_{2,\infty} dt&\le CT\, ,\\
 \int_0^T(\mathbb E[\|\partial_{X}^{\alpha}\lambda_\eta''(X_t)\|^{2n}_2])^{1/(2n)} \, dt&\le C_\alpha T\delta^{-|\alpha|}\, ,\\
\end{split}
\end{equation}
for  a constant $c_\delta\le 1/\log \delta^{-1}$ and  constants $C_\alpha$ independent of $\delta$, $M$, $N$ and $\epsilon$.
Then  for each multi index $\beta$ %
there is a constant $\hat C$, independent of $M$, $\delta$, $N$ and $\epsilon$, such that
\begin{equation}\label{ctbound}
\sup_{Z_0\in\rset^{6N}} \| \partial_{Z_0}^\beta\mathbb E[g\circ S_T(Z_0)\, |\ Z_0]\|_2
\le e^{\hat C (T+1)}\delta^{\min(0,-|\beta|+1)}\, .
\end{equation}
\end{lemma}

\begin{proof}
Remark \ref{c_rem} below motives the assumptions.
First we study solutions $w:[0,T]\times \Omega\rightarrow \rset^{6N}$ to the linearized  Ito equation
\[
dw_t=a'(Z_t)w_t dt+b'(Z_t)w_td\tilde W_t +\alpha_tdt+\beta_td\tilde W_t
\]
where $Z_t$ solves \eqref{zekvation}=\eqref{ito_eq} and 
the stochastic processes $\alpha:[0,T]\times \Omega\rightarrow\rset^{6N}$
and $\beta:[0,T]\times\Omega \rightarrow\rset^{6N\times 6N}$ satisfy for each $n\in \mathbb N\setminus 0$ and
some positive number $m$ the bound 
\begin{equation}\label{alfabeta}
\begin{split}
\int_0^t \mathbb E[\|\beta_s\|_{2n}^{2n}\ |\ Z_0] ds &=\mathcal O(e^{Ct})\, ,\\
\int_0^t \big(\mathbb E[\|\alpha_s\|_2^{2n}\ |\ Z_0]\big)^{1/(2n)} ds &=\mathcal O(e^{Ct}\delta^{-m})\, .\\
\end{split}
\end{equation}
Ito's formula implies
\[
\begin{split}
d\mathbb E[\|w_t\|_2^2\ | Z_0]
&= 2\mathbb E[w_t\cdot dw_t\ | Z_0] + \mathbb E[dw_t\cdot dw_t\ | Z_0]\\
&\le 
2\mathbb E[\|a'(Z_t)\|_2 \|w_t\|_2^2\ | Z_0] dt + \mathbb E[\|b'(Z_t)\|_2^2  \|w_t\|_2^2\ | Z_0] dt\\
&\qquad +2\mathbb E[ \|w_t\|_2\|\alpha_t\|_2 \ | Z_0] dt
+\mathbb E[\|\beta\|_2^2+2\|\beta_t\|_2\|b'_t\|_2\|w_t\|_2\ |\ Z_0]dt\\
&\le (2\|a'(Z_t)\|_{2,\infty}+ \|b'(Z_t)\|_{2,\infty}^2)\mathbb E[ \|w_t\|_2^2\ | Z_0] dt+ 
\mathbb E[2\|w_t\|_2\|\alpha_t\|_2+\|\beta_t\|_2^2\ |\ Z_0]dt
\end{split}
\]
where by assumption \eqref{assump_aprim} %
we have $2\|a'_s\|_{2,\infty}+ \|b'_s\|_{2,\infty}^2=\mathcal O(\delta^{-1} + \epsilon\delta^{-2})=\mathcal O(\delta^{-1} )$ and by \eqref{cdeupp} there holds 
$\int_0^t2\|a'_s\|_{2,\infty}+ \|b'_s\|_{2,\infty}^2ds=\mathcal O(t)$.
Here $\|b'\|_{2,\infty}^2=\sum_i\|b'_{i\cdot\cdot}\|_{2,\infty}^2$.
We use $\epsilon\ll \delta$, with the arbitrarily small diffusion constant $\epsilon$ in \eqref{sde_proj},  so that the contribution from $b'$ and $\beta$ are negligible compared to those from $a'$ and $\alpha$.
Gronwall's inequality shows that
\begin{equation}\label{gronwall}
\begin{split}
\mathbb E[ \|w_t\|_2^2\ | Z_0]
&\le \big(\mathbb E[ \|w_0\|_2^2\ | Z_0]+ \int_0^t\mathbb E[2\|w_s\|_2\|\alpha_s\|_2+\|\beta_s\|_2^2\ |\ Z_0]ds\big) e^{\int_0^t2\|a'_s\|_{2,\infty}+ \|b'_s\|_{2,\infty}^2ds}\\
&\le \big(\mathbb E[ \|w_0\|_2^2\ | Z_0]+ \int_0^t 
2\sqrt{\mathbb E[\|w_t\|_2^2\ |\ Z_0]}\sqrt{\mathbb E[\|\alpha_s\|_2^2\ |\ Z_0]}
+\mathbb E[\|\beta_s\|_2^2\ |\ Z_0]ds\big)e^{Ct}
\end{split}
\end{equation}
which by \eqref{alfabeta} implies 
\begin{equation}\label{wtbound}
\begin{split}
\sup_{0<s<t}\sqrt{\mathbb E[ \|w_s\|_2^2\ | Z_0]}
&\le e^{Ct} \big(\sqrt{\mathbb E[ \|w_0\|_2^2\ | Z_0]}+ \int_0^t 
2\sqrt{\mathbb E[\|\alpha_s\|^2_2\ |\ Z_0]} ds+ \mathcal O(e^{Ct})\big) \\
&=e^{Ct}\mathcal O(\delta^{-m})
\end{split}
\end{equation}
provided this holds initially for $t=0$.
Analogously for $n>1$, we obtain by H\"older's inequality
\[
\begin{split}
d\mathbb E[\|w_t\|_2^{2n}\ | Z_0] &= n \mathbb E[\|w_t\|_2^{2(n-1)} w_t\cdot dw_t  + n(n-1)\|w_t\|_2^{2(n-2)} \frac{(w_t\cdot dw_t)^2}{2}
+ n\|w_t\|_2^{2(n-1)} \frac{dw_t\cdot dw_t}{2}\ | Z_0]\\
&= n\mathbb E[\|w_t\|_2^{2(n-1)} (w_t\cdot a_t'w_t dt + w_t\cdot\alpha_t dt) + h.o.t.\ |\ Z_0]\\
&\le n \mathbb E[\|w_t\|_2^{2n} \|a'_t\|_2  \ |\ Z_0]dt + n \mathbb E[\|w_t\|_2^{2n-1} \|\alpha_t\|_2  \ |\ Z_0] dt+ h.o.t.\\
&\le n \mathbb E[\|w_t\|_2^{2n} \|a'_t\|_2  \ |\ Z_0] dt 
+n (\mathbb E[\|w_t\|_2^{2n}  \ |\ Z_0])^{1-1/(2n)} (\mathbb E[\|\alpha_t\|_2^{2n}  \ |\ Z_0])^{1/(2n)} dt  + h.o.t.\\
\end{split}
\]
so that
\begin{equation}\label{wtbound2}
\begin{split}
\sup_{0<s<t}(\mathbb E[ \|w_s\|_2^{2n}\ | Z_0])^{1/(2n)}
&\le e^{Ct} \big((\mathbb E[ \|w_0\|_2^{2n}\ | Z_0])^{1/(2n)}+ \int_0^t 
2(\mathbb E[\|\alpha_t\|_2^{2n}\ |\ Z_0])^{1/(2n)} ds+ \mathcal O(e^{Ct})\big) \\
&=e^{Ct}\mathcal O(\delta^{-m})
\end{split}
\end{equation}
provided this holds initially for $t=0$.

The first variation $Z'_t=\partial Z_t/\partial Z_0$ satisfies
\[
dZ'_t=a'(Z_t)Z'_tdt + b'(Z_t)Z'_t d\tilde W_t\, ,
\]
so that \eqref{wtbound}, for $m=0$, implies 
\[
\begin{split}
 \mathbb E[\| Z'_t\|_{2}^{2n}\ |\ Z_0]&=\mathcal O(e^{Ct})\, .\\
\end{split}
\]

The second variation satisfies
\[
dZ_t''=a'(Z_t)Z_t''dt + b'(Z_t)Z''_td\tilde W_t + a''(Z_t)Z_t'Z_t'dt + b''(Z_t)Z_t'Z_t'd\tilde W_t\, .
\]
We have for $\alpha=a''Z'Z'$ that
\begin{equation}\label{ijkrq}
\begin{split}
\|a''Z'Z'\|_2^2&=\sum_{ijk}(\sum_{rq} a''_{irq}Z'_{rj}Z'_{qk})^2\\
 &\le \sum_{ijk}\sum_{rq} (a''_{irq})^2\sum_{rq} (Z'_{rj}Z'_{qk})^2\\
&=\sum_{irq} (a''_{irq})^2\sum_{rqjk} (Z'_{rj}Z'_{qk})^2\\
&=\|a''\|_2^2\|Z'\|_2^4\, ,
\end{split}
\end{equation}
and similarly for $\beta$,
which combined with
Ito's formula and Gronwall's inequality  as above imply
\[
\sqrt{\mathbb E[\|Z''_t\|_{2}^2\ | Z_0]}=\mathcal O(\delta^{-1}e^{Ct}) %
\]
using \eqref{cdeupp},
that \eqref{wtbound2} yields $(\mathbb E[\|Z'_t\|_2^8\ |\ Z_0])^{1/4}=\mathcal O(e^{Ct})$
and
\[
\int_0^t \sqrt{\mathbb E[\|a''_sZ_s'Z_s'\|_2^2\ |\ Z_0]} ds\le
\int_0^t (\mathbb E[\|a''_s\|_2^4\ |\ Z_0])^{1/4} (\mathbb E[\|Z'_s\|_2^8\ |\ Z_0])^{1/4}ds\le e^{Ct}\mathcal O(\delta^{-1})\, .
\]
Also the higher variations satisfy
\[
dZ_t^{(m)}= a'(Z_t)Z_t^{(m)}dt + b'(Z_t)Z_t^{(m)}d\tilde W_t +\alpha_tdt +\beta_td\tilde W_t
\]
where the stochastic processes $\alpha$ and $\beta$ have the bound \eqref{alfabeta},
since  following 
\eqref{ijkrq} we have 
\[\begin{split}
\|a''Z''Z'\|_2^2&\le \|a''\|_2^2\|Z''\|_2^2\|Z'\|_2^2\, , \\
\|a'''Z'Z'Z'\|_2^2&\le \|a'''\|_2^2\|Z'\|_2^6\, ,
\end{split}
\]
and similarly for the higher order derivatives.
Therefore \eqref{wtbound} and \eqref{wtbound2} yield 
\begin{equation}\label{z'''}
\begin{split}
(\mathbb E[\|Z^{(m)}_t\|_{2}^{2n}\ |\ Z_0])^{1/(2n)}&=\mathcal O(e^{Ct}\delta^{(1-m)})\, .\\
\end{split}
\end{equation}
We have by Jensens inequality, for $|\beta|=n$ ,
 \[
 \begin{split}
\sup_{Z_0} \|\mathbb{E} [\partial_{Z_0}^\beta g(Z_t)\ |\  Z_0]\|_2 &\le
\sup_{Z_0}\mathbb{E}[\|\partial_{Z_0}^\beta g(Z_t)\|_2\ |\  Z_0] \\
  &\le \sup_{Z_0}\sqrt{\mathbb E[\|\partial_{Z_0}^\beta g(Z_t)\|_2^2\ |\  Z_0] }\\
  &=  \left\{
  \begin{array}{cc}
  \sup_{Z_0}\sqrt{\mathbb E[\|g'_tZ'_t\|_2^2\ |\  Z_0]}\, ,  & n=1\\
 \sup_{Z_0}\sqrt{\mathbb E[\|g'_tZ''_t+ g''_tZ_t'Z'_t\|_2^2\ |\  Z_0]}\, ,  & n=2\\
  \sup_{Z_0}\sqrt{\mathbb E[\|g_t'Z_t'''+ 2g_t''Z_t''Z_t'+g_t'''Z_t'Z_t'Z_t'\|_2^2\ |\  Z_0]}\, ,  & n=3\\
 \end{array} \right.
 \\  
 \end{split}
 \]
which  by \eqref{z'''} proves \eqref{ctbound}.
\end{proof}

\begin{remark}\label{c_rem}{\rm
We motivate assumption \eqref{cdeupp} by the behavior for the example of the avoided conical intersection 
 eigenvalues $\lambda_0(X)=\pm c_\delta\sqrt{|X|^2+\delta}$ for $X\in\rset^2$ or
$X\in\rset^1$ discussed in Section \ref{sec:2d-example} and \eqref{tilde_v_def}.
Although $\max_{X\in K}\|\lambda_0''(X_t)\|_2=\mathcal O(\delta^{-1})$, we assume
the near crossing is weak, i.e. that $\|\lambda_0''(X)\|_2\le c_\delta(\mbox{dist}(X,a)^2+\delta^2|)^{-1/2}$ for some constant 
$c_\delta\le 1/\log\delta^{-1}$ and a codimension $2$ set $a\subset\rset^{3N}$. Then the exponent in Gronwall's inequality has the bound
\[
\int_t^T\|\lambda_0''(X_s)\|_2ds\le C(T-t) + c_\delta\int_0^1 \frac{ds}{\sqrt{c^2s^2+\delta^2}}\le C(T-t)+c^{-1/2}c_\delta\log\delta^{-1}\]
provided the near crossing of potential surfaces is located away from the part of $\Sigma_E$ where $|P|=0$.
Here we assume  that the velocity through the avoided crossing domain is bounded from below by $c>0$.
Similarly, an assumption $\|\partial_X^\alpha \lambda_0''(X)\|_2\le C(\mbox{dist}(X,a)^2+\delta^2|)^{-(|\alpha|+1)/2}$ satisfies the last estimate in \eqref{cdeupp}.}
\end{remark}

We will use the $s$-quantization and the notation
\[
\int_{\rset^{3N}} \langle \Phi(X), \OP^s[ h]\Phi(X)\rangle dX
=\int_{\rset^{6N}} h(Z)\cdot W^{(s)}(Z) dZ
\]
for any smooth symbol $h$ and the Wigner functions $W=W^{(1/2)}$ and $W^{(s)}$ defined in \eqref{w_s_def}.
\begin{lemma}\label{u_lemma} 
Assume that $h:\rset^{6N}\to \mathbb C$ and $f:\rset^{6N}\to \mathbb C$
can be written as $h=\bar h*\phi_\eta$ and $f=\bar f*\phi_\eta$,
where $\phi_\eta(z):= (2\pi\eta)^{-3N} e^{-|z|^2/(2\eta)}$,
and  that $\bar h$ and $\bar f$ have continuous second derivatives with polynomial growth, i.e. for some $n\in \mathbb Z$ and $m\in\mathbb Z$ 
and $|\alpha|\le 2$ there is a constant $C$ such that 
 $\partial_z^\alpha\bar h(z)\le C(1+|z|^2)^m$ and $\partial_z^\alpha\bar f(z)\le C(1+|z|^2)^n$ uniformly in $z$, then for every $s\in[0,1]$
\[
 |\int_{\rset^{3N}} \langle \Phi(X), \OP^s[ h(z) f(z)]\Phi(X)\rangle dX|=
  \sup_{z\in\rset^{6N}}|h(z)f(z)| + \mathcal O(\eta)\, ,
  \]
  uniformly in $N$.
\end{lemma}

\begin{proof} The proof  to estimate the $L^2(\rset^{3N})$ operator norm for a symbol $r$ has three steps: to determine a representation of the solution of $\int_{\rset^{3N}} \langle \Phi,\OP^{s}[r]\Phi\rangle dX$ using the FBI transform $T\Phi$, 
to calculate a Fourier multiplicator  applied to $r$  using the representation and in the third step to 
verify that the product of two FBI transformed functions is the convolution of
the Wigner function with a Gaussian.

{\bf 1.} 
We will use that the components of the Wigner function $W=W^{(1/2)}$ convolved with the Gaussian 
$\phi_{M^{-1/2}}=(2\pi M^{-1/2})^{-3N}e^{-(|X|^2+|P|^2)M^{1/2}/2}$
is the product of the FBI transforms \eqref{fbi_def} of the components $\Phi_i$ and $\Phi_j$
\begin{equation}\label{tensor}
W_{ij}*\phi_{M^{-1/2}}(X,P)= T\Phi_i(X,P) \overline{T\Phi_j(X,P)}
\end{equation}
as verified in Step 3 below.
Consequently, the diagonal entries are non negative
\[
W_{ii}*\phi_{M^{-1/2}}(X,P)= |T\Phi_i^*(X,P)|^2\ge 0\, .
\]

An $s_*$-quantization remainder $\hat r^{s_*}$ can be related to a remainder 
$\hat r=\OP[h]$
in the Wigner quantization (with $s=1/2$)
by 
\[
\int_{\rset^{6N}} r(X,P)\cdot W^{(s_*)}(X,P) dXdP
=\int_{\rset^{6N}} r^{s_*}(X,P)\cdot W^{(1/2)}(X,P) dXdP
\]
where for $s_*\in[0,1]$
\[
r^{s_*}(X,P)= e^{-iM^{-1/2}(s_*-\frac{1}{2}) \nabla_X\cdot\nabla_P} r(X,P)\, ,
\]
which is proved  in \cite{martinez_book} Remark 2.7.3 and related to the expansion of the exponential in \eqref{cab}.

We will introduce the convolution with the Gaussian in the estimate of the remainder using the Taylor expansion of the exponential
\[
\begin{split}
\int_{\rset^{3N}} \langle \Phi, \hat r^{s_*} \Phi\rangle dX
&=\int_{\rset^{6N}} r^{s_*}(X,P)\cdot W^{(1/2)}(X,P) dXdP\\
&=\int_{\rset^{6N}} \mathcal F r^{s_*}(\xi_X,\xi_P)\cdot \mathcal F W^{(1/2)}(\xi_X,\xi_P) d\xi_Xd\xi_P\\
&=\int_{\rset^{6N}} \big( \mathcal F \phi_{M^{-1/2}}(\xi_X,\xi_P)\big)^{-1} \mathcal F r^{s_*}(\xi_X,\xi_P)\cdot \mathcal F W^{(1/2)}(\xi_X,\xi_P) \mathcal F \phi_{M^{-1/2}}(\xi_X,\xi_P) d\xi_X d\xi_P\\
&=\int_{\rset^{6N}}   r^{*}(X,P)
\cdot  W^{(1/2)}* \phi_{M^{-1/2}}(X,P) dXdP
\, ,\\
 \end{split}
 \] 
 where the function $r^*$ is defined by $ r^*= \mathcal F^{-1}\{\big( \mathcal F \phi_{M^{-1/2}}(\xi_X,\xi_P)\big)^{-1} \mathcal F r^{s_*}(\xi_X,\xi_P)\}$.
 Here we use the notation $r\cdot W:=\sum_{ij} r_{ij}W_{ij}$.
The tensor product property \eqref{tensor} of $W*\phi_{M^{-1/2}}$ implies the matrix norm estimate
\begin{equation}\label{r*est}
\begin{split}
\int_{\rset^{3N}} \langle \Phi(X), \hat r^{s_*} \Phi(X)\rangle dX &=
\int_{\rset^{6N}}   r^*(X,P)
\cdot  W^{(1/2)}* \phi_{M^{-1/2}}(X,P) dXdP\\
&=
\int_{\rset^{6N}}   \langle T\Phi(X,P), r^*(X,P)T\Phi(X,P)\rangle dXdP\\
&\le \int_{\rset^{6N}}   \| r^*(X,P)\|_2 |T\Phi(X,P)|^2  dXdP\\
&\le   \sup_{(X,P)\in\rset^{6N}}\| r^*(X,P)\|_2 \int_{\rset^{6N}}  |T\Phi(X,P)|^2  dXdP\\
&=  \sup_{(X,P)\in\rset^{6N}}\| r^*(X,P)\|_2 \int_{\rset^{3N}}  \langle T^*T\Phi,\Phi\rangle  dX\\
&=  \sup_{(X,P)\in\rset^{6N}}\| r^*(X,P)\|_2  
\underbrace{\int_{\rset^{3N}}  \langle \Phi(X),\Phi(X)\rangle  dX}_{=1}\\
\end{split}
\end{equation}
using the $L^2$ identity $T^*T=I$  in \eqref{fbi_iso} for the FBI transform
and the $\ell^2$ matrix norm $\|r^*(X,P)\|_2$ (or the less sharp Frobenius norm).

{\bf 2.} We have $\mathcal Fr^*=e^{|\omega|^2M^{-1/2}/2}e^{i\omega_x\cdot \omega_pM^{-1/2}}\mathcal F r(\omega)$. The goal is to absorb this exponentially growing pre-factor
$e^{|\omega|^2M^{-1/2}/2}e^{i\omega_x\cdot \omega_pM^{-1/2}}$ in $\mathcal F r(z)$. 

Consider first the case $r=\bar r*\phi_\eta$, then $\mathcal Fr^*=e^{|\omega|^2M^{-1/2}/2}e^{i\omega_x\cdot \omega_pM^{-1/2}}\mathcal F \phi_\eta(\omega)\mathcal F \bar r(\omega)$.
We have $\mathcal F \phi_\eta(\omega)=e^{-|\omega|^2\eta/2}$ so that
\[
e^{|\omega|^2M^{-1/2}/2} \ e^{-|\omega|^2\eta/2}= e^{-|\omega|^2{\nu}/{2}}=\mathcal F \phi_\nu(\omega)
\]
where $\nu=\eta-M^{-1/2}>0$. 
The  inverse Fourier transform of the pre factor  
$e^{is\omega_x\cdot\omega_p M^{-1/2}}\mathcal F\phi_\nu(\omega) $ is
\[\begin{split}
&(2\pi)^{-6N}\int_{\rset^{6N}} e^{-\nu|\omega_x|^2/2 -\nu|\omega_p|^2/2 +i\omega_x\cdot\omega_p sM^{-1/2} +i\omega_x\cdot x+i\omega_p\cdot p} d\omega_x d\omega_p\\
&=e^{-\frac{1}{2\nu(1+s^2M^{-1}\nu^{-2})} 
(|p|^2 + |x|^2)}   e^{-ip
\cdot x sM^{-1/2}\frac{1}{s^2M^{-1}+\nu^2}} (2\pi \nu)^{-3N} (1+\frac{s^2}{\nu^2 M})^{-3N/2}
=:\phi_{\nu,s}\, ,
\end{split}
\]
so that $r^*=\bar r*\phi_{\nu,s}$.
The estimate \eqref{r*est} implies then in the case $r=\bar r*\phi_\eta$
\[
\begin{split}
\int_{\rset^{3N}} \langle \Phi(X), \hat r^{s_*} \Phi(X)\rangle dX &\le  \sup_{z\in\rset^{6N}}\| r^*(z)\|_2\\
&=\sup_{z\in\rset^{6N}} \|\bar r * \phi_{\nu,s}(z)\|_2  \, .
\end{split}
\]
If $\bar r$ is uniformly bounded in $\mathcal C^2(\rset^{6N})$ 
we have by the two moments $\int_{\rset^{6N}} \phi_\eta(z)dz=1$ and $\int_{\rset^{6N}} z\phi_{\eta}(z)dz=0$ that
\[
r-\bar r=\mathcal O(\eta)
\]
and by the properties $\int_{\rset^{6N}} \phi_{\nu,s}(z)dz=1$ and $\int_{\rset^{6N}} z\phi_{\nu,s}(z)dz=0$ we obtain
\[
\bar r*\phi_{\nu,s}-\bar r=\mathcal O(\eta^{})
\]
so that 
\begin{equation}\label{r_en}
\begin{split}
 \int_{\rset^{3N}} \langle \Phi(X), \hat r^{s_*} \Phi(X)\rangle dX 
&= \sup_{z\in\rset^{6N}} |r(z)| + \mathcal O(\nu^{})\, .\\
\end{split}
\end{equation}

In the case $r=(\bar h*\phi_\eta) \, (\bar f*\phi_\eta)$ we have similarly
\[
\begin{split}
&(2\pi)^{-12N}\int_{\rset^{12N}} \mathcal F\bar h(\omega')\mathcal F \bar f(\omega-\omega') e^{-|\omega'|^2\eta/2-|\omega-\omega'|^2\eta/2}
e^{i\omega_x\cdot\omega_p sM^{-1/2} +|\omega|^2/(2M^{-1/2}) }e^{i\omega\cdot z} d\omega'd\omega\\
&=(2\pi)^{-12N}\int_{\rset^{12N}} \mathcal F \bar h(\omega')\mathcal F \bar f(\omega-\omega') e^{-|\omega|^2\eta/4-|\omega'-\omega/2|^2\eta}
e^{i\omega_x\cdot\omega_p sM^{-1/2} +|\omega|^2/(2M^{-1/2}) }e^{i\omega\cdot z} d\omega'd\omega\\
&=(2\pi)^{-9N}\eta^{-3N}\int_{\rset^{12N}}\mathcal F^{-1}\big\{\mathcal F \bar h\, \mathcal F\{\bar f(\omega-\cdot)\}\big\}(z') e^{-|z'|^2/(4\eta) 
+i\omega\cdot z'/2}
e^{i\omega_x\cdot\omega_p sM^{-1/2} -\nu|\omega|^2/2 }e^{i\omega\cdot z} d\omega dz'\\
&=\underbrace{(2\pi\eta)^{-3N}(2\pi\nu)^{-3N}(1+\frac{s^2}{\nu^2M})^{-3N/2}}_{=:\sigma_{N,M}}
\times\\
&\qquad\times
\int_{\rset^{12N}} \bar h(z'-z'') \bar f(-z'') e^{-|z'|^2/(4\eta)} 
e^{-|z-z''+z'/2|^2/(4\nu(s))} e^{i\frac{(z_x-z''_x+z'_x/2)\cdot (z_p-z''_p+z'_p/2)  sM^{-1/2}}{\nu^2+s^2M^{-1}}} dz'dz''\\
&= \sigma_{N,M} %
\int_{\rset^{12N}} \bar h(v) \bar f(w) e^{-|w-v|^2/(4\eta)} 
e^{-|z+(w+v)/2|^2/(4\nu(s))} e^{\frac{i(z_x+(w_x+v_x)/2)\cdot (z_p+(w_p+v_p)/2) sM^{-1/2}}{\nu^2+s^2M^{-1}}} dvdw\\
&= \sigma_{N,M} %
\int_{\rset^{12N}} \bar h(z-v) \bar f(z-w) e^{-|w-v|^2/(4\eta)} 
e^{-|(w+v)/2|^2/(4\nu(s))} e^{\frac{i(w_x+v_x)\cdot (w_p+v_p) sM^{-1/2}}{4(\nu^2+s^2M^{-1})}}\, dvdw\\
&=:Q(\bar h\bar f)(z)\, . %
\end{split}
\]
We have by the Fourier transform
\[
\begin{split}
\sigma_{N,M}
\int_{\rset^{12N}}  e^{-|w-v|^2/(4\eta)} 
e^{-|(w+v)/2|^2/(4\nu(s))}
e^{\frac{i(w_x+v_x)\cdot (w_p+v_p) sM^{-1/2}}{4(\nu^2+s^2M^{-1})}} dvdw&=1\\
\int_{\rset^{12N}} v e^{-|w-v|^2/(4\eta)} 
e^{-|(w+v)/2|^2/(4\nu(s))}
e^{\frac{i(w_x+v_x)\cdot (w_p+v_p) sM^{-1/2}}{4(\nu^2+s^2M^{-1})}} dvdw&=0\\
\int_{\rset^{12N}} w e^{-|w-v|^2/(4\eta)} 
e^{-|(w+v)/2|^2/(4\nu(s))}
e^{\frac{i(w_x+v_x)\cdot (w_p+v_p) sM^{-1/2}}{4(\nu^2+s^2M^{-1})}} dvdw&=0\\
\end{split}
\] so that
\[
|Q(\bar h\bar f)(z)- \bar h(z)\bar f(z)|=|Q\big(\bar h(\cdot)\bar f(\cdot)-\bar h(z)\bar f(z)\big)|\le C 
(\|\frac{\bar h(z)}{1+|z|^m}\|^2_{\mathcal C^2(\rset^{6N})}+\|\bar f(z)(1+|z|^m) \|^2_{\mathcal C^2(\rset^{6N})}){\nu(s)}\, .
\]
Similarly using  the moments  $\int_{\rset^{6N}}\phi_\eta(z)dz=1$ and $\int_{\rset^{6N}}z\phi_\eta(z)dz=0$ and the second
derivatives $\|\bar h\|^2_{\mathcal C^2(\rset^{6N})} + \|\bar f\|^2_{\mathcal C^2(\rset^{6N})}\le C$ we obtain
\[
h(z)f(z)-\bar h(z)\bar f(z)=\mathcal O(\eta)\, ,
\]
and we conclude that $|Q(\bar h\bar f)(z)-  h(z) f(z)|= \mathcal O(\eta)$ which combined with \eqref{r*est}, as in \eqref{r_en},
proves the theorem.

{\bf 3.} This step verifies \eqref{tensor} following the proof of Proposition 1 in \cite{lasser_keller}.
We have 
\[
\begin{split}
&W_{ij}*\phi_{M^{-1/2}}(X,P)\\
&=(\pi M^{1/2}\times 2\pi M^{-1/2})^{-3N} \int_{\rset^{9N}} e^{iM^{1/2} Y\cdot P'} \Phi_i(X'-\frac{Y}{2})\Phi_j^*(X'+\frac{Y}{2})
e^{-(|X-X'|^2+|P-P'|^2)M^{1/2}}dY dX'dP'\\
&=(\pi M^{1/2})^{-3N/2}(2\pi M^{-1/2})^{-3N} \int_{\rset^{6N}}
e^{iM^{1/2} Y\cdot P} \Phi_i(X'-\frac{Y}{2})\Phi_j^*(X'+\frac{Y}{2})
e^{-(|X-X'|^2+\frac{1}{4} |Y|^2)M^{1/2}}dY dX'\, .
\end{split}
\]
The change of variables $X'-Y/2=v$ and $X'+Y/2=w$ implies
\[
M^{1/2}(|X-X'|^2 + \frac{1}{4}|Y|^2)= \frac{M^{1/2}}{2}(|v-X|^2+|w-X|^2)
\]
and we conclude that the Wigner function and the FBI transform have the relation
\[
\begin{split}
W_{ij}*\phi_{M^{-1/2}}(X,P)&=(2^{1/3}\pi)^{-3N/4}M^{9N/8}\int_{\rset^{3N}} e^{iM^{1/2}(X-v)\cdot P}e^{-|X-v|^2M^{1/2}/2} \Phi_i(v) dv\\
&\qquad\times(2^{1/3}\pi)^{-3N/4}M^{9N/8}
\int_{\rset^{3N}} e^{-iM^{1/2}(X-w)\cdot P}e^{-|X-w|^2M^{1/2}/2} \Phi_j^*(w) dw\\
&=T\Phi_i(X,P) (T\Phi_j(X,P))^*\, .
\end{split}
\]
\end{proof}

\begin{lemma}\label{lem_cn}
{ Assume that $r:\rset^{6N}\to \mathbb C$ is smooth and for each $n\in \mathbb N$ and $|\alpha|\le n$ 
there is a constant $C_n$ such that $\sup_{X\in\rset^{3N}}\int_{\rset^{3N}} |\partial_P^\alpha r(X,P)|dP\le C_{n}$, 
then there is a constant $C$ such that
\begin{equation}\label{est_r1}
\begin{split}
|\int_{\rset^{3N}} \langle \Phi(X), \OP^s[ r]\Phi(X)\rangle dX|
&\le C\sup_{n\le 3N+3}C_n\, .
\end{split}
\end{equation}
If for $|\alpha|\le n$ there is a constant $C_n$ such that
 $\sup_{z\in\rset^{6N}} |\partial_z^\alpha r(z)|\le C_{n}\delta^{\min(0,-n+1)}$,
  then there is a constant $C$ such that
\begin{equation}\label{est_r2}
\begin{split}
|\int_{\rset^{3N}} \langle \Phi(X), \OP^s[ r]\Phi(X)\rangle dX|
&\le C(\sup_{z\in\rset^{6N}} |r(z)| +M^{-1/4}\delta^{-1}N^{1/2}\sup_{n\le CN}C_n)\, .
\end{split}
\end{equation}
}
\end{lemma}
\begin{proof} A proof if this is given e.g. in \cite{zworski} (Theorem 4.21) and presented here for completeness. 

{%
The $s$-symbol $r$ %
has the integral kernel 
\[
K(X,Y)=(\frac{1}{2\pi M^{-1/2}})^{3N}\int_{\rset^{3N}}e^{iM^{1/2}(X-Y)\cdot P} r\big((1-s)X+sY,P\big) dP
\]
so that
\[
\OP^{s}[ r]\Phi(X)=\int_{\rset^{3N}} K(X,Y) \Phi(Y) dY 
\]
and 
\begin{equation}\label{rk_est1}
\begin{split}
\|\OP^{s}[r]\Phi(X)\|_{L^2(\rset^{3N})}^2 &=\int_{\rset^{9N}} \bar K(Z,Y) K(Z,X) \bar \Phi(Y)\Phi(X) dX dY dZ\\
&\le \frac{1}{2}\int_{\rset^{9N}} | K(Z,Y)| |K(Z,X)| (| \Phi(Y)|^2 + |\Phi(X)|^2) dX dY dZ\, .\\
\end{split}
\end{equation}
Here 
\begin{equation}\label{rk_est2}
\begin{split}
&\int_{\rset^{9N}} | K(Z,Y)| |K(Z,X)| | \Phi(Y)|^2  dX dY dZ\\
&\le \sup_{Y\in\rset^{3N}}(\int_{\rset^{6N}} | K(Z,Y)| |K(Z,X)|  dX  dZ ) \underbrace{\int_{\rset^{3N}}|\Phi(Y)|^2 dY}_{=1}\\
&\le \sup_{Y\in\rset^{3N}}\int_{\rset^{3N}} | K(Z,Y)| dZ  
\sup_{Z\in\rset^{3N}}\int_{\rset^{3N}}|K(Z,X)|  dX 
\end{split}
\end{equation}
and similarly for the term with $Y$ and $X$ replaced. 
We have 
\[
\begin{split}
K(X,Y)&=( M^{1/2})^{-3N}  \mathcal F^{-1}\{  r\big((1-s)X+sY, \cdot \big)\}\big((X-Y)M^{1/2}\big) \\
\end{split}
\]
and
\begin{equation*}\label{rk_est}
\begin{split}
\int_{\rset^{3N}} |K(X,Y)| dX &= \int_{\rset^{3N}} | \mathcal F^{-1}\{  r\big((1-s)X'M^{-1/2}+sY, \cdot\big)\}(X') | dX'\\
&= \int_{\rset^{3N}} (1+|X'|^{3N+1} )|\mathcal F^{-1}\{  r\big((1-s)X'M^{-1/2}+sY, \cdot\big)\}(X') | \frac{1}{1+|X'|^{3N+1}}dX'\\
&= \int_{\rset^{3N}} | \mathcal F^{-1}\{ (1+\Delta_{\cdot}^{(3N+1)/2}) r\big((1-s)X'M^{-1/2}+sY, \cdot\big)\}(X') | \frac{1}{1+|X'|^{3N+1}}dX'\\
&\le \|\mathcal F^{-1} \{(1+\Delta_{\cdot}^{(3N+1)/2}) r\big((1-s)X'M^{-1/2}+sY, \cdot\big)\}(X') \|_{L^\infty}
\int_{\rset^{3N}} \frac{1}{1+|X'|^{3N+1}}dX'\\
&\le C \sup_{Y'}\int_{\rset^{3N}} |(1+\Delta_{P}^{(3N+1)/2})r\big(Y', P\big)|dP\, ,
\end{split}
\end{equation*}
where $\Delta_\cdot$ means the Laplacian with respect to the second variable in $r$.
The function $r$ satisfies $\sup_{Y'\in\rset^{3N}}\int_{\rset^{3N}} |\partial_P^\alpha r(Y',P)|dP\le C_{n}$
which together with  \eqref{rk_est1} and \eqref{rk_est2} proves \eqref{est_r1}. %

In the case $\sup_z |\partial_z^{\alpha}r(z)|\le C_n\delta^{-|\alpha|}$
 Theorem 5.1 in \cite{zworski} applied to the symbol $r$ proves \eqref{est_r2}.
 }
\end{proof}

\begin{lemma}\label{moyal_lemma}
 The composition
$\hat C= \hat A \hat B$ of  two Fourier integral operators,  with smooth 
symbols $A(X)$ %
and $B(X,P)$ %
in the Schwartz space,
has the Weyl symbol
\begin{equation}\label{cab}
\begin{split}
C(X,P) &= e^{\frac{i}{2}M^{-1/2} \sum_k( -\partial_{X_k}\partial_{P'_k})} A(X) B(X',P')
\Big|_{(X,P)=(X',P')}\\
&= 
\sum_{n=0}^2 \frac{(iM^{-1/2})^{n}}{2^n n!} ( -\nabla_{X}\cdot \nabla_{P'})^n A(X) B(X',P')
\Big|_{(X,P)=(X',P')}
+ M^{-3/2} r_2\, .\\
\end{split}
\end{equation}
The remainder $r_2$ is smooth and
if $B(X_0,P_0)=\mathbb E[\bar g\circ S_t(\cdot)\, |\, X_0,P_0]*\phi_\eta(X_0,P_0)$, and $A(X_0)= \bar A*\phi_\eta(X_0)$
the remainder satisfies
\begin{equation}\label{lemma5_3_1}
\int_{\rset^{3N}}\langle \Phi,\hat r_2\Phi\rangle dX%
=\mathcal O(e^{\hat Ct}\delta^{-2})\, .
\end{equation}

If $A\in\mathcal C^\infty(\rset^{6N})$ and $B\in\mathcal C^\infty(\rset^{6N})$  are in  the Schwartz space there holds
\begin{equation}\label{weyl_part2}
\begin{split}
C(X,P)  &=
e^{\frac{i}{2}M^{-1/2} \sum_k(\partial_{P_k}\partial_{X'_k} -\partial_{X_k}\partial_{P'_k})} A(X,P) B(X',P')
\Big|_{(X,P)=(X',P')}\\
& = \sum_{n=0}^m \frac{(iM^{-1/2})^{n}}{2^n n!} (\nabla_{P}\cdot\nabla_{X'} -\nabla_{X}\cdot \nabla_{P'})^n A(X,P) B(X',P')
\Big|_{(X,P)=(X',P')} + M^{-(m+1)/2} r_m\, \\
\end{split}
\end{equation}
where for $B(X_0,P_0)=\mathbb E[\bar g\circ S_t(\cdot)\, |\, X_0,P_0]*\phi_\eta(X_0,P_0)$ and $A(X_0,P_0)= \bar A*\phi_\eta(X_0,P_0)$
\begin{equation}\label{r_11}
\int_{\rset^{3N}}\langle \Phi,\hat r_m\Phi\rangle dX
=\mathcal O(e^{\hat Ct}\delta^{-m})\, ,
\end{equation}
for $A(X)=\nabla V(X)-\nabla\lambda_\eta*\phi_\nu(X)$ and 
$B(X_0,P_0)=\partial_P\mathbb E[\bar g\circ S_t(\cdot)\, |\, X_0,P_0]*\phi_\eta(X_0,P_0)$
\begin{equation}\label{rem_2_1_3}
\int_{\rset^{3N}}\langle \Phi,\hat r_1\Phi\rangle dX
=\mathcal O(e^{\hat Ct} \delta^{-4})\, ,
\end{equation}
and for $A(Z)=H_\eta(Z)$ and $B(Z)=2\int_0^1 \tilde g''_{\MD}\big(sH_\eta(Z) + (1-s)E\big)(1-s)ds$ %
\begin{equation}\label{rem_2_2}
\int_{\rset^{3N}}\langle \Phi,\hat r_1\Phi\rangle dX
=\mathcal O(e^{\hat Ct} \delta^{-4})\, .
\end{equation}
\end{lemma}

\begin{proof}
The proof has five steps. The first step uses the definition of the Weyl quantization 
to define an integral kernel for the product $\hat C=\hat A\hat B$ following H\"ormander's work \REF{hormander}.
The next step formulates the Moyal expansion and 
identifies the remainder  as an average of $s-$ quantizations \eqref{s_kvant},
 which proves \eqref{cab}.
 The third step writes the $s$-quantized remainder as a Wigner quantization (with $s=1/2$)
 and estimates the remainder, using that the product of two FBI transformed
 functions is the convolution of the Wigner function with a Gaussian. %
Step 4 proves \eqref{weyl_part2} and \eqref{r_11}.
Step 5 proves \eqref{rem_2_2}.  %

{\bf 1.} To verify \eqref{cab}  we start with the definition of the Weyl operator
\[
\hat C\Phi(X) = \underbrace{(2\pi M^{-1/2})^{-3N}}_{=:\gamma} \int_{\rset^{6N}} C(\frac{X+Y}{2}, P) e^{iM^{1/2}(X-Y)\cdot P} \Phi(Y) dY dP\, .
\]
The $L^2$ inner product
\[
\begin{split}
\int_{\rset^{3N}} \hat C\Phi(X) \Psi^*(X) dX
&=\gamma \int_{\rset^{9N}}  C(\frac{X+Y}{2}, P) e^{iM^{1/2}(X-Y)\cdot P} \Phi(Y) \Psi^*(X) dP dYdX\\
&=\int_{\rset^{6N}}  C_K(X,Y) \Phi(Y) \Psi^*(X) dY dX\\
\end{split}
\]
defines the kernel
\begin{equation}\label{kernel}
C_K(X,Y):= \gamma \int_{\rset^{3N}} C(\frac{X+Y}{2}, P) e^{iM^{1/2}(X-Y)\cdot P}  dP\, ,
\end{equation}
and the inverse Fourier transform implies
\begin{equation}\label{fourierc}
C(U,P)= \int_{\rset^{3N}} C_K(U+ \frac{Z}{2}, U- \frac{Z}{2}) e^{-iM^{1/2}Z\cdot P}  dZ\, .
\end{equation}
Our examples have $A$ independent of $P$, i.e. $A(X,P)=A(X)$,
and the  kernel of the composition $\hat A \hat B$ becomes a multiplication
\[
C_K(X,Y)=\gamma \int_{\rset^{3N}} A(X) B(\frac{X+Y}{2}, P') 
e^{ iM^{1/2}(X-Y)\cdot P'}   dP'\, .
\]
The definition $\hat C=\hat A\hat B$ 
 and  \eqref{fourierc} yields
\begin{equation*} %
\begin{split}
C(U,P) &= (2\pi M^{-1/2})^{-3N}\int_{\rset^{6N}} 
A(U+\frac{Z}{2}) B(U,P+ P')e^{iM^{1/2}P'\cdot Z} dZ   dP'\\
&= (\pi M^{-1/2})^{-3N}\int_{\rset^{6N}} 
A(U+Z') B(U,P+ P')e^{2iM^{1/2}P'\cdot Z'} dZ'   dP'\, .
\end{split}
\end{equation*}

The final step in H\"ormander's derivation uses the  standard Fourier transforms $\mathcal F(f)=\hat f$ of a Schwartz function $f(X,P)$ and of $e^{iM^{1/2}X\cdot P}$ combined with $L^2(\rset^{6N})$ Fourier transform isometry and Taylor expansion function of the exponential. Here we modify this  by identifying the remainder in the Moyal expansion 
for $f(X+X',P+P'):=A(X+X')B(X,P+P')$
as
an $s-$quantization \eqref{s_kvant}.

{\bf  2.} We use in this step the  Fourier transform $\mathcal{F}_{P}$ in the $P$-direction 
defined for $f(X_0,P_0)\in \rset$ by
\[
\mathcal{F}_P \{f\}(X_0, Y):= %
\int_{\rset^{3N}}f(X_0,\xi)e^{-iY\cdot \xi} d\xi\, .
\]
The remainder is based on a Taylor expansion of the exponential function as follows
\begin{equation}\label{exp_def}
\begin{split}
&\int_{\rset^{6N}} f(X+X',P+P') e^{2iM^{1/2}X'\cdot P'} dX'dP'\\
&= %
(\frac{1}{2\pi})^{6N} (\frac{\pi}{\sqrt M})^{3N}
\int_{\rset^{6N}} \mathcal F f(\xi_X,\xi_P) e^{-\frac{i}{2}M^{-1/2}\xi_X\cdot\xi_P} d\xi_Xd\xi_P\\
&= %
(\frac{1}{2\pi})^{6N} (\frac{\pi}{\sqrt M})^{3N}
\int_{\rset^{6N}} \mathcal F f(\xi_X,\xi_P) \Big(\sum_{n=0}^m(\frac{-i\xi_X\cdot\xi_P}{2M^{1/2}})^n \frac{1}{n!} \\
&\qquad +(\frac{-i\xi_X\cdot\xi_P}{2M^{1/2}})^{m+1} \frac{1}{m!} \int_0^1(1-s)^m e^{-\frac{is}{2}M^{-1/2}\xi_X\cdot\xi_P} ds\Big) d\xi_Xd\xi_P\\
&= %
(\frac{\pi}{\sqrt M})^{3N}\Big(\sum_{n=0}^m\frac{1}{n!} (\frac{-i\nabla_X\cdot\nabla_P}{2M^{1/2}})^n  f(X,P) \\
&\qquad + (\frac{1}{2M^{1/2}})^{m+1}
\int_0^1\int_{\rset^{6N}} \mathcal F\big((-i\nabla_X\cdot\nabla_P)^{m+1} f\big)(\xi_X,\xi_P) e^{-\frac{is}{2}M^{-1/2}\xi_X\cdot\xi_P}  \frac{(1-s)^m}{(2\pi)^{6N} m! }d\xi_Xd\xi_Pds\Big)\\
&= (\frac{\pi}{ M^{1/2}})^{3N}  \Big(\sum_{n=0}^m\frac{1}{n!} (\frac{-i\nabla_X\cdot\nabla_P}{2M^{1/2}})^n  f(X,P) \\
&\qquad + (\frac{1}{2M^{1/2}})^{m+1}
\int_0^1\int_{\rset^{3N}} \mathcal F_P\big((-i\nabla_X\cdot\nabla_P)^{m+1} f(X+\cdot, P+\cdot)\big)(\frac{s\xi_P}{2M^{1/2}},\xi_P) 
  \frac{(1-s)^m}{(2\pi)^{3N} m! }d\xi_Pds\Big)\, .\\
\end{split}
\end{equation}
The remainder term can  be written
\[
\begin{split}
&(\frac{1}{2M^{1/2}})^{m+1}
\int_0^1\int_{\rset^{3N}} \mathcal F_P\big((-i\nabla_X\cdot\nabla_P)^{m+1} f(X+\cdot, P+\cdot)\big)(\frac{s\xi_P}{2M^{1/2}},\xi_P) 
  \frac{(1-s)^m}{(2\pi)^{3N} m! }d\xi_Pds\\
  &=
  (\frac{1}{2M^{1/2}})^{m+1}
\int_0^1\int_{\rset^{6N}} (-i\nabla_X\cdot\nabla_P)^{m+1} f(X+\frac{s\xi_P}{2M^{1/2}}, P+P')e^{iP'\cdot \xi_P} 
  \frac{(1-s)^m}{(2\pi)^{3N} m! }dP'd\xi_Pds\\
  &=: (\frac{1}{2M^{1/2}})^{m+1}r_m\\
\end{split}
\]
where $r_m$ is smooth %
since $f$ is a Schwartz function.
The change of variables $P+P'=P''$ and \eqref{kernel} shows that kernel of the remainder becomes 
\[
\begin{split}
r_{m,K}(X,Y) &=\gamma
\int_0^1\int_{\rset^{9N}} (-i\nabla_X\cdot\nabla_P)^{m+1} f(\frac{X+Y}{2}+\frac{s\xi_P}{2M^{1/2}}, P'')\\
&\qquad\times e^{iM^{1/2}(X-Y)\cdot P}
e^{i(P''-P)\cdot \xi_P} 
\frac{(1-s)^m}{(2\pi)^{3N} m! }dP''d\xi_PdP ds\\
&=\gamma
\int_0^1\int_{\rset^{3N}} (-i\nabla_X\cdot\nabla_P)^{m+1} f(\frac{X+Y}{2}+\frac{s(X-Y)}{2}, P'')
e^{iM^{1/2}(X-Y)\cdot P''}
\frac{(1-s)^m}{ m! }dP'' ds\\
&=\gamma
\int_0^1\mathcal F_P\{(-i\nabla_X\cdot\nabla_P)^{m+1} f \}
\big(\frac{X+Y}{2}+\frac{s(X-Y)}{2}, M^{1/2}(X-Y)\big)
\frac{(1-s)^m}{ m! } ds\\
\end{split}
\]
We note that the kernel $A_K^{s_*}(X,Y)$ for any symbols $A$ in the $s_*-$quantization \eqref{s_kvant}
is 
\[\begin{split}
A^{s_*}_K(X,Y)&=\gamma\int_{\rset^{3N}} A\big(X+s_*(Y-X),P\big)e^{iM^{1/2}(X-Y)\cdot P} dP\\
&=\gamma\mathcal F_P\{A\}\big(X+s_*(Y-X),M^{1/2}(X-Y)\big)=\gamma\mathcal F_P\{A\}\big(\frac{X+Y}{2}+\frac{s}{2}(Y-X),M^{1/2}(X-Y)\big)\end{split}
\]
for $s_*=(1-s)/2$. Consequently we  have %
\begin{equation}\label{kvant_est}\begin{split}
\int_{\rset^{3N}} \langle \Phi, \hat r_m\Phi\rangle dX
&=\int_0^1\int_{\rset^{6N}} (\nabla_X\cdot\nabla_P)^{m+1} f(X,P) W^{(s_*)}(X,P)dX dP\frac{(1-s)^m}{ m! } ds\, \\
\end{split}
\end{equation}
and the next step shows %
the final bound
\begin{equation}\label{rms}
\int_{\rset^{3N}} \langle \Phi, \hat r_2\Phi\rangle dX=\mathcal O( e^{\hat CT}\delta^{-2} ).
\end{equation}

{\bf 3.} 

It remains to estimate the remainder  for $r=r_2$ in \eqref{rms}.
Lemma \ref{u_lemma} and Lemma \ref{lemma_time} imply that  we have
\[
\sup_{(X,P)\in\rset^{6N}}\|r^*(X,P)\|_2\le  C \sup_{(X,P)\in\rset^{6N}}\| r_2(X,P)\|_\infty +\mathcal O(\eta)
=\mathcal O(e^{\hat CT}\delta^{-2})\, ,
\]
uniformly in $N$, which proves \eqref{rms}.

{\bf 4.}
The general case
$\hat C= \hat A \hat B$ of  two Fourier integral operators,  with smooth 
symbols $A,B\in \mathcal C^\infty(\rset^{6N})$ in the Schwartz space
yields
\[
\begin{split}
C(U,P)&= \gamma^2 \int_{\rset^{12N}} 
A(\frac{U+Z}{2} + \frac{T}{4}, P'') B(\frac{U+ Z}{2} - \frac{T}{4}, P') e^{iM^{1/2}F}dP'' dP' dZ dT\\
& =
\gamma^2 \int_{\rset^{12N}} 
A(\frac{U+Z}{2} + \frac{T}{4}, P'') B(\frac{U+ Z}{2} - \frac{T}{4}, P') e^{iM^{1/2}F}dP'' dP' dZ dT
\, ,
\end{split}
\]
where 
\[
\begin{split}
F&:=(U-Z+ \frac{T}{2})\cdot P'' + (Z-U + \frac{T}{2})\cdot P' -  T\cdot P\\
&=(U-Z+ \frac{T}{2})\cdot (P''-P) + (Z-U + \frac{T}{2})\cdot (P'-P)\, .
\end{split}
\]
The change of variables $(P''-P, P'-P, (Z-U+ \frac{T}{2})/2, (Z-U - \frac{T}{2})/2)$ replacing $(P'',P', Z, T)$ has
the Jacobian $2^{6N}$ and implies
\begin{equation*} %
C(U,P)= (\pi M^{-1/2})^{-6N}\int_{\rset^{12N}} 
A(U+Z, P+P'') B(U+T,P+ P')e^{iM^{1/2}(P'\cdot Z-T\cdot P'')} dZ dP'' dT dP'\, .
\end{equation*}
The same expansion as in \eqref{exp_def} shows that
\begin{equation*}%
\begin{split}
C(X,P) &= e^{\frac{i}{2}M^{-1/2} \sum_k(\partial_{P_k}\partial_{X'_k} -\partial_{X_k}\partial_{P'_k})} A(X,P) B(X',P')
\Big|_{(X,P)=(X',P')}\\
&= \sum_{n=0}^m \frac{(iM^{-1/2})^{n}}{2^n n!} (\nabla_{P}\cdot\nabla_{X'} -\nabla_{X}\cdot \nabla_{P'})^n A(X,P) B(X',P')
\Big|_{(X,P)=(X',P')} + M^{-(m+1)/2} r_m\\
\end{split}
\end{equation*}
based on the representation %
\[
\begin{split}
&\int_{\rset^{6N}} f(X+X',P+P', X+X'',P+P'') e^{2iM^{1/2}(X'\cdot P'-X''\cdot P'')} dX'dP' dX''dP''\\
&= (\frac{\pi}{ M^{1/2}})^{3N}  \Big(\sum_{n=0}^m\frac{1}{n!} \big(i(\nabla_{X''}\cdot\nabla_{P'}-\nabla_{X'}\cdot\nabla_{P''})\big)^{m+1}  f(X+X',P+P',X+X'',P+P'')\Big|_{
{\tiny \begin{array}{c}
X'=X''=0\\
P'=P''=0
\end{array}} }\\
&\qquad + (\frac{1}{2M^{1/2}})^{m+1}r_m\Big)\COMMA\\
&r_m := 
\int_0^1\int_{\rset^{12N}} \big(i(\nabla_{X''}\cdot\nabla_{P'}-\nabla_{X'}\cdot\nabla_{P''}\big)^{m+1} 
f(X+\frac{s\xi_{P'}}{2M^{1/2}}, P+P', X-\frac{s\xi_{P''}}{2M^{1/2}}, P+P'')\\
&\qquad\qquad\times e^{i(P'\cdot \xi_{P'}-P''\cdot \xi_{P''}} 
  \frac{(1-s)^m}{(2\pi)^{6N} m! }dP'd\xi_{P'}dP''d\xi_{P''}ds\, .\\
  \end{split}
  \]
    The kernel for $m=1$ becomes
  \[
\begin{split}
&r_{1,K}(X,Y)\\
&=\gamma
\int_0^1 \mathcal F_{P'P''}\{\big(i(\nabla_{X''}\cdot\nabla_{P'}-\nabla_{X'}\cdot\nabla_{P''})\big)^{2} 
f\}\\
&\qquad \Big(\frac{X+Y}{2}-\frac{s(X-Y)}{2M^{1/2}}, M^{1/2}(X-Y),\frac{X+Y}{2}-\frac{s(X-Y)}{2M^{1/2}}, M^{1/2}(X-Y)\Big)
(1-s)ds\, .\\
\end{split}
\]
and as in \eqref{kvant_est}
we obtain
\[\begin{split}
&\int_{\rset^{3N}} \langle \Phi, \hat r_1\Phi\rangle dX\\
&=\int_0^1\int_{\rset^{6N}} (\nabla_{X''}\cdot\nabla_{P'}-\nabla_{X'}\cdot\nabla_{P''})\big)^{2}f(X+X',P+P',X+X'',P+P'')\big|_{X'=X''=P'=P''=0}\\
&\qquad\times W^{(s_*)}(X,P)dX dP(1-s) ds\\
\end{split}
\]
We apply this to the special case $C(X,P)=H_\eta(X,P)B(X,P)$  
where $B(Z)= \mathbb E[\bar g\circ S_t(\cdot)\, |\, X_0,P_0]*\phi_\eta(Z)$ is a Schwartz function. %
Then we have  $f(X',P',X'',P'')=H_\eta(X',P')B(X'',P'')$
and for $B=\partial_{P} \mathbb E[\bar g\circ S_T(\cdot)\, |\, X_0,P_0]*\phi_\eta(X,P) $ Lemma \ref{u_lemma} yields as in Step 3 
the  bound
\begin{equation*}\label{shurll}
\int_{\rset^{3N}} \langle \Phi, \hat r_1\Phi\rangle dX=\mathcal O( e^{\hat CT}\delta^{-4} ).
\end{equation*}

 {\bf 5.} %
 If $A=H_\eta$ and $B(Z)=2\int_0^1 \tilde g''_{\MD}\big(sH_\eta(Z) + (1-s)E\big)(1-s)ds$
 the symbol $D:=AB$ is a smooth product of two  convolutions, by \eqref{gmd_bar}, 
 and $D(Z)$ tends to zero fast for large $|Z|$.
 We obtain %
 as in Step 3
 \[
 \int_{\rset^{3N}} \langle \Phi, \hat r_1\Phi\rangle dX=\mathcal O( e^{\hat CT}\delta^{-4} ).
 \]

 \end{proof}

\section{Appendix: Fourier integral WKB states including caustics}\label{sec:caustics}

\subsection{A preparatory example with the simplest caustic}
\begin{figure}[htbp]
\centering
\includegraphics[height=5cm]{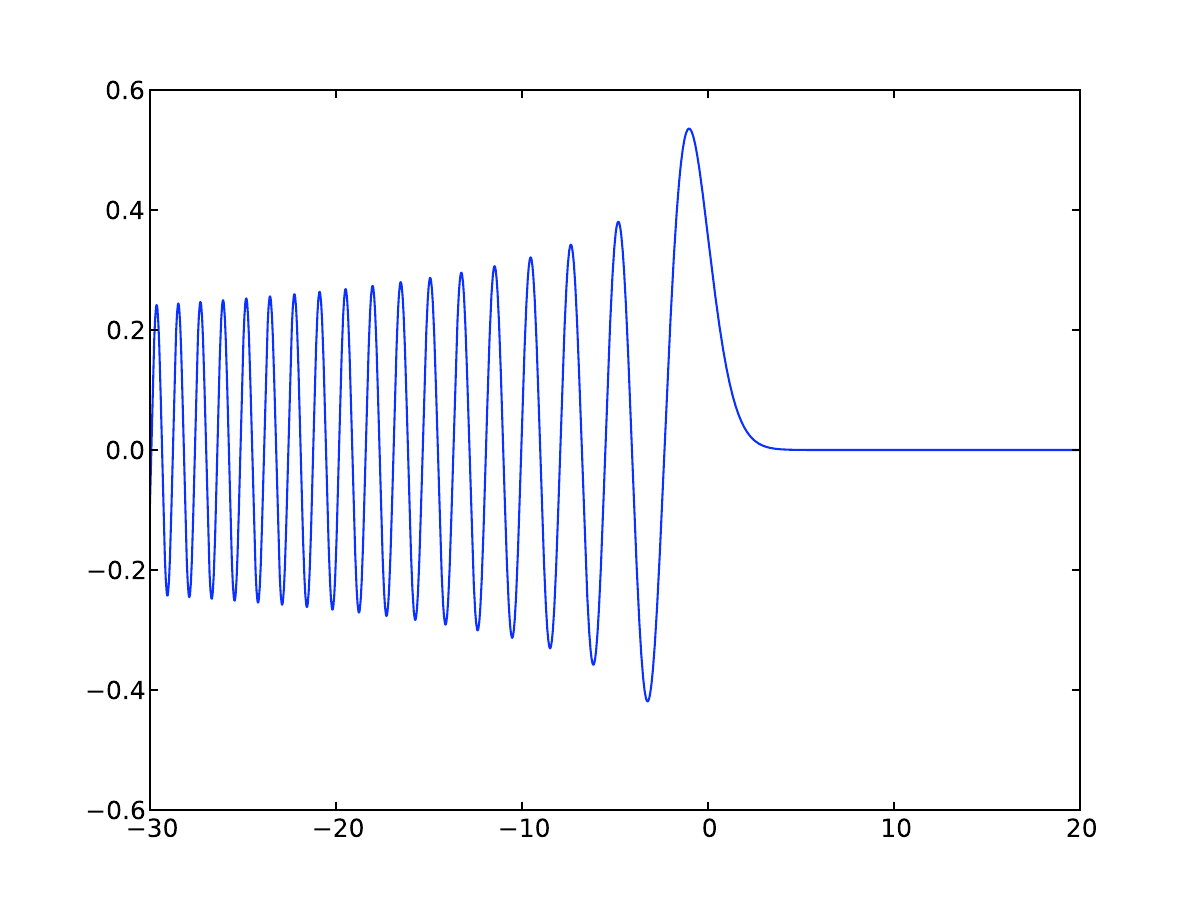}
\caption{The Airy function.}
\label{airy_fig}
\end{figure}

As an example of a caustic, we study first the simplest example of a fold caustic based on the Airy function
$\AIRY:\rset\rightarrow\rset$ which solves
\begin{equation}\label{eq:Airy_ODE}
  -\partial_{xx}{\AIRY}(x)  +x \AIRY(x)=0\PERIOD
\end{equation}
The scaled Airy function
\[
  u(x)= C\,\AIRY(M^{1/3} x)
\]
solves the Schr\"odinger equation
\begin{equation}\label{eq:airy_m}
-\frac{1}{M}\partial_{xx}{u}(x) + x u(x)=0\COMMA
\end{equation}
for any constant $C$. In our context an important property of the Airy function is the fact
that it is the inverse Fourier transform of the function
\[
 \hat\AIRY(p) = \sqrt{\frac{2}{\pi}}\EXP{\Iunit p^3/3}\COMMA
\]
i.e.,
\begin{equation}\label{airy_fourier}
  \AIRY(x)= \frac{1}{\pi} \int_\rset \EXP{\Iunit (xp+p^3/3)}\, dp\PERIOD
\end{equation}
In the next section, we will consider a general Schr\"odinger equation and determine
a WKB Fourier integral corresponding to \eqref{airy_fourier} for the
Airy function; as an introduction to the general case we show how the derive \eqref{airy_fourier}:  by taking the Fourier transform of the ordinary differential equation~\eqref{eq:Airy_ODE}
\begin{equation}
    0 = \int_\rset \left(-\partial_{xx}  +x\right){\AIRY}(x) \EXP{-\Iunit xp}\, dx
      = (p^2+\Iunit \partial_p) \hat\AIRY(p)\COMMA %
\end{equation}
we obtain an ordinary differential equation for the Fourier transform $ \hat\AIRY(p)$ %
with the solution $ \HATAIRY(p)= C \EXP{\Iunit p^3}$, for any constant $C$.
Then, by differentiation, it is clear that the scaled Airy function $u$ solves~\eqref{eq:airy_m}.
Furthermore, the stationary phase method, cf. Section~\ref{stat_phase_sec}, shows that to
the leading order $u$ is approximated by 
\begin{equation}\label{airy_stat_fas}
u(x)\simeq C\left(-x M^{1/3}\right)^{-1/4} \cos \big(M^{1/2} (-x)^{3/2} -\pi/4\big)\COMMA\;\;
  \mbox{ for } x < 0\COMMA
\end{equation}
and $u(x)\simeq 0$ to any order (i.e., $\BIGO(M^{-K})$ for any positive $K$) when $x>0$. 
The behaviour of the Airy function is illustrated in Figure~\ref{airy_fig}.

\subsubsection{Molecular dynamics for the Airy function}
The eikonal equation  corresponding to \eqref{eq:airy_m} is 
\[
  p^2+x=0
\]
with solutions for $x\le 0$, which leads to the phase
\begin{equation}\label{legendre_airy}
  p=\theta'(x)=\pm (-x)^{1/2}\COMMA\;\mbox{and}\;\;\theta(x)= \mp \frac{2}{3} (-x)^{3/2}\PERIOD
\end{equation}
We have the Legendre transform
\[
  \LFT{\theta}(p)=\min_y\big( y p -\theta(y)\big)
\]
with the solution $p=\theta'(y)$. Consequently,  we obtain by  \eqref{legendre_airy} that $x=y$
and
\[
   \LFT{\theta}(p)=xp-\theta(x)=- p^2 p +\frac{2}{3} p^3 = -\frac{p^3}{3}\PERIOD
\]
We note that this solution is also obtained from the eikonal equation
\[
  p^2 + \partial_p\LFT{\theta}(p)=0\COMMA
\]
which is solved by 
\[
   \LFT{\theta}(p)=-p^3/3\PERIOD
\]
Thus we recover the relation for the Legendre transform 
$-xp +\LFT{\theta}(p)=-\theta(x)$.

\subsection{A general Fourier integral ansatz}\label{sec:caustic_general}
In order to treat a more general case with a caustic of the dimension $d$ we use  the Fourier integral ansatz
\begin{equation}\label{caustic_ansats}
  \Phi(X)= \int_{\rset^d} \hpsi(X) 
             \EXP{-\Iunit M^{1/2}\PHASE}\, d\CHECKP
\end{equation}
where
\[
\begin{split}
\phi(X)&\in \mathbb C^J\COMMA\\
  X &= (\HATX, \CHECKX)\in \rset^{3N}\COMMA\;\;\; P = (\HATP,\CHECKP)\in \rset^{3N}\\
  \EPROD{\CHECKX}{\CHECKP} &=\sum_{j=1}^d\CHECKX^j\CHECKP^j\COMMA\;\;\;
        \EPROD{\HATX}{\HATP}      =\sum_{j=d+1}^N\HATX^j\HATP^j\\
  \PHASE & = \EPROD{\CHECKX}{\CHECKP} -\LFT{\theta}(\HATX,\CHECKP)\COMMA
\end{split}
\]
based on the Legendre transform relationships (cf. \cite{evans})
\[
\begin{split}
  \LFT{\theta}(\HATX,\CHECKP)&= \min_{\CHECKX}\left(\EPROD{\CHECKX}{\CHECKP} - \theta(\HATX, \CHECKX)\right)\COMMA\\
  {\theta}(\HATX,\CHECKX)&= \min_{\CHECKP}\left(\EPROD{\CHECKX}{\CHECKP} - \theta^*(\HATX, \CHECKP)\right)\, .
  \end{split}
\]
If the function $\LFT{\theta}(\HATX,\CHECKP)$ is not defined for all $\CHECKP\in\rset^d$,
but only for $\CHECKP\in \NEIGH\subset \rset^d$ we replace the integral over $\rset^d$ by integration over $\NEIGH$
using a smooth cut-off function $\chi(\CHECKP)$.
The cut-off function is zero outside $\NEIGH$ and equal to one in a large part
of the interior of $\NEIGH$. %
The  ansatz \eqref{caustic_ansats} is inspired by Maslov's work \REF{maslov}, although
it is not the same since our amplitude function $\hpsi$
depends on $(\HATX,\CHECKX)$ but not on $\CHECKP$.

\subsubsection{Making the ansatz for a Schr\"odinger solution}
In this section we construct a solution to the Schr\"odinger equation from the ansatz \eqref{caustic_ansats}.
The constructed solution will be an {\it actual} solution and not only an asymptotic solution as in 
\REF{maslov}.
We consider first the case when the integration is over $\rset^d$ and
then conclude in the Remark \ref{chi_rem} that the cut-off function $\chi(\CHECKP)$ can be included in all
integrals without changing the property of the Fourier integral ansatz being a solution
in the $\CHECKX$-domain  where $\CHECKX=\GRAD_{\CHECKP}\LFT{\theta}(\HATX,\CHECKP)$ for some $\CHECKP$ satisfying
$\chi(\CHECKP)=1$.

The requirement to be a solution means that there should hold
\begin{equation}\label{H_ekv_caustic}
\begin{split}
0 &=(\HOPER-E)\Phi \\
  &=\int_{\rset^d} \left(\frac{1}{2}|\nabla_{\HATX}\LFT{\theta}(\HATX,\CHECKP)|^2 
       +\frac{1}{2}|\CHECKP|^2 +V_0(X) -E\right)\hpsi(X) \EXP{-\Iunit M^{1/2}\PHASE}\, d\CHECKP\\
  &\quad {}-\int_{\rset^d}\left(\Iunit M^{-1/2}(\EPROD{\nabla_{\HATX}\hpsi}{\nabla_{\HATX}\LFT{\theta}} 
                - \EPROD{\nabla_{\CHECKX}\hpsi}{\CHECKP}+\frac{1}{2} \hpsi\LAP_{\HATX}\LFT{\theta})
                -(\VOPER-V_0)\hpsi +\frac{1}{2M}\LAP_{X}\hpsi\right)
            \EXP{-\Iunit M^{1/2}\PHASE}\, d\CHECKP .
\end{split}
\end{equation}
Comparing this expression to the previously discussed case of a single WKB-mode we see that the zero order term is
now $\LAP_{\HATX}\LFT{\theta}$ instead of $\LAP_{X}\theta$
and that we have $-\EPROD{\nabla_{\CHECKX}\hpsi}{\CHECKP}$ instead of 
$\EPROD{\nabla_{\CHECKX}\hpsi}{\nabla_{\CHECKX}\theta}$.
However, the main difference is that the first integral is 
not zero (only the leading order term of its
stationary phase expansion is zero, cf. \eqref{caustic_expansion}). Therefore, the first integral contributes to the second integral. The
goal is now to determine a function $\SFUN(\HATX,\CHECKX, \CHECKP)$ satisfying
\begin{equation}\label{v-v_int}
\begin{split}
   \int_{\rset^d} & \left(\frac{1}{2}|\nabla_{\HATX}\LFT{\theta}|^2 +\frac{1}{2}|\CHECKP|^2
        +V_0(X) -E\right)\EXP{-\Iunit M^{1/2}\PHASE}\, d\CHECKP\\
        & =\Iunit M^{-1/2}\int_{\rset^d}  \SFUN(\HATX,\CHECKX, \CHECKP)\,\EXP{-\Iunit M^{1/2}\PHASE}\, d\CHECKP\COMMA
\end{split}
\end{equation}
and verify that it is bounded.
\begin{lemma}\label{F_lem}
There holds
$F=F_0+F_1$
where
\[
\begin{split}
F_0 &=\frac{1}{2}\sum_{i,j}  \partial_{\CHECKX^i\CHECKX^j } V_0\left(\DPLFT{\theta}(\CHECKP)\right) 
                    \partial_{\CHECKP^j\CHECKP^i}\LFT{\theta}(\CHECKP)\COMMA \\
F_1 &= \Iunit M^{-1/2} \int_0^1\int_0^1\int_{\rset^d}  \sum_{i,j,k} t(1-t)
     \partial_{\CHECKP^k} \left[\partial_{\CHECKX^i\CHECKX^j\CHECKX^k} 
          V_0\left(\DPLFT{\theta}(\CHECKP)  + s\, t\,\DTHETA\right)
           \partial_{\CHECKP_j\CHECKP_i}\DPLFT{\theta}(\CHECKP)\right]\,dt\,ds\PERIOD\\
\end{split}
\]
\end{lemma}
\begin{proof}
 The function
$\LFT{\theta}(\HATX,\CHECKP)$ is defined as a solution to
the Hamilton-Jacobi (eikonal) equation
\begin{equation}\label{caustic_hj}
    \frac{1}{2}|\nabla_{\HATX}\LFT{\theta}(\HATX,\CHECKP)|^2 +\frac{1}{2}|\CHECKP|^2 
       + V_0\left(\HATX,\nabla_{\CHECKP}\LFT{\theta} (\HATX,\CHECKP)\right) - E = 0
\end{equation}
for all $(\HATX,\CHECKP)$.
Consequently, the integral on the left hand side of \eqref{v-v_int} is 
\[
    \int_{\rset^d} \left(V_0(\HATX,\CHECKX) -V_0(\HATX,\nabla_{\CHECKP}\LFT{\theta}(\HATX,\CHECKP)\right)
    \EXP{-\Iunit M^{1/2}\left(\EPROD{\CHECKX}{\CHECKP} -\LFT{\theta}(\HATX,\CHECKP)\right)}\, d\CHECKP\PERIOD
\]
Let  $\CHECKP_0(\CHECKX)$ be any solution to the stationary phase equation
$\CHECKX=\nabla_{\CHECKP}\LFT{\theta}(\HATX,\CHECKP_0)$
and introduce the notation 
\[
 \PHASED  :=\EPROD{\GRAD_{\CHECKP}\LFT{\theta}(\HATX,\CHECKP_0)}{\CHECKP} 
           -\LFT{\theta}(\HATX,\CHECKP)\PERIOD
\]
Then by writing a difference as $V(y_1)-V(y_2)=\int_0^1\partial_y V(y_2+t(y_1-y_2))dt\cdot(y_1-y_2)$, 
identifying a derivative $\partial_{\CHECKP_i} $ and integrating by parts
the integral can be written
\[
\begin{split}
    &\int_{\rset^d} \left(V_0(\HATX,\nabla_{\CHECKP}\LFT{\theta}(\HATX,\CHECKP_0))
           -V_0(\HATX,\nabla_{\CHECKP}\LFT{\theta}(\HATX,\CHECKP)\right)
             \EXP{-\Iunit M^{1/2}\PHASED}\, d\CHECKP\\
    &=  \int_0^1\int_{\rset^d}  \sum_i\partial_{\CHECKX^i} V_0\left(\DPLFT{\theta}(\CHECKP) 
                    +t\,\left[\DPLFT{\theta}(\CHECKP_0)-\DPLFT{\theta}(\CHECKP)\right]\right)\times\\
    & \qquad\times        \left(\partial_{\CHECKP^i}\LFT{\theta}(\CHECKP_0)-\partial_{\CHECKP^i}\LFT{\theta}(\CHECKP)\right) 
        \EXP{-\Iunit M^{1/2}\PHASED}\,d\CHECKP\, dt \\ 
    &= -\Iunit M^{-1/2} \int_0^1\int_{\rset^d}  \sum_i\partial_{\CHECKX^i} V_0
           \left(\DPLFT{\theta}(\CHECKP) +t\,\left[\DPLFT{\theta}(\CHECKP_0)-\DPLFT{\theta}(\CHECKP)\right]\right)
           \partial_{\CHECKP_i} \EXP{-\Iunit M^{1/2}\PHASED} \,d\CHECKP\, dt \\ 
   &= \Iunit M^{-1/2} \int_0^1\int_{\rset^d}  \sum_i\partial_{\CHECKP_i} \partial_{\CHECKX^i} 
                V_0\left(\DPLFT{\theta}(\CHECKP) +t\,\left[\DPLFT{\theta}(\CHECKP_0)-\DPLFT{\theta}(\CHECKP)\right]\right)
            \EXP{-\Iunit M^{1/2}\PHASED}\,d\CHECKP\, dt\PERIOD
\end{split}
\]
Therefore the  leading order term  in $\SFUN=: \SFUN_0 + \SFUN_1$ is
\[
\begin{split}
  \SFUN_0&:=\int_0^1\sum_{i,j}  (1-t) \partial_{\CHECKX^i\CHECKX^j } V_0\left(\DPLFT{\theta}(\CHECKP)\right)
                                 \partial_{\CHECKP^j\CHECKP^i}\LFT{\theta}(\CHECKP)\, dt\\
     &=\frac{1}{2}\sum_{i,j}  \partial_{\CHECKX^i\CHECKX^j } V_0\left(\DPLFT{\theta}(\CHECKP)\right) 
                    \partial_{\CHECKP^j\CHECKP^i}\LFT{\theta}(\CHECKP)\PERIOD
\end{split}
\]
Denoting $\DTHETA := \DPLFT{\theta}(\CHECKP_0)-\DPLFT{\theta}(\CHECKP)$  the remainder becomes
\[
\begin{split}
    & -\Iunit M^{-1/2} \int_0^1\int_{\rset^d}  \sum_{i,j} \left[
           \partial_{\CHECKX^i\CHECKX^j} V_0\left(\DPLFT{\theta}(\CHECKP) \right)
         - \partial_{\CHECKX^i\CHECKX^j} V_0\left(\DPLFT{\theta}(\CHECKP)  + t\,\DTHETA\right) \right]\\
    &\qquad \times (1-t)\partial_{\CHECKP^j\CHECKP^i}\LFT{\theta}(\CHECKP) \,\EXP{-\Iunit M^{1/2}\PHASED} \,d\CHECKP\, dt\\
    &= \Iunit M^{-1/2} \int_0^1\int_0^1\int_{\rset^d}  \sum_{i,j,k} t(1-t)
          \partial_{\CHECKX^i\CHECKX^j\CHECKX^k} V_0\left(\DPLFT{\theta}(\CHECKP)  +  s\, t\,\DTHETA\right) 
           \partial_{\CHECKP^j\CHECKP^i}\LFT{\theta}(\CHECKP) \\
    &\qquad \times \left(\partial_{\CHECKP^k}\LFT{\theta}(\CHECKP_0)-\partial_{\CHECKP^k}\LFT{\theta}(\CHECKP)\right)
      \EXP{-\Iunit M^{1/2}\PHASED} \,d\CHECKP\,dt\, ds\\
    &= - \frac{1}{M} \int_0^1\int_0^1\int_{\rset^d}  \sum_{i,j,k} t(1-t)
           \partial_{\CHECKP^k} \left[\partial_{\CHECKX^i\CHECKX^j\CHECKX^k} 
             V_0\left(\DPLFT{\theta}(\CHECKP)  + s\, t\, \DTHETA\right)   
             \partial_{\CHECKP^j\CHECKP^i}\LFT{\theta}(\CHECKP)\right] \\
    &\qquad\times \EXP{-\Iunit M^{1/2}\PHASED}\,d\CHECKP\,dt\,ds\COMMA
\end{split}
\]
hence the function $\SFUN_1$ is purely imaginary and  small
\[
\begin{split}
  &\SFUN_1= \Iunit M^{-1/2} \int_0^1\int_0^1\int_{\rset^d}  \sum_{i,j,k} t(1-t)
     \partial_{\CHECKP^k} \left[\partial_{\CHECKX^i\CHECKX^j\CHECKX^k} 
          V_0\left(\DPLFT{\theta}(\CHECKP)  + s\, t\,\DTHETA\right)
           \partial_{\CHECKP_j\CHECKP_i}\DPLFT{\theta}(\CHECKP)\right]\,dt\,ds\COMMA
\end{split}
\]
and
\begin{equation}\label{Re_s}
2\REAL{\SFUN}=\sum_{i,j}  \partial_{\CHECKX^i\CHECKX^j } V_0\left(\DPLFT{\theta}(\CHECKP)\right) 
            \partial_{\CHECKP^j\CHECKP^i}\LFT{\theta}(\CHECKP)\PERIOD
\end{equation}
\end{proof}

The eikonal equation \eqref{caustic_hj} and the requirement that
$(\HOPER-E)\Phi =0$ in \eqref{H_ekv_caustic} then imply that
\begin{equation}\label{before_hj}
\begin{split}
   0&= \int_{\rset^d}\left[\Iunit M^{-1/2}\left(\EPROD{\nabla_{\HATX}\hpsi}{\nabla_{\HATX}\LFT{\theta}} 
         - \EPROD{\nabla_{\CHECKX}\hpsi}{\CHECKP}
         +\frac{1}{2} \hpsi \left( \LAP_{\HATX}\LFT{\theta} - 2 \SFUN(X,\CHECKP) \right)\right)\right.\\
    &\qquad \left. -(\VOPER-V_0)\hpsi +\frac{1}{2M}\LAP_{X}\hpsi\right]
            \EXP{-\Iunit M^{1/2}\PHASE}\, d\CHECKP\PERIOD
\end{split}
\end{equation}
The Hamilton-Jacobi eikonal equation \eqref{caustic_hj}, in the primal variable $(\HATX,\CHECKP)$
with the corresponding dual variable $(\hat P, \check X)$,
can be 
solved locally by  the characteristics
\begin{equation}\label{caustic_char}
\begin{split}
   \dot{\HATX} & = \HATP\\
   \dot{\HATP} &= -\GRAD_{\HATX} V_0(\HATX,\CHECKX)\\
   \dot{\CHECKX} & = -\CHECKP\\
   \dot{\CHECKP} &= \GRAD_{\CHECKX} V_0(\HATX,\CHECKX)\COMMA
\end{split}
\end{equation}
using the definition
\[
\begin{split}
    \GRAD_{\HATX}\LFT{\theta}(\HATX,\CHECKP)   &= \HATP         \\
    \GRAD_{\CHECKP}\LFT{\theta}(\HATX,\CHECKP) &= \CHECKX\PERIOD
\end{split}
\]
The characteristics give 
\[
   \frac{d}{dt} \hpsi = \EPROD{\GRAD_{\HATX}\hpsi}{\GRAD_{\HATX}\LFT{\theta}} - 
                       \EPROD{\GRAD_{\CHECKX}\hpsi}{\CHECKP}\COMMA
\]
so that the Schr\"odinger transport equation becomes, as in  \eqref{schrod_first},
\begin{equation}\label{caustic_transport}
    \Iunit M^{-1/2} \left(\dot{\hpsi} +  \hpsi \frac{\dot G}{G}\right) 
                = (\VOPER - V_0) \hpsi - \frac{1}{2M} \LAP_{X}\hpsi
\end{equation}
and for $\tpsi=G\phi$
\begin{equation}\label{caustic_transport_psi}
    \Iunit M^{-1/2} \dot{\tpsi}= (\VOPER - V_0) \tpsi - \frac{G}{2M} \LAP_{X}\frac{\tpsi}{G}
\end{equation}
with  the complex valued weight function $G$ defined by
\begin{equation}\label{g_2}
  \frac{d}{dt} \log G_t= \frac{1}{2}\LAP_{\HATX}\LFT{\theta}(\hat X_t,\CHECKP_t) - \SFUN(\hat X_t,\CHECKP_t)\PERIOD
\end{equation}
This transport equation is of the same form as the transport equation  for a single WKB-mode, with
a modification of the weight function $G$.

Differentiation of the second equation in the 
Hamiltonian system \eqref{caustic_char} implies that the first variation 
$\partial \CHECKP_t/\partial \CHECKX_0$
satisfies
\[
  \frac{d}{dt}\left(\frac{\partial\CHECKP_t^i}{\partial \CHECKX_0}\right) = 
      \sum_{j,k}\partial_{\CHECKX^i\CHECKX^j} V_0(\HATX,\CHECKX_t) 
           \partial_{\CHECKP^j\CHECKP^k}\LFT{\theta}(\CHECKP) \frac{\partial\CHECKP_t^k}{\partial \CHECKX_0} \COMMA
\]
which by the Liouville formula \eqref{liouville_form} and the equality 
\[
2\REAL{\SFUN}=\sum_{i,j} \partial_{\check X^i\check X^j}V_0\partial_{\CHECKP^j\CHECKP^i}\LFT{\theta}
=\TRACE(\sum_j\partial_{\CHECKX^i\CHECKX^j}V_0\partial_{\CHECKP^j\CHECKP^k}\LFT{\theta})\]
in \eqref{Re_s} yields the relation, 
\begin{equation}\label{re_s}
   \EXP{-2\int_0^t\REAL{\SFUN}\, dt'} = C\, \left|\DET \frac{\partial \CHECKP_t}{\partial \CHECKX_0}\right|\COMMA
\end{equation}
for the constant $C:=|\DET \frac{\partial \CHECKX_0}{\partial \CHECKP_0}|$.  

\begin{remark}\label{chi_rem} {\rm
  The conclusion in this section holds also when all integrals over $d\CHECKP$ in $\rset^d$
  are replaced by integrals with the measure $\chi(\CHECKP)\, d\CHECKP$. Then there holds
  $2\REAL{\SFUN}=\sum_{ij} \partial_{\CHECKX^i\CHECKX^j}\VOPER \partial_{\CHECKP^i}(\chi \partial_{\CHECKP^j}\LFT{\theta})$. }
\end{remark}

\subsubsection{Liouville's formula}\label{liouville}
Here we state Liouville's formula
\begin{equation}\label{liouville_form}
   \frac{G^2_0}{G_t^2}=\EXP{-\int_0^t \TRACE\left(\GRADX P(X_t)\right) \,dt}= 
   \left|\DET \frac{\partial (X_0)}{\partial (X_t)}\right| \COMMA
\end{equation}
given in \REF{maslov}.

\section{Appendix: The stationary phase expansion}\label{stat_phase_sec}
Consider the phase function
$\EPROD{\CHECKX}{\CHECKP} - \LFT\theta(\HATX,\CHECKP)$
and let  $\CHECKP_0(\HATX)$ be any solution to the stationary phase equation
$\CHECKX=\GRAD_{\CHECKP}\LFT \theta(\HATX,\CHECKP_0)$. We rewrite the phase function
\[
   \EPROD{\CHECKX}{\CHECKP} - \LFT\theta(\CHECKX,\CHECKP)= 
   \underbrace{\EPROD{\CHECKX}{\CHECKP_0} - \LFT\theta(\CHECKX,\CHECKP_0)}_{=\theta(\HATX,\CHECKX)} +
                         \EPROD{(\CHECKP-\CHECKP_0)}{\int_0^1(1-t)
                         \partial_{PP}\LFT\theta\left(\CHECKP_0+ t[\CHECKP-\CHECKP_0]\right)\, dt}\,
                         [\CHECKP-\CHECKP_0]\PERIOD
\]
The relation
\[
 \frac{1}{2} \EPROD{Y}{\partial_{PP}\bar\theta(\CHECKP_0) Y}= 
    \EPROD{(\CHECKP-\CHECKP_0)}{\int_0^1 (1-t)
      \partial_{PP}\bar\theta\left(\CHECKP_0 + t[\CHECKP-\CHECKP_0]\right)\, dt}\, [\CHECKP-\CHECKP_0]
 \]
 defines the function $Y(\CHECKP)$, and its inverse $\CHECKP(Y)$, so that the phase is a quadratic function in $Y$.
 The stationary phase expansion of an integral takes the form, see \REF{zworski},
 \begin{equation}\label{caustic_expansion}
   \begin{split}
     &\int_{\rset^d} w(\CHECKP)\, \EXP{-\Iunit M^{1/2}(\EPROD{\CHECKX}{\CHECKP} 
                     -\LFT\theta(\HATX,\CHECKP))}\, d\CHECKP\\
     &\sim \sum_{\{\CHECKP_0\ : \ \nabla_{\check P}\LFT\theta(\CHECKP_0) =\CHECKX\}} (2\pi M^{-1/2})^{d/2} 
            \left| \DET \frac{\partial(\CHECKP_0)}{\partial(\CHECKX)}\right|^{1/2} \,
            \EXP{\Iunit\frac{\pi}{4}\SGN(\partial_{PP}\LFT\theta(\CHECKP_0))}\,
             \EXP{-\Iunit M^{1/2}\theta(\HATX,\CHECKX;\check P_0)}\\
     &\qquad \times \sum_{k=0}^\infty \frac{M^{-k/2}}{k!} 
             \left(\sum_{l,j} \Iunit (\partial_{P^l P^j}\LFT\theta)^{-1}(\CHECKP_0) 
                              \partial_{Y^l Y^j}\right)^k
                  \left(w(\CHECKP(Y))\,\left|\DET \frac{\partial(\CHECKP)}{\partial(Y)}\right|\right)\Big|_{Y=0}\PERIOD
   \end{split}
 \end{equation}

\section{Appendix: An alternative motivation for assumption \eqref{g_bar}}\label{motiv_indep}
\begin{figure}[htbp]
  \includegraphics[width=0.32\textwidth,height=0.5\textwidth]{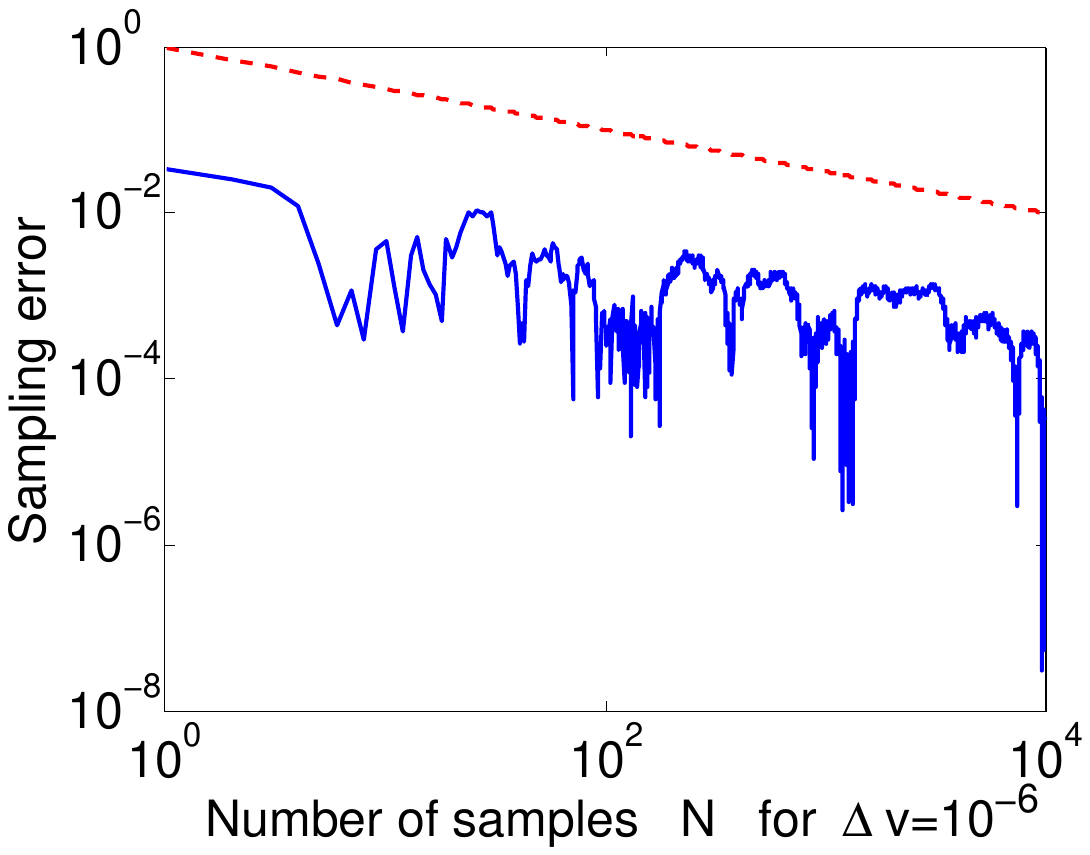}
  \includegraphics[width=0.32\textwidth,height=0.5\textwidth]{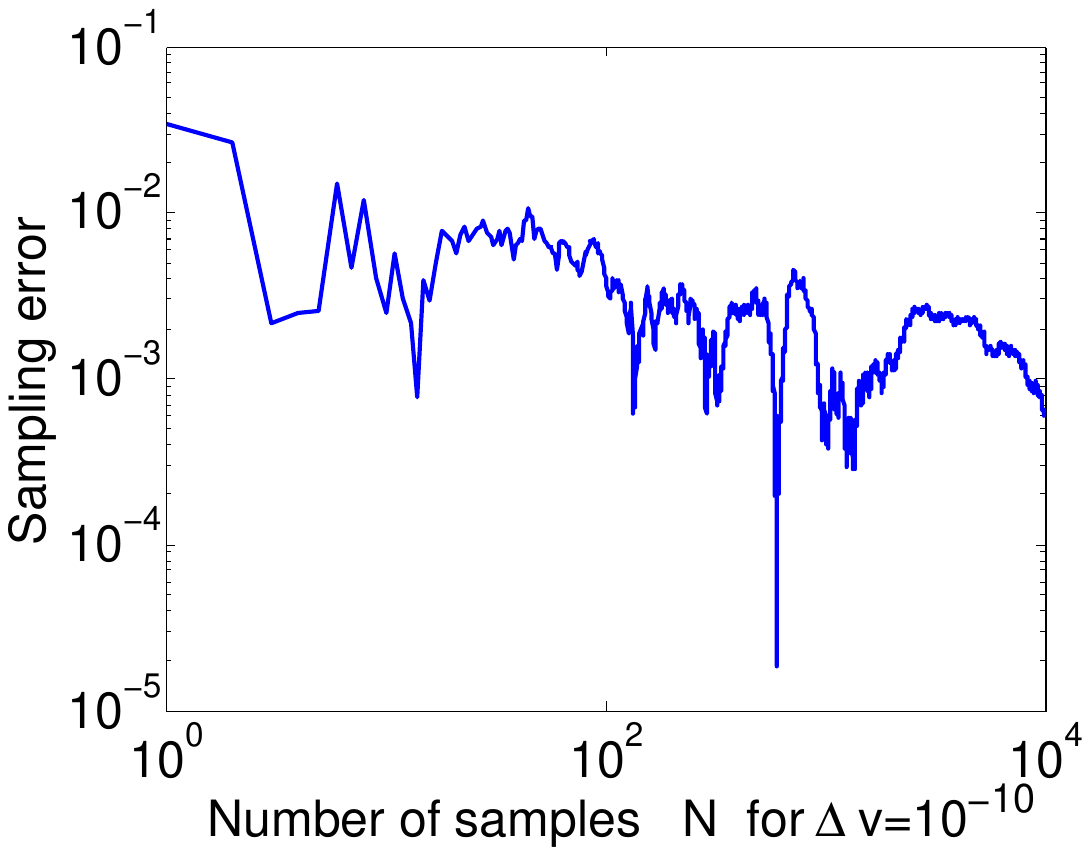}
  \includegraphics[width=0.32\textwidth,height=0.5\textwidth]{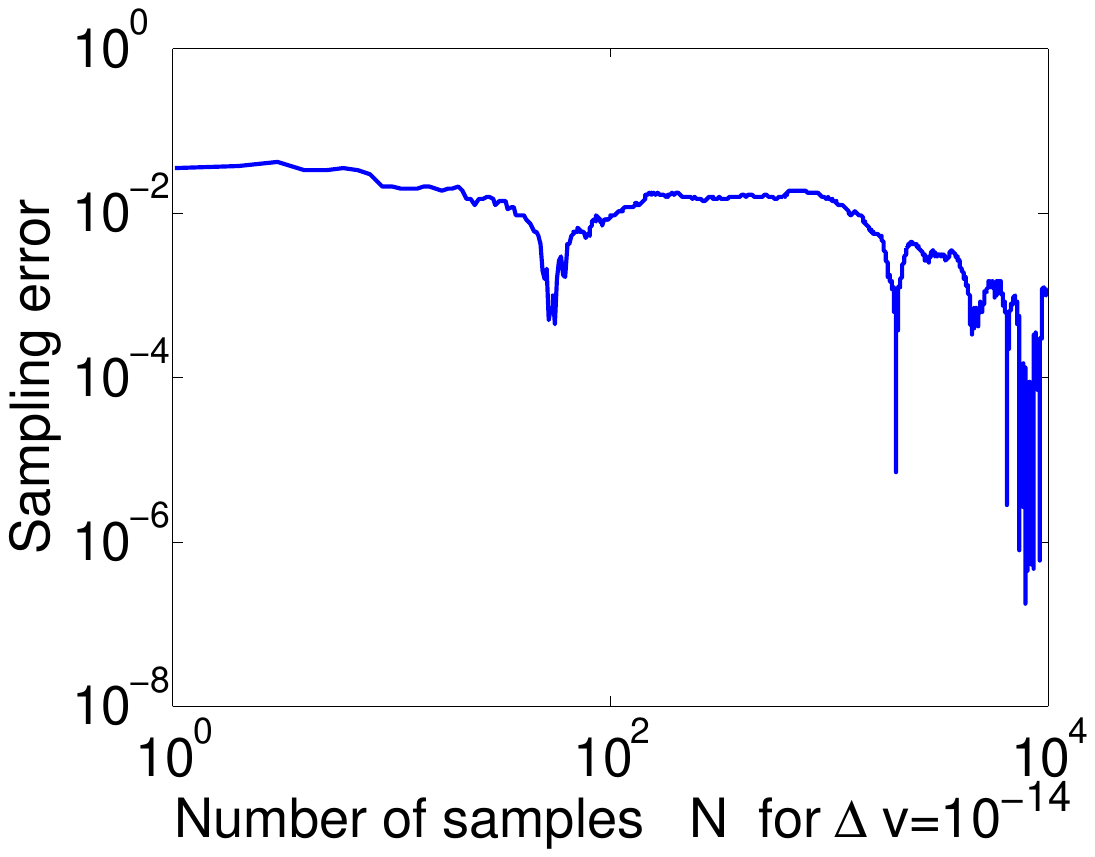}
   \caption{
  Ergodicity test results for Born-Oppenheimer molecular 
  dynamics simulations}
\label{fig45}
\end{figure}

The double average \eqref{g_bar} in time and phase-space is in Monte Carlo methods approximated by sampling several paths with random initial points $(X_0,P_0)$. For a single path one expects $\mathcal O(T^{-1/2})$ convergence rate,
in the case of bounded correlation times. Many paths correspond in some sense to longer averaging time for a single path and one may ask how long. This answer is related to how the dynamics maps an initial phase-space neighborhood $D$ to $S_{T/2}(D)$. The exponential $e^{\hat CT}$ growth of the first variation $|\partial_{P_0}^\alpha g\circ S_T(X_0,P_0)|$, for $|\alpha|=1$,
implies that an initial domain $D$, on the constant energy hyper-surface $F:=\{(X,P)\, |\, H_0(X,P)=\mbox{constant}\}$, is stretched but mapped to the bounded energy surface $F$. If the dynamics is ergodic the sampling density becomes 
uniform on this energy hyper-surface, asymptotically in time. 
Paths that are very close are not acting as independent samples for finite $T$. The assumption \eqref{g_bar} means that the exponential growth of the 
first variation measures how small the initial distance, $e^{-\hat CT/2}$, %
of the paths can be in order to lead to approximately independent paths after time $T/2$. 
The paths that are initially $e^{-\hat CT/2}$
close will after time 
$T/2$ have a distance of order one. 
The remaining $T/2$ time this $\BIGO(e^{\hat CT/2})$ 
number of paths would act as approximately independent samples (over time $T/2$),  as tested in Figure \ref{fig45},
and would give the expected decay $e^{-\hat CT/4}$,
if the Lyapunov exponent is positive in only one direction, $d=1$.
If $d$ Lyapunov exponents are similar to the maximal, we expect the decay $e^{-\hat CTd/4}$. 

 Figure \ref{fig45} tests ergodicity for Born-Oppenheimer molecular 
  dynamics simulations with the potential $\lambda_s(X)=X_1^2/2 + X_2^2/\sqrt{2}+ 2\sin(X_1X_2)$ and the double average over time and $N$ initial points
 $ \{(X_0[n],P_0[n])=(0,0,\sqrt{2(E-\lambda_s(0))} \cos(1.2+n\Delta v), \sqrt{2(E-\lambda_s(0))}\sin(1.2+n\Delta v)\}_{n=1}^N$ 
 in \eqref{g_bar}. The vertical axis shows the  sampling error
  $\sum_{n=1}^N \frac{2}{T}\int_{T/2}^T \big(g(X_t) -g_{\MD}\big)(X_0[n],P_0[n])dt/N$
  and the horizontal axis shows the number of samples $N$ for the symplectic Euler method.
   The convergence rate $\mathcal O(N^{-1/2})$ for $\Delta v=10^{-6}$ (left figure)
  and $\Delta v=10^{-10}$ (middle figure) %
  indicate independent sample paths while the closer paths with $\Delta v=10^{-14}$ in (right figure) %
  are not acting independently since the convergence is slower. The maximal Lyapunov exponent is in this case computed to be roughly $0.35$ 
  (the others are $(-0.35,0,0)$), so that the expansion of the path distance from time zero to time $70$ becomes $e^{0.35\times 70}\approx 10^{10}$ and the 
  approximate $N^{-1/2}$ convergence rate for $\Delta v\ge 10^{-10}$ indicates that $\gamma\approx 1/4$ as suggested in the paragraph following \eqref{g_bar}.
  Here $T=140$, $\Delta t=0.01$,  $E=1.5$, $g(X)=\sin(X_1X_2)$.
The ergodic limit $g_{\MD}$ is approximated in two different ways: by a single symplectic Euler path to $-0.4389$ and by the Monte Carlo samples \eqref{eq:g-MD} to $-0.4388$. 
The dashed line shows $N^{-1/2}$.

\section*{Acknowledgment}
The research was partially supported by the U.S.
National Science Foundation under the grant
NSF-DMS-0813893, U.S. DOD-ARO Grant Award W911NF-14-1-0247, 
Swedish Research Council grants 621-2010-5647 and 621-2014-4776
and the Swedish e-Science Research Center. P.P. also thanks KTH 
and Nordita for their hospitality
during his visit when the presented research was initiated.

\providecommand{\bysame}{\leavevmode\hbox to3em{\hrulefill}\thinspace}
\providecommand{\MR}{\relax\ifhmode\unskip\space\fi MR }
\providecommand{\MRhref}[2]{%
  \href{http://www.ams.org/mathscinet-getitem?mr=#1}{#2}
}
\providecommand{\href}[2]{#2}

\end{document}

%% file: ppmacros_2.tex
\def\LPROD#1#2{\langle {#1},{#2}\rangle}
     
\def\EPROD#1#2{{#1}\cdot {#2}}

\def\CHECKP{{\check P}}                       
\def\CHECKX{{\check X}}                       
                      
\def\HATP{{\hat P}}                          
\def\HATX{{\hat X}}

\def\MD{{\mathrm{MD}}}
\def\BO{0}
\def\SCH{{\mathrm{S}}}
\def\OP{{\mathrm{Op}}}

\def\HOPER{\hat{H}}
\def\VOPER{{V}}

\def\PHASE{\Theta(\CHECKX,\HATX,\CHECKP)}
\def\PHASED{\Theta'(\CHECKX,\HATX,\CHECKP)}
\def\DTHETA{\delta\LFT{\theta}(\CHECKP)}

\def\SFUN{F}

\def\PSIA{\phi} 
\def\EXPWKB{\,\EXP{\Iunit M^{1/2}\theta(X)}}

\def\AIRY{\mathrm{A}}          
\def\HATAIRY{\hat{\mathrm{A}}} 
\def\LFT#1{{#1}^*}             
\def\DPLFT#1{{\nabla_{\CHECKP}\LFT{#1}}}

\def\NEIGH{\mathcal{U}}        

\def\Oo{\mathcal{O}}

\def\C{\mathbb{C}}

\def\R{\mathbb{R}}

\def\EXP#1{e^{#1}}
\def\DET{\det}
\def\TRACE{\mathrm{Tr}\,}

\def\GRAD{\nabla}
\def\GRADX{\nabla} 

\def\LAP{\Delta}

\def\DT#1{{\dot{#1}}}
\def\DDT#1{{\ddot{#1}}}

\def\SGN{\mathrm{sgn}}

\def\COMMA{\,,}             
\def\PERIOD{\,.}            

\def\VIZ#1{(\ref{#1})}      

\def\Iunit{i}
\def\REAL{\mathrm{Re}\,}

\def\BIGO{\Oo}

